\documentclass{amsart}
\pdfoutput=1
\usepackage[margin=1.5in]{geometry}
\usepackage{amsmath}
\usepackage{amsfonts}
\usepackage{amsthm}
\usepackage{mathrsfs}
\usepackage{enumerate}
\usepackage[shortlabels]{enumitem}
\usepackage{dsfont}
\usepackage{amssymb}
\usepackage{pifont}
\usepackage{mathtools}
\usepackage{import}
\usepackage{graphicx}
\usepackage{color}
\usepackage[dvipsnames, hyperref]{xcolor}
\usepackage{thmtools}
\usepackage{hyperref, cleveref}
\usepackage{import}
\usepackage{subfiles}
\usepackage{microtype}
\usepackage{svg}
\usepackage{calc}
\usepackage[final]{showkeys}

\hypersetup{
    colorlinks=true,
    linkcolor=Blue,
    citecolor=Blue
  }

\newcommand{\parbdry}{\ensuremath{\dee_{\mathrm{par}}}}
\newcommand{\conbdry}{\ensuremath{\dee_{\mathrm{con}}}}
\newcommand{\pargps}{\mathcal{P}}
\newcommand{\pargroup}{\ensuremath{P^+}}
\newcommand{\oppgroup}{\ensuremath{P^-}}
\newcommand{\stdpargroup}{\ensuremath{P^+_\theta}}
\newcommand{\stdoppgroup}{\ensuremath{P^-_\theta}}
\newcommand{\pflags}{\ensuremath{G/P^+}}
\newcommand{\oppflags}{\ensuremath{G/P^-}}
\newcommand{\flags}{\ensuremath{G/P}}
\newcommand{\cay}{\ensuremath{\mathrm{Cay}}}
\newcommand{\conedgraph}{\ensuremath{\cay(\Gamma, S, \mc{P})}}
\newcommand{\sigmamod}{\ensuremath{\sigma_{\mathrm{mod}}}}
\newcommand{\taumod}{\ensuremath{{\tau_{\mathrm{mod}}}}}

\newcommand{\ctr}{\ensuremath{\mathrm{ctr}}}
\newcommand{\bgamh}{\ensuremath{\partial(\Gamma, \mathcal{H})}}

\newcounter{dummy}
\makeatletter
\newcommand\custitem[1][]{\item[#1]\refstepcounter{dummy}\def\@currentlabel{#1}}
\makeatother

\declaretheorem[numberwithin=section]{theorem}
\declaretheorem[sibling=theorem, name=Lemma]{lem}
\declaretheorem[sibling=theorem, name=Proposition]{prop}

\declaretheorem[sibling=theorem, name=Corollary]{cor}

\declaretheorem[sibling=theorem, style=definition, name=Assumption]{assumption}

\declaretheorem[sibling=theorem, style=definition,
name=Definition]{definition}
\declaretheorem[sibling=theorem, style=definition]{example}
\declaretheorem[sibling=theorem, style=definition]{remark}

\declaretheorem[sibling=theorem, style=definition]{notation}

\newenvironment{aside}[1]{\par\smallskip\noindent \textit{#1}.}{\par}

\newcommand{\eps}{\varepsilon}
\newcommand{\dee}{\ensuremath{\partial}}

\newcommand{\mf}[1]{\ensuremath{\mathfrak{#1}}}
\newcommand{\mc}[1]{\ensuremath{\mathcal{#1}}}

\newcommand{\mr}[1]{\ensuremath{\mathrm{#1}}}

\newcommand{\N}{\ensuremath{\mathbb{N}}}
\newcommand{\Z}{\ensuremath{\mathbb{Z}}}

\newcommand{\R}{\ensuremath{\mathbb{R}}}
\newcommand{\C}{\ensuremath{\mathbb{C}}}
\newcommand{\K}{\ensuremath{\mathbb{K}}}

\renewcommand{\H}{\ensuremath{\mathbb{H}}}
\renewcommand{\P}{\ensuremath{\mathbb{P}}}

\newcommand{\RP}{\ensuremath{\R\mathrm{P}}}

\DeclareMathOperator{\Hom}{Hom}
\DeclareMathOperator{\Homeo}{Homeo}

\DeclareMathOperator{\Opp}{Opp}

\DeclareMathOperator{\Sym}{Sym}

\DeclareMathOperator{\Gr}{Gr}

\DeclareMathOperator{\Stab}{Stab}
\DeclareMathOperator{\PSL}{PSL}
\DeclareMathOperator{\SL}{SL}

\DeclareMathOperator{\PGL}{PGL}

\DeclareMathOperator{\PO}{PO}

\title[An extended definition of relative Anosov representation]{An
  extended definition of Anosov representation for relatively
  hyperbolic groups}

\author{Theodore Weisman}

\address{Department of Mathematics, University of Michigan, Ann Arbor
  MI 48109, USA}

\email{tjwei@umich.edu}

\date{\today}

\begin{document}

\begin{abstract}
  We define a new family of discrete representations of relatively
  hyperbolic groups which unifies many existing definitions and
  examples of geometrically finite behavior in higher rank.  The
  definition includes the relative Anosov representations defined by
  Kapovich-Leeb and Zhu, and Zhu-Zimmer, as well as holonomy
  representations of various different types of ``geometrically
  finite'' convex projective manifolds. We prove that these
  representations are all stable under deformations whose restriction
  to the peripheral subgroups satisfies a dynamical condition, in
  particular allowing for deformations which do not preserve the
  conjugacy class of the peripheral subgroups.
\end{abstract}

\maketitle

\tableofcontents

\section{Introduction}

\subsection{Overview}

Classically, the best understood discrete subgroups of rank one Lie
groups have been those that are \emph{convex cocompact}, or
equivalently quasi-isometrically embedded. Only slightly less
well-behaved, however, are the geometrically finite subgroups. While
geometrically finite subgroups are not always quasi-isometrically
embedded, all of their distortion is confined to isolated ``cuspidal''
regions of their orbits. Thus, these subgroups can still be understood
using hyperbolic geometry, by piecing together the behavior of their
peripheral subgroups.

Since the practice of isolating non-hyperbolic behavior has been
succesful in rank one, it is reasonable to try and extend the
technique when investigating the still mysterious world of discrete
subgroups of higher-rank Lie groups. The natural generalization of
convex cocompact groups in higher rank is given by \emph{Anosov
  subgroups}. Originally defined by Labourie \cite{labourie2006anosov}
and Guichard-Wienhard \cite{gw2012anosov}, Anosov subgroups are a
class of Gromov-hyperbolic discrete subgroups of semisimple Lie
groups, which have many dynamical and geometric properties in common
with convex cocompact groups. They have allowed many tools and ideas
from rank-one geometry to be applied fruitfully in higher rank.

The main goal in this paper is therefore to provide a higher-rank
version of geometrical finiteness by defining a new class of subgroups
which have isolated ``non-Anosov'' behavior. As in rank one, in this
situation we can hope to understand the entire subgroup by stitching
together what happens on the isolated pieces.

Formally, we introduce the theory of \emph{extended geometrically
  finite} (or \emph{EGF}) representations as a unifying framework for
studying ``relativized Anosov'' groups. An EGF representation is
always a discrete finite-kernel representation of a relatively
hyperbolic group $\Gamma$ into some semisimple Lie group $G$. An
important feature of the definition is that it places essentially no
restrictions on the induced representations of the peripheral
subgroups of $\Gamma$. This strongly contrasts with previous (less
general) notions of relative Anosov representation due to
Kapovich-Leeb \cite{kl2018relativizing} and Zhu
\cite{zhu2019relatively} (see also \cite{zz1}). These definitions
unavoidably impose strict requirements on peripherals, and thus cannot
capture many interesting examples of intrinsically higher-rank
behavior.

The starting point for our framework is a definition of Anosov
representation in terms of topological dynamics: if $\Gamma$ is a
hyperbolic group, $G$ is a semisimple Lie group, and $P \subset G$ is
a symmetric parabolic subgroup, a representation $\rho:\Gamma \to G$
is $P$-Anosov if there is a $\rho$-equivariant embedding
$\xi:\dee \Gamma \to G/P$ of the Gromov boundary $\dee \Gamma$ of
$\Gamma$ satisfying certain dynamical properties.

Existing notions of relative Anosov representation simply replace
$\xi$ with an embedding of the Bowditch boundary of a relatively
hyperbolic group. Our idea is to instead \emph{reverse the direction}
of the boundary map: we characterize geometrical finiteness via the
existence of an equivariant map \emph{from} a closed subset of a flag
manifold \emph{to} the Bowditch boundary of a relatively hyperbolic
group, rather than the other way around.

This ``backwards'' boundary map does not need to be a homeomorphism,
which makes our definition much more flexible. For instance, the
relative Anosov subgroups of Kapovich-Leeb and Zhu can have no
subgroup which is isomorphic to a lattice in a higher-rank Lie
group---a limitation which is \emph{not} shared by EGF
representations. The definition also makes EGF representations better
suited for studying the rich world of \emph{convex projective
  structures} on manifolds, which provide a number of examples of
discrete subgroups of Lie groups displaying an intriguing mix of
``higher rank'' and ``rank one'' phenomena. In particular, EGF
representations interact well with the theory of \emph{convex
  cocompact} projective actions developed by
Danciger-Gu\'eritaud-Kassel \cite{dgk2017convex}; see
\Cref{sec:examples} below.

\subsubsection*{Relative stability}

Another significant advantage of the flexibility inherent in our
definition is that it allows for a strong \emph{stability} property
(see \Cref{thm:cusp_stable_stability} below), which generalizes the
stability of Anosov representations originally demonstrated by
Labourie \cite{labourie2006anosov}. It is not true that an arbitrary
(sufficiently small) deformation of an EGF representation is still
EGF; indeed, it is possible to find small deformations of
geometrically finite representations in rank one which are not even
discrete. However, we prove that any EGF representation
$\rho:\Gamma \to G$ is \emph{relatively} stable: any small deformation
of $\rho$ in $\Hom(\Gamma, G)$ which satisfies a condition on the
peripheral subgroups is also EGF.

This peripheral condition---which we call \emph{peripheral
  stability}---is very general, and can hold even in the absence of a
topological conjugacy between the actions of the original peripheral
subgroups and their deformations. As a result, EGF representations
provide a unifying framework for understanding the transitions between
discrete relatively hyperbolic groups with qualitatively different
cuspidal behavior (see e.g. \cite[Sec. 5]{WeismanExamples}). In
particular, the notion of peripheral stability can even describe
transitions between ``Anosov'' and ``non-Anosov'' behavior of cusp
groups. This means that our framework can be used to understand
situations where non-Anosov subgroups occur as limits of Anosov
subgroups.

In fact, even in the rank-one case (where Anosov subgroups are
precisely same as convex cocompact groups), such transitions are still
not completely understood, and in upcoming work \cite{GW2024} we plan
to explore applications of the theory of EGF representations in this
setting specifically.

We note also that the main stability theorem for EGF representations
can be used to deduce new results about less dramatic transitions
between representations. For instance, since EGF representations
generalize the relative Anosov representations of Kapovich-Leeb and
Zhu, one can apply the theorem to obtain a relative stability result
for relative Anosov representations (see \Cref{sec:relative_anosov}
below for a comparison of this result with independent work of
Zhu-Zimmer). In another direction, see \cite[Section
6]{WeismanExamples} for an application of the theorem towards an
\emph{absolute} (i.e. non-relative) stability property for a family of
discrete non-hyperbolic groups.

\subsubsection*{Relative automata}

To prove the general stability theorem, we introduce a general tool: a
\emph{relative automaton} giving a ``relative coding'' for points in
the Bowditch boundary of an arbitrary relatively hyperbolic group. The
construction is loosely based on Sullivan's symbolic coding of points
in the limit set of a convex cocompact group in rank one (recently
adapted and generalized by Kapovich-Kim-Lee \cite{kkl2019structural}),
as well as an ``automatic'' description of Anosov representations due
to Bochi-Potrie-Sambarino \cite{bps2019anosov}. We expand on the basic
idea in two different ways simultaneously.

First, rather than coding \emph{points} in the ``limit set'' of some
discrete faithful representation, we essentially code \emph{fibers} in
the invariant set $\Lambda$ surjecting onto the Bowditch boundary of
our relatively hyperbolic group. Second, we provide a way to code
parabolic points (or more accurately, ``parabolic fibers'') in a way
which is compatible with the coding for generic points in the Bowditch
boundary.

We remark again that even though both of these innovations were
originally developed to better understand the behavior of discrete
subgroups of higher-rank Lie groups, they are also directly useful
outside of this context. For example, in \cite{MMW1}, \cite{MMW2}, we
show how to apply the ideas appearing in this paper towards the theory
of abstract (relatively) hyperbolic groups, and in future work
\cite{GW2024} we will explore applications to geometrically finite
subgroups in rank-one. It seems possible that there are additional
applications beyond even these, for instance towards subgroups of
mapping class groups.

\subsection{Main definition}

In the rest of the introduction we give some more detail regarding the
central definition of this paper and the main results surrounding it.

If $\Gamma$ is a relatively hyperbolic group, relative to a collection
$\mc{H}$ of peripheral subgroups, then $\Gamma$ acts as a
\emph{convergence group} on the Bowditch boundary $\bgamh$. We recall
the definition here.

\begin{definition}
  \label{defn:convergence_group}
  Let $\Gamma$ act on a topological space $M$. The group $\Gamma$ is
  said to act as a \emph{convergence group} if for every infinite
  sequence of distinct elements $\gamma_n \in \Gamma$, there exist
  points $a, b \in M$ and a subsequence $\gamma_m \in \Gamma$ such
  that $\gamma_m$ converges uniformly on compacts in $M - \{a\}$ to
  the constant map $b$.
\end{definition}

When $\gamma_n$ is a sequence of distinct elements in a relatively
hyperbolic group $\Gamma$, then $\gamma_n$ converges to $b$ uniformly
on compacts in $\bgamh - \{a\}$ if and only if $\gamma_n$ converges to
$b$ in the compactification $\overline{\Gamma} = \Gamma \sqcup \bgamh$
and the inverse sequence $\gamma_n^{-1}$ converges to $a$ (see
\cite[Theorem 3A]{tukia1994convergence}).

Recall that if a group $\Gamma$ acts by homeomorphisms on a Hausdorff
space $X$, the pair $(\Gamma, X)$ is called a \emph{topological
  dynamical system}. We say that an \emph{extension} of $(\Gamma, X)$
is a topological dynamical system $(\Gamma, Y)$ together with a
$\Gamma$-equivariant surjective map $\phi:Y \to X$.

In this paper, when $(\Gamma, \mc{H})$ is a relatively hyperbolic pair
(i.e. $\Gamma$ is hyperbolic relative to a collection $\mc{H}$ of
peripheral subgroups), we will consider \emph{embedded extensions} of
the topological dynamical system $(\Gamma, \bgamh)$. We want these
embedded extensions to respect the convergence group action of
$(\Gamma, \bgamh)$ in some sense, so we introduce the following
definition:
\begin{restatable}{definition}{dynamicsPreserving}
  \label{defn:extends_convergence_dynamics}
  Let $(\Gamma, \mc{H})$ be a relatively hyperbolic pair, with
  $\Gamma$ acting on a connected compact metrizable space $M$ by
  homeomorphisms. Let $\Lambda \subset M$ be a closed
  $\Gamma$-invariant set.
  
  We say that a continuous equivariant surjective map
  $\phi:\Lambda \to \bgamh$ \emph{extends the convergence action of
    $\Gamma$} if for each $z \in \bgamh$, there exists an open set
  $C_z \subset M$ satisfying the following:
  \begin{enumerate}[(C1)]
  \item\label{item:strata_interior} For every compact
    $K \subset \bgamh$, the set $\bigcap_{z \in K}C_z$ contains an
    open neighborhood of $\phi^{-1}(\bgamh - K)$.
  \item\label{item:extended_convergence} For every sequence $\gamma_n$
    in $\Gamma$ with $\gamma_n^{\pm 1} \to z_\pm$ for
    $z_\pm \in \bgamh$, then for any compact set $K \subset C_{z_-}$
    and any open set $U$ containing $\phi^{-1}(z_+)$, for sufficiently
    large $n$, $\gamma_n \cdot K$ lies in $U$.
  \end{enumerate}
  In this situation, we also say that the action of $\Gamma$ on $M$ is
  an \emph{extended convergence action} (or that $\Gamma$ acts on $M$
  as an \emph{extended convergence group}). The map $\phi$ is called a
  \emph{boundary extension} for the action, the closed invariant set
  $\Lambda$ is called a \emph{boundary set}, and the sets $C_z$ are
  called \emph{repelling strata}.
\end{restatable}

Now let $G$ denote a semisimple Lie group with no compact factor. The
central definition of the paper is the following:
\begin{restatable}{definition}{EGF}
  \label{defn:boundary_dynamics}
  Let $\Gamma$ be a relatively hyperbolic group, let
  $\rho:\Gamma \to G$ be a representation, and let $P \subset G$ be a
  symmetric parabolic subgroup. We say that $\rho$ is \emph{extended
    geometrically finite} (EGF) with respect to $P$ if there exists a
  closed $\rho(\Gamma)$-invariant set $\Lambda \subset G/P$ and a
  continuous $\rho$-equivariant surjective antipodal map
  $\phi:\Lambda \to \bgamh$ extending the convergence action of
  $\Gamma$.
\end{restatable}

We refer to \Cref{sec:egf_representations} for the definition of
``antipodal map'' in this context.

\subsection{Main results}

Like (relative) Anosov representations, extended geometrically finite
representations are always discrete with finite kernel (see
\ref{sec:discrete_finite_kernel}). The central result of this paper
says that EGF representations have a \emph{relative stability}
property: if $\rho$ is an EGF representation, then certain small
relative deformations of $\rho$ must also be EGF.

To state the theorem, we define a notion of a \emph{peripherally
  stable} subspace of $\Hom(\Gamma, G)$. The precise definition is
given in \Cref{sec:relative_stability}, but roughly speaking, a
subspace $\mathcal{W} \subseteq \Hom(\Gamma, G)$ is \emph{peripherally
  stable} if the large-scale dynamical behavior of the peripheral
subgroups of $\Gamma$ is in some sense preserved by small deformations
inside of $\mathcal{W}$. We emphasize again that the action of a
deformed peripheral subgroup does \emph{not} need to be even
topologically conjugate to the action of the original peripheral
subgroup.

We prove the following:
\begin{restatable}{theorem}{cuspStableStability}
  \label{thm:cusp_stable_stability}
  Let $\rho:\Gamma \to G$ be EGF with respect to $P$, let
  $\phi:\Lambda \to \bgamh$ be a boundary extension, and let
  $\mc{W} \subseteq \Hom(\Gamma, G)$ be peripherally stable at
  $(\rho, \phi)$. For any compact subset $Z$ of $\bgamh$ and any open
  set $V \subset \flags$ containing $\phi^{-1}(Z)$, there is an open
  subset $\mc{W}' \subset \mc{W}$ containing $\rho$ such that each
  $\rho' \in \mc{W}'$ is EGF with respect to $P$, and has an EGF
  boundary extension $\phi'$ satisfying $\phi'^{-1}(Z) \subset V$.
\end{restatable}

\begin{remark}
  When $\rho:\Gamma \to G$ is a $P$-Anosov representation of a
  hyperbolic group $\Gamma$, then the associated boundary embedding
  $\dee \Gamma \to \flags$ also varies continuously with $\rho$ in the
  compact-open topology on maps $\dee\Gamma \to \flags$. Since EGF
  representations come with boundary \emph{extensions} (rather than
  embeddings), \Cref{thm:cusp_stable_stability} only gives us a
  \emph{semicontinuity} result.

  We expect that it is possible to extend the methods of this paper to
  prove stronger continuity results when the original EGF
  representation $\rho:\Gamma \to G$ satisfies additional assumptions;
  see the discussion following \Cref{cor:rel_anosov_cusp_stable}.
\end{remark}

While the peripheral stability condition in
\Cref{thm:cusp_stable_stability} is mildly technical, we can also
apply it to yield more concrete results.
\begin{definition}
  Let $(\Gamma, \mc{H})$ be a relatively hyperbolic pair, and let
  $\rho:\Gamma \to G$ be a representation. The space of
  \emph{cusp-preserving} representations
  $\Hom_{\mr{cp}}(\Gamma, G, \mc{H}, \rho)$ is the set of
  representations $\rho':\Gamma \to G$ such that for each peripheral
  subgroup $H \in \mc{H}$, we have
  $\rho'|_H = g \cdot \rho|_H \cdot g^{-1}$ for some $g \in G$ (which
  may depend on $H$).
\end{definition}

\begin{cor}
  \label{cor:cusp_preserving_stability}
  Let $\rho:\Gamma \to G$ be an EGF representation. Then there is a
  neighborhood of $\rho$ in $\Hom_{\mr{cp}}(\Gamma, G, \mc{H}, \rho)$
  consisting of EGF representations.
\end{cor}

\Cref{cor:cusp_preserving_stability} gives a very restrictive example
of a peripherally stable subspace of $\Hom(\Gamma, G)$. But, in
general the peripheral stability condition is flexible enough to allow
peripheral subgroups to deform in nontrivial ways.

In particular, it is possible to find peripherally stable deformations
of an EGF representation $\rho:\Gamma \to \PGL(d, \R)$ which change
the Jordan block decomposition of elements in the peripheral
subgroups. For instance, one can deform an EGF representation in
$\PGL(d, \R)$ with unipotent peripheral subgroups into an EGF
representation with diagonalizable peripheral subgroups---see
\Cref{ex:peripheral_stability}.

\subsection{Examples}
\label{sec:examples}

The related paper \cite{WeismanExamples} is focused on examples of EGF
representations. For illustrative purposes, however, we briefly
describe some examples here.

\subsubsection{Convex projective structures}
\label{sec:convex_projective}
A host of examples of Anosov representations arise from the theory of
\emph{convex projective structures}; see
e.g. \cite{benoist04convexes}, \cite{benoist2006convexeshyp},
\cite{kapovich2007}, \cite{DGKgeomded}, \cite{dgklm2021convex}. In
fact, work of Danciger-Gu\'eritaud-Kassel \cite{dgk2017convex} and
Zimmer \cite{zimmer2021projective} implies that Anosov representations
can be essentially characterized as holonomy representations of
\emph{convex cocompact} projective orbifolds with hyperbolic
fundamental group. However, convex projective structures also yield a
number of interesting examples of discrete non-Anosov subgroups of
$\PGL(d, \R)$. In many cases, the groups in question are relatively
hyperbolic, and appear to have ``geometrically finite'' properties.

The theory of EGF representations is well-suited to these
examples. For instance, in \cite{WeismanExamples}, we apply some of
our previous work \cite{weisman2023dynamical} together with work of
Islam-Zimmer \cite{iz2022structure} to see that whenever a subgroup
$\Gamma \subset \PGL(d, \R)$ is relatively hyperbolic and
\emph{projectively convex cocompact} in the sense of
\cite{dgk2017convex}, then the inclusion
$\Gamma \hookrightarrow \PGL(d, \R)$ is EGF with respect to the
parabolic subgroup stabilizing a flag of type $(1, d-1)$ in $\R^d$. If
$\Gamma$ is not hyperbolic, then these examples are \emph{not} covered
by other definitions of relative Anosov representations (see
\cite[Remark 1.14]{weisman2023dynamical}); such non-hyperbolic
examples have been constructed in e.g. \cite{benoist2006convexes},
\cite{bdl2015convex}, \cite{clm2020convex}, \cite{CLM2},
\cite{dgklm2021convex}, \cite{BV}.

In \cite{cm2014finitude}, Crampon-Marquis introduced several
definitions of ``geometrical finiteness'' for \emph{strictly} convex
projective manifolds. Not all of their definitions are equivalent (and
in fact some of the originally claimed equivalences were incorrect);
recently Blayac--Marquis \cite{BCM12} gave a correction clarifying the
situation greatly. In \cite{zhu2019relatively}, Zhu proved that
manifolds satisfying one of their definitions have ``relative Anosov''
holonomy in the sense of \cite{kl2018relativizing},
\cite{zhu2019relatively}, \cite{zz1} (see \Cref{sec:relative_anosov}),
which means they also have EGF holonomy by
\Cref{thm:bdry_dynamics_homeo_equiv_rel_asymp_embedded}
below. Examples can be found by deforming geometrically finite groups
in $\PO(d, 1)$ into $\PGL(d+1, \R)$ while keeping the conjugacy
classes of cusp groups fixed (see \cite{ballas2014},
\cite{bm2020properly}, \cite{clt2018deforming}), or via Coxeter
reflection groups \cite{CLM2}.

There are, however, more general notions of ``geometrically finite''
convex projective structures. In \cite{clt2018deforming},
Cooper-Long-Tillmann considered the situation of a convex projective
manifold $M$ (with strictly convex boundary) which is a union of a
compact piece and finitely many ends homeomorphic to
$N \times [0, \infty)$, where $N$ is a compact manifold with virtually
nilpotent fundamental group. The ends of such a manifold are called
``generalized cusps,'' and the possible ``types'' of generalized cusps
were later classified by Ballas-Cooper-Leitner
\cite{bcl2020generalized}. Examples of projective manifolds with
generalized cusps have been produced by Ballas \cite{ballas2021},
Ballas-Marquis \cite{bm2020properly}, and Bobb
\cite{bobb2019convex}. In general the holonomy representations of
these manifolds are \emph{not} ``relative Anosov'' in the sense of
\cite{zz1}, but in \cite{WeismanExamples} we prove that they do
provide additional examples of EGF representations. The proof is an
application of \Cref{thm:cusp_stable_stability}: it turns out that
peripheral stability is actually flexible enough to allow for
deformation between the different Ballas-Cooper-Leitner generalized
cusp types.

\begin{remark}
  In the EGF examples described above, $M$ is a convex projective
  manifold which decomposes into a compact piece and finitely many
  generalized cusps, and the universal cover of the compact part is a
  subset of $\tilde{M} \subset \P(\R^d)$ with no nontrivial segments
  in its ideal boundary. There are indications that, if $M$ is
  \emph{any} convex projective manifold satisfying these conditions,
  then $M$ should have EGF holonomy; see for example \cite{choi2011},
  \cite{wolf2020} and the general setup in \cite{iz2022structure},
  \cite{BV}. However, we do not (yet) have a general theorem for a
  statement along these lines.
\end{remark}

After a version of this paper originally appeared as a preprint,
Blayac-Viaggi \cite{BV} also produced still more general examples of
convex projective $n$-manifolds which decompose into a compact piece
and several projective ``cusps.'' In these examples (which can arise
as limits of convex cocompact representations), each cusp is finitely
covered by a product $N \times S^1 \times [0, \infty)$, where $N$ is a
closed hyperbolic manifold of dimension $n - 2$. Consequently, these
manifolds do not have ``generalized cusps'' in the sense of
Cooper-Long-Tillmann, and their fundamental groups cannot even admit
relative Anosov representations. Nevertheless, Blayac-Viaggi showed
that the holonomy representations of their examples are always EGF.

\subsubsection{Other examples}

In \cite{WeismanExamples}, we construct additional examples of EGF
representations by considering \emph{compositions} of projectively
convex cocompact representations $\rho:\Gamma \to \PGL(V)$ with the
symmetric representation $\tau_k:\PGL(V) \to \PGL(\Sym^kV)$. We show
that, assuming the peripheral subgroups in $\Gamma$ are all virtually
abelian, then the composition $\tau_k \circ \rho$ is still EGF; this
holds even though the compositions are \emph{not} believed to be
convex cocompact.

We are also able to prove that the entire space
$\Hom(\Gamma, \PGL(\Sym^kV))$ is peripherally stable about
$\tau_k \circ \rho$. Via \Cref{thm:cusp_stable_stability}, this gives
a new source of examples of stable discrete subgroups of higher-rank
Lie groups.

\subsection{Comparison with relative Anosov representations}
\label{sec:relative_anosov}

Previously, Kapovich--Leeb \cite{kl2018relativizing} and Zhu
\cite{zhu2019relatively} independently introduced several notions of a
\emph{relative $P$-Anosov} representation, where $P$ is a parabolic
subgroup of a semisimple Lie group $G$. Later work of Zhu--Zimmer
\cite{zz1} showed that Zhu's definition (that of a \emph{relatively
  dominated} representation) is equivalent to essentially the most
general Kapovich--Leeb definition (specifically, of a \emph{relatively
  asymptotically embedded} representation). Zhu--Zimmer further proved
that several of the Kapovich--Leeb definitions are equivalent in the
case where $G = \SL(d, \mathbb{K})$ ($\mathbb{K} = \R$ or $\C$). In
the special case where the domain group is isomorphic to a Fuchsian
group, these definitions also all agree with a notion of relative
Anosov representation for Fuchsian groups introduced by
Canary--Zhang--Zimmer \cite{czz2021cusped}.

For the rest of this paper we refer to representations satisfying any
of these equivalent definitions as \emph{relative Anosov}
representations. Note that we will never simply say ``Anosov
representation'' when we mean to refer to a \emph{relative} Anosov
representation (this differs from the convention in
\cite{czz2021cusped}).

\begin{remark}
  For a parabolic subgroup $P$ in a general semisimple Lie group $G$,
  the different definitions for a ``relatively $P$-Anosov
  representation'' suggested by Kapovich--Leeb are not all the
  same. See \Cref{ex:nonuniform_example}.
\end{remark}

The theorem below says that extended geometrically finite
representations are a strict generalization of relative Anosov
representations. Additionally, it provides a precise characterization
of when an EGF representation is also relative Anosov.
\begin{restatable}{theorem}{relAsympBdryDynamics}
  \label{thm:bdry_dynamics_homeo_equiv_rel_asymp_embedded}
  Let $(\Gamma, \mc{H})$ be a relatively hyperbolic pair, and let
  $P \subset G$ be a symmetric parabolic subgroup. A representation
  $\rho:\Gamma \to G$ is relatively $P$-Anosov if and only if $\rho$
  is EGF with respect to $P$, and has an injective boundary extension
  $\phi:\Lambda \to \bgamh$.
\end{restatable}

We emphasize again that almost all of the EGF examples mentioned in
\Cref{sec:examples} are \emph{not} relative Anosov, so the theorem
tells us that in many cases it is actually \emph{not} possible to
construct an injective boundary extension for a given EGF
representation.

\begin{remark}
  By \Cref{prop:conical_pts_singletons}, any EGF representation has a
  boundary extension which is injective on preimages of conical limit
  points. So, in the case where the peripheral structure $\mc{H}$ is
  trivial (meaning that $\Gamma$ is a hyperbolic group and $\bgamh$ is
  identified with the Gromov boundary $\dee \Gamma$ of $\Gamma$),
  \Cref{thm:bdry_dynamics_homeo_equiv_rel_asymp_embedded} implies that
  EGF representations are precisely the same as (non-relative) Anosov
  representations.

  This actually gives a new characterization of Anosov
  representations, since a priori the EGF boundary extension $\phi$
  surjecting onto the Gromov boundary of a hyperbolic group does not
  need to be a homeomorphism; the theorem tells us that if such a
  boundary extension exists, then it is possible to replace $\phi$
  with an injective boundary extension, whose inverse is the Anosov
  boundary map.
\end{remark}

\subsubsection{Stability for relative Anosov representations}

In \cite{kl2018relativizing}, Kapovich--Leeb suggested that a relative
stability result should hold for relative Anosov representations, but
did not give a precise statement. By applying
\Cref{thm:cusp_stable_stability}, \Cref{prop:conical_pts_singletons},
and \Cref{thm:bdry_dynamics_homeo_equiv_rel_asymp_embedded}, we obtain
the following stability theorem:
\begin{theorem}
  \label{thm:relative_anosov_peripheral_stable}
  Let $(\Gamma, \mc{H})$ be a relatively hyperbolic pair, let
  $\rho:\Gamma \to G$ be a relative $P$-Anosov representation, and let
  $\mc{W} \subset \Hom(\Gamma, G)$ be a peripherally stable subspace,
  such that for each $H \in \mc{H}$ and each $\rho' \in \mc{W}$, the
  restriction $\rho'|_H$ is $P$-divergent with $P$-limit set a
  singleton. Then an open neighborhood of $\rho$ in $\mc{W}$ consists
  of relative $P$-Anosov representations of $\Gamma$.
\end{theorem}

If we restrict the allowable peripheral deformations to conjugacies,
this result reduces to:
\begin{cor}
  \label{cor:rel_anosov_cusp_stable}
  Let $(\Gamma, \mc{H})$ be a relatively hyperbolic pair, and let
  $\rho:\Gamma \to G$ be a relative $P$-Anosov representation. There
  is an open neighborhood of $\rho$ in
  $\Hom_{\mr{cp}}(\Gamma, G, \mc{H}, \rho)$ consisting of relative
  $P$-Anosov representations.
\end{cor}

\begin{remark}
  In the special case where $\Gamma$ is isomorphic to a Fuchsian
  group, \Cref{cor:rel_anosov_cusp_stable} follows from previous work
  of Canary-Zhang-Zimmer
  \cite{czz2021cusped}. \Cref{cor:rel_anosov_cusp_stable} itself was
  also proved independently by Zhu-Zimmer \cite{zz1}, who showed
  further that the associated relative boundary maps vary continuously
  (in fact, analytically). The methods used in \cite{czz2021cusped}
  and \cite{zz1} are considerably different from those in this paper,
  and they do \emph{not} imply \Cref{thm:cusp_stable_stability} or
  \Cref{thm:relative_anosov_peripheral_stable} (and they do not apply
  to most of the examples of deformations considered in the companion
  paper \cite{WeismanExamples}).

  As Zhu-Zimmer observe, it seems unlikely that their techniques can
  be easily adapted for the general study of EGF representations, or
  for the general class of relative deformations considered in this
  article. On the other hand, while we expect that the methods in this
  paper could be used to generalize the Zhu-Zimmer result regarding
  \emph{continuously} varying boundary embeddings for relative Anosov
  representations, it does not seem simple to use our approach to
  study analytic variation of the boundary maps.
\end{remark}

\begin{remark}
  After this paper was originally posted, Wang \cite{wang2023b} showed
  that the situation of an EGF representation $\rho$ with
  $P$-divergent image can be interpreted in terms of \emph{restricted
    Anosov} representations, i.e. representations which are Anosov
  ``along a subflow'' of a certain flow space associated to the
  representation (see \cite{wang2023a}). Using these ideas, Wang
  proves a version of \Cref{cor:cusp_preserving_stability} for this
  special class of EGF representations.
\end{remark}

\subsection{Further applications, and potential future applications}

\subsubsection{Anosov relativization}

When $\Gamma$ is a relatively hyperbolic group, the Bowditch boundary
of $\Gamma$ (and thus, the definition of an EGF representation)
depends on the choice of peripheral structure $\mc{H}$ for
$\Gamma$. In general, there might be more than one possible choice:
for instance, if a group $\Gamma$ is hyperbolic relative to a
collection $\mc{H}$ of hyperbolic subgroups, then $\Gamma$ is itself
hyperbolic, relative to an \emph{empty} collection of peripheral
subgroups (see \cite[Corollary 1.14]{ds2005tree}).

In this paper, we prove the following \emph{Anosov relativization
  theorem}:
\begin{theorem}
  \label{thm:changing_peripheral_structure}
  Let $(\Gamma, \mc{H})$ be a relatively hyperbolic pair, and suppose
  that each $H \in \mc{H}$ is hyperbolic. If $\rho:\Gamma \to G$ is an
  EGF representation with respect to $P$ for the peripheral structure
  $\mc{H}$, and $\rho$ restricts to a $P$-Anosov representation on
  each $H \in \mc{H}$, then $\rho$ is a $P$-Anosov representation of
  $\Gamma$.
\end{theorem}

See \cite{GW2024} for one application of this theorem. Another
potential application is the construction of new examples of Anosov
representations: one could start with an EGF representation
$\rho:\Gamma \to G$ which is \emph{not} Anosov, and then attempt to
find a peripherally stable deformation of $\rho$ which restricts to an
Anosov representation on peripheral
subgroups. \Cref{thm:cusp_stable_stability} and
\Cref{thm:changing_peripheral_structure} would then imply that the
original representation $\rho$ can be realized as a non-Anosov limit
of Anosov representations in the peripherally stable deformation
space. We plan to explore applications of this type in upcoming work
with Danciger.

\subsubsection{Limits of Anosov representations}

In \cite{LLS}, Lee-Lee-Stecker considered the deformation space of
Anosov representations $\rho:\Gamma_{p,q,r} \to \SL_\pm(3, \R)$, where
$\Gamma_{p,q,r}$ is a triangle reflection group, and showed that
certain components of this space have representations in their
boundary which are \emph{not} Anosov. Interestingly, these limiting
representations still have equivariant injective boundary maps from
$\dee \Gamma_{p,q,r}$ into the space of full flags in $\R^3$, but they
fail to be Anosov because the boundary maps fail to be transverse.

The limiting representations in question are \emph{not} relatively
Anosov, but it turns out that they are EGF. More precisely, we have
the following:
\begin{prop}
  Let $\Gamma_{p,q,r}$ be a hyperbolic triangle reflection group with
  $p,q,r$ odd, so that (by \cite{LLS}) the space of Anosov
  representations $\Gamma_{p,q,r} \to \SL_\pm(3,\R)$ has a component
  whose closure in the character variety is a compact interval. The
  representations in the boundary of this interval are not relatively
  Anosov, but they are EGF, with respect to the peripheral structure
  corresponding to the cyclic subgroup generated by the Coxeter
  element $\gamma \in \Gamma_{p,q,r}$.
\end{prop}
\begin{proof}[Proof sketch]
  Let $\rho$ be one of the representations in the boundary of the
  interval; as mentioned above, \cite{LLS} shows that there is a
  (non-antipodal) $\rho$-equivariant embedding from
  $\dee \Gamma_{p,q,r}$ to the space of full flags in $\R^3$, with
  image $\Lambda$. There is a continuous equivariant surjective map
  $\dee \Gamma_{p,q,r} \to \dee(\Gamma_{p,q,r}, \mc{H})$, where
  $\mc{H}$ is the collection of conjugates of the cyclic subgroup
  $\langle \gamma \rangle$ (see \Cref{sec:relativization}). This
  implies that the boundary embedding determines a continuous
  equivariant map $\Lambda \to \dee(\Gamma_{p,q,r},
  \mc{H})$. \cite[Lem. 7.1]{LLS} implies that this map is actually
  antipodal, so this is the candidate boundary extension.

  One may then check (say using the result in the appendix to this
  paper) that the general construction in \cite[Sec. 4]{LLS} implies
  that $\rho$ is also \emph{$P_{1,2}$-divergent} along quasi-geodesic
  sequences, with $P_{1,2}$-limit points lying in the image of the
  boundary embedding; here $P_{1,2}$ is the stabilizer of a full flag
  in $\R^3$. This makes it straightforward to verify the hypotheses of
  several of our alternative characterizations of EGF representations
  (either \Cref{prop:conical_peripheral_implies_egf} or
  \Cref{prop:weakest_egf}). That $\rho$ is not relatively Anosov may
  be seen from
  \Cref{thm:bdry_dynamics_homeo_equiv_rel_asymp_embedded}, together
  with the fact that the $P_{1,2}$-limit set of
  $\rho(\langle \gamma \rangle)$ is not a singleton.
\end{proof}

Together with the Anosov relativization theorem mentioned above, the
proposition provides evidence that EGF representations can serve as a
useful tool in the study of boundaries of spaces of Anosov
representations. In addition, it gives a source of examples of EGF
representations which do not directly derive from convex projective
structures.

\subsubsection{Deformations in rank one}
\label{sec:rank_one}

Even in rank one, the deformation theory of geometrically finite
representations is not completely understood. In
\cite{bowditch1998spaces}, Bowditch described circumstances which
guarantee that a small deformation of a geometrically finite group
$\Gamma \subset \PO(d, 1)$ is still geometrically finite, but his
criteria do not have an obvious analog in other rank one Lie
groups. Moreover, the conditions Bowditch gives are too strict to
allow for deformations which change the homeomorphism type of the
limit set $\Lambda(\Gamma)$.  Such deformations exist and are often
peripherally stable, meaning that the EGF framework could be used to
understand them further. It even seems possible that a version of the
theory could be applied in circumstances where the isomorphism type of
$\Gamma$ is allowed to change.

\subsection{Outline of the paper}

We begin by providing some background in
\Cref{sec:rel_hyp_background,sec:lie_theory_background}, and then give
the full formal definition of EGF representations in
\Cref{sec:egf_definition}. In that section we also prove
\Cref{thm:bdry_dynamics_homeo_equiv_rel_asymp_embedded} (giving the
connection between EGF representations and relative Anosov
representations) and \Cref{thm:changing_peripheral_structure} (the
Anosov relativization theorem). Some of these proofs assume the
results of later sections, but they are not relied upon anywhere else
in the paper.

The rest of the paper is devoted to the proof of our main stablity
theorem for EGF representations (\Cref{thm:cusp_stable_stability}). In
\Cref{sec:relative_automata} and \Cref{sec:extended_convergence}, we
develop the main technical tool needed for the proof, which involves
using the notion of an extended convergence group action to construct
the \emph{relative quasigeodesic automaton} alluded to
previously. Then, in \Cref{sec:contracting_paths}, we prove a key
result (\Cref{prop:caratheodory_metric_shrinks}) regarding a metric on
certain open subsets of flag manifolds $\flags$, which we use to
relate the results of the previous sections to relatively hyperbolic
group actions on $\flags$. Then, we use all of these tools to develop
an alternative characterization of EGF representations in
\Cref{sec:weak_egf}, and finally prove our main theorem in
\Cref{sec:relative_stability}.

\subsection{Acknowledgements}

The author thanks his PhD advisor, Jeff Danciger, for encouragement
and many helpful conversations---without which this paper could not
have been written. The author also thanks Katie Mann and Jason Manning
for assistance simplifying some of the arguments in
\Cref{sec:relative_automata,sec:extended_convergence}. Further thanks
are owed to Daniel Allcock, Dick Canary, Fanny Kassel, Max
Riestenberg, Anna Wienhard, Feng Zhu, and Andy Zimmer for providing
feedback on various versions of this project.

This work was supported in part by NSF grants DMS-1937215 and
DMS-2202770.


\section{Relative hyperbolicity}
\label{sec:rel_hyp_background}

In this section we discuss some of the basic theory of relatively
hyperbolic groups, mostly to establish the notation and conventions we
will use throughout the paper. We refer to \cite{bh1999metric},
\cite{bowditch2012relatively}, \cite{ds2005tree} for background on
hyperbolic groups and relatively hyperbolic groups. See also section 3
of \cite{kl2018relativizing} for an overview (which we follow in part
here).

\begin{notation}
  Throughout this paper, if $X$ is a metric space, $A$ is a subset of
  $X$, and $r \ge 0$, we let $N_X(A, r)$ denote the open
  $r$-neighborhood in $X$ about $A$. For a point $x \in X$, we let
  $B_X(x, r)$ denote the open $r$-ball about $x$.

  When the metric space $X$ is implied from context, we will often
  just write $N(A, r)$ or $B(x, r)$.
\end{notation}

\subsection{Geometrically finite actions}

Recall that a finitely generated group $\Gamma$ is \emph{hyperbolic}
(or \emph{word-hyperbolic} or \emph{$\delta$-hyperbolic} or
\emph{Gromov-hyperbolic}) if and only if it acts properly
discontinuously and cocompactly on a $\delta$-hyperbolic proper
geodesic metric space $Y$.

A \emph{relatively hyperbolic group} is also a group with an action by
isometries on a $\delta$-hyperbolic proper geodesic metric space $Y$,
but instead of asking for the action to cocompact, we ask for the
action to be in some sense ``geometrically finite.''

To be precise, this means that $Y$ has a $\Gamma$-invariant
decomposition into a \emph{thick part} $Y_{\mr{th}}$ and a countable
collection $\mc{B}$ of \emph{horoballs}. For a horoball $B$, we let
$\ctr(B)$ denote the center of $B$ in $\dee Y$, and we let
$\Stab_\Gamma(p)$ denote the stabilizer of any $p \in \dee Y$.

\begin{definition}
  \label{defn:relatively_hyperbolic_group}
  Let $\Gamma$ be a finitely generated group acting on a hyperbolic
  metric space $Y$, and let $\mc{B}$ be a countable collection of
  horoballs in $Y$, invariant under the action of $\Gamma$ on $Y$. If:
  \begin{enumerate}
  \item The action of $\Gamma$ on the closure of
    $Y_{\mr{th}} = Y - \bigcup_{B \in \mc{B}} B$ is cocompact, and
  \item for each $B \in \mc{B}$, the stabilizer of $\ctr(B)$ in
    $\Gamma$ is infinite,
  \end{enumerate}
  then we say that $\Gamma$ is a \emph{relatively hyperbolic group},
  relative to the collection
  $\mc{H} = \{\Stab_\Gamma(p) : p = \ctr(B) \textrm{ for } B \in
  \mc{B}\}$.
\end{definition}

\begin{definition}
  Let $\Gamma$ be a relatively hyperbolic group, relative to a
  collection of subgroups $\mc{H}$.
  \begin{itemize}
  \item The centers of the horoballs in $\mc{B}$ are called
    \emph{parabolic points} for the $\Gamma$-action on $\dee Y$. The
    set of parabolic points in $\dee Y$ is denoted $\parbdry Y$.
  \item The parabolic point stablizers
    $\mc{H} = \{\Stab_\Gamma(p) : p \in \parbdry Y\}$ are called
    \emph{peripheral subgroups}. We often write $\Gamma_p$ for
    $\Stab_\Gamma(p)$.
  \end{itemize}
\end{definition}

A group $\Gamma$ might be hyperbolic relative to different collections
$\mc{H}$, $\mc{H}'$ of peripheral subgroups. The collection $\mc{H}$
of peripheral subgroups is sometimes called a \emph{peripheral
  structure} for $\Gamma$.

\begin{definition}
  Let $\Gamma$ be a finitely generated group, and let $\mc{H}$ be a
  collection of subgroups. We say that $(\Gamma, \mc{H})$ is a
  \emph{relatively hyperbolic pair} if $\Gamma$ is hyperbolic relative
  to $\mc{H}$.
\end{definition}

\subsection{The Bowditch boundary}

\begin{definition}
  Let $(\Gamma, \mc{H})$ be a relatively hyperbolic pair, so that
  $\mc{H}$ is the set of stabilizers of parabolic points for an action
  of $\Gamma$ on a metric space $Y$ as in
  \Cref{defn:relatively_hyperbolic_group}. We say that $Y$ is a
  \emph{Gromov model} for the pair $(\Gamma, \mc{H})$.
\end{definition}

In general there is \emph{not} a unique choice of Gromov model for a
given relatively hyperbolic pair $(\Gamma, \mc{H})$, even up to
quasi-isometry. There are various ``canonical'' constructions for a
preferred quasi-isometry class of Gromov model, with certain desirable
metric properties (see e.g. \cite{bowditch2012relatively},
\cite{gm2008dehn}).

Given \emph{any} two Gromov models $Y$, $Y'$ for $(\Gamma, \mc{H})$,
there is always a $\Gamma$-equivariant homeomorphism
$\dee Y \to \dee Y'$ \cite{bowditch2012relatively}. The $\Gamma$-space
$\dee Y$ is the \emph{Bowditch boundary} of $(\Gamma, \mc{H})$. We
will denote it by $\bgamh$, or sometimes just
$\dee \Gamma$ when the collection of peripheral subgroups is
understood from context. Since a Gromov model $Y$ is a proper
hyperbolic metric space, $\bgamh$ is always compact and
metrizable.

\begin{definition}
  We say a relatively hyperbolic pair $(\Gamma, \mc{H})$ is
  \emph{elementary} if $\Gamma$ is finite or virtually cyclic, or if
  $\mc{H} = \{\Gamma\}$.
\end{definition}

Whenever $(\Gamma, \mc{H})$ is nonelementary, its Bowditch boundary
contains at least three points. The convergence properties of the
action of $\Gamma$ on $\bgamh$ (see below) imply that in
this case, $\bgamh$ is \emph{perfect} (i.e. contains no
isolated points).

\subsubsection{Cocompactness on pairs}

Let $Y$ be a Gromov model for a relatively hyperbolic pair
$(\Gamma, \mc{H})$. Since $Y$ is hyperbolic, proper, and geodesic, for
any compact subset $K \subset Y$, the space of bi-infinite geodesics
passing through $K$ is compact.

Given any distinct pair of points $u, v \in \dee Y$, there is a
bi-infinite geodesic $c$ in $Y$ joining $u$ to $v$. Since a horoball
in a hyperbolic metric space has just one point in its ideal boundary,
this geodesic must pass through the thick part $Y_{\mathrm{th}}$ of
$Y$, so up to the action of $\Gamma$ it passes through a fixed compact
subset $K \subset Y_{\mathrm{th}}$.

This implies:
\begin{prop}
  \label{prop:pairs_cocompact}
  The action of $\Gamma$ on the space of distinct pairs in
  $\bgamh$ is cocompact.
\end{prop}

\subsection{Convergence group actions}
\label{subsec:convergence_groups}

If a group $\Gamma$ acts on a proper geodesic hyperbolic metric space
$Y$, we can characterize the geometrical finiteness of the action
entirely in terms of the topological dynamics of the action on
$\dee Y$. In particular, we can understand geometrical finiteness by
studying properties of \emph{convergence group actions}. See
\cite{tukia1994convergence}, \cite{tukia1998conical},
\cite{bowditch1999convergence} for further detail on such actions, and
justifications for the results stated in this section.

\begin{definition}
  \label{defn:convergence_group_properties}
  Let $\Gamma$ act as a convergence group (see
  \Cref{defn:convergence_group}) on a topological space $Z$.
  \begin{itemize}
  \item A point $z \in Z$ is a \emph{conical limit point} if there
    exists a sequence $\gamma_n \in \Gamma$ and distinct points $a, b
    \in Z$ such that $\gamma_n z \to a$ and $\gamma_n y \to b$ for any
    $y \ne z$.
  \item An infinite subgroup $H$ is a \emph{parabolic subgroup} if it
    fixes a point $p \in Z$, and every infinite-order element of $H$
    fixes exactly one point in $Z$.
  \item A point $p \in Z$ is a \emph{parabolic point} if it is the
    fixed point of a parabolic subgroup.
  \item A parabolic point $p$ is \emph{bounded} if its stabilizer
    $\Gamma_p$ acts cocompactly on $Z - \{p\}$.
  \end{itemize}
\end{definition}

The name ``conical limit point'' makes more sense in the context of
convergence group actions on boundaries of hyperbolic metric spaces.

\begin{definition}
  \label{defn:geometric_conical_limit}
  Let $Y$ be a hyperbolic metric space, and let $z \in \dee Y$. We say
  that a sequence $y_n \in Y$ \emph{limits conically to $z$} if there
  is a geodesic ray $c:\R^+ \to Y$ limiting to $z$ and a constant
  $D > 0$ such that
  \[
    d_Y(y_n, c(t_n)) < D
  \]
  for some sequence $t_n \to \infty$.
\end{definition}

A bounded neighborhood of a geodesic in a hyperbolic metric space
looks like a ``cone,'' hence ``conical limit.''

\begin{prop}[\cite{tukia1994convergence}, \cite{tukia1998conical}]
  \label{prop:conical_limit_point_interp}
  Let $\Gamma$ be a group acting properly discontinuously by
  isometries on a proper geodesic hyperbolic metric space $Y$, and fix
  a basepoint $y_0 \in Y$.

  Then $\Gamma$ acts on $\dee Y$ as a convergence group. Moreover, a
  point $z \in \dee Y$ is a conical limit point (in the dynamical
  sense of \Cref{defn:convergence_group_properties}) if and only if
  there is a sequence $\gamma_n \cdot y_0$ limiting conically to $z$
  (in the geometric sense of \Cref{defn:geometric_conical_limit}). In
  this case, there are distinct points $a, b \in \dee Y$ such that
  $\gamma_n^{-1} \cdot z \to a$ and $\gamma_n^{-1} \cdot z' \to b$ for
  any $z' \ne z$ in $\dee Y$.
\end{prop}

If $\gamma_n \cdot y_0$ limits conically to a point $z \in \dee Y$ for
some (hence any) basepoint $y_0 \in Y$, we just say that $\gamma_n$
limits conically to $z$.

\begin{theorem}[\cite{bowditch2012relatively}]
  Let $\Gamma$ be a finitely generated group acting by isometries on a
  hyperbolic metric space $Y$. Then $\Gamma$ is a relatively
  hyperbolic group, acting on $Y$ as in
  \Cref{defn:relatively_hyperbolic_group}, if and only if:
  \begin{enumerate}
  \item The induced action of $\Gamma$ on $\dee Y$ is a convergence
    group action.
  \item Every point $z \in \dee Y$ is either a conical limit point or
    a bounded parabolic point.
  \end{enumerate}
\end{theorem}

Whenever a group $\Gamma$ acts as a convergence group on a perfect
compact metrizable space $Z$, every point in $Z$ is either a conical
limit point or a bounded parabolic point, and the stabilizer of each
parabolic point is finitely generated, we say the $\Gamma$-action on
$Z$ is \emph{geometrically finite}. This is justified by a theorem of
Yaman \cite{yaman2004topological}, which says that any such group
action is induced by the action of a relatively hyperbolic group on a
Gromov model $Y$ whose boundary is equivariantly homeomorphic to
$Z$. We can then identify the space $Z$ with the Bowditch boundary
$\bgamh$. The set of parabolic points in $Z$ coincides exactly with
the set of fixed points of peripheral subgroups.

\begin{definition}
  Let $(\Gamma, \mc{H})$ be a relatively hyperbolic pair. We write
  \[
    \bgamh = \conbdry(\Gamma, \mc{H}) \sqcup \parbdry
    (\Gamma, \mc{H}),
  \]
  where $\conbdry(\Gamma, \mc{H})$ and $\parbdry(\Gamma, \mc{H})$ are
  respectively the conical limit points and parabolic points in
  $\bgamh$.
\end{definition}

\subsubsection{Compactification of $\Gamma$ and divergent sequences}

When $(\Gamma, \mc{H})$ is a relatively hyperbolic pair, there is a
natural topology on the set
\[
  \overline{\Gamma} = \Gamma \sqcup \bgamh
\]
making it into a \emph{compactification} of $\Gamma$
(i.e. $\overline{\Gamma}$ is compact, $\bgamh$ and $\Gamma$ are both
embedded in $\overline{\Gamma}$, and $\Gamma$ is an open dense subset
of $\overline{\Gamma}$). Specifically, we view $\Gamma$ as a subset of
(any) Gromov model $Y$, via an orbit map
$\gamma \mapsto \gamma \cdot y_0$ for some basepoint $y_0 \in
Y$. Since $\Gamma$ acts properly on $Y$, this is a proper embedding,
so if we compactify $Y$ by adjoining its visual boundary $\bgamh$, we
compactify $\Gamma$ as well; this does not depend on the choice of
basepoint $y_0$ or even the choice of space $Y$.

\begin{definition}
  A sequence $\gamma_n \in \Gamma$ is \emph{divergent} if it leaves
  every bounded subset of $\Gamma$ (equivalently, if a subsequence of
  it consists of pairwise distinct elements).
\end{definition}
Up to subsequence, a divergent sequence $\gamma_n \in \Gamma$
converges to a point $z \in \bgamh$. When $(\Gamma, \mc{H})$ is
non-elementary, the point $z$ is determined solely by the action of
$\Gamma$ on $\bgamh$: we have $\gamma_n \to z$ if and only if there is
at most one $x \in \bgamh$ such that $\gamma_n \cdot x$ does not
converge to $z$ (this is a consequence of \cite[Theorem
3A]{tukia1994convergence}).

\subsection{The coned-off Cayley graph}

Whenever $(\Gamma, \mc{H})$ is a relatively hyperbolic pair, there are
only finitely many conjugacy classes of groups in $\mc{H}$. We can fix
a finite set $\mc{P}$ of conjugacy representatives for the groups in
$\mc{H}$. The set $\mc{P}$ corresponds to a finite set
$\Pi \subset \parbdry \Gamma$ of parabolic points, such that
\[
  \mc{P} = \{\Gamma_p : p \in \Pi\}.
\]
Then $\Pi$ contains exactly one point in each $\Gamma$-orbit in
$\parbdry\Gamma$.

\begin{definition}
  \label{defn:conedgraph}
  Let $(\Gamma, \mc{H})$ be a relatively hyperbolic pair, and fix a
  finite generating set $S$ for $\Gamma$ and finite collection of
  conjugacy representatives $\mc{P}$ for $\mc{H}$.

  The \emph{coned-off Cayley graph} $\conedgraph$ is a metric space
  obtained from the Cayley graph $\cay(\Gamma, S)$ as follows: for
  each coset $gP_i$ for $P_i \in \mc{P}$, we add a vertex
  $v(gP_i)$. Then, we add an edge of length 1 from each $h \in gP_i$
  to $v(gP_i)$.
\end{definition}

The quasi-isometry class of $\conedgraph$ is independent of the choice
of generating set $S$. When $(\Gamma, \mc{H})$ is a relatively
hyperbolic pair, $\conedgraph$ is a hyperbolic metric space. It is
\emph{not} a proper metric space if $\mc{H}$ is nonempty.


\section{Lie theory notation and background}
\label{sec:lie_theory_background}

For the rest of the paper, we let $G$ be a connected semisimple Lie
group with no compact factor and finite center. We will be concerned
with representations $\rho:\Gamma \to G$, where $\Gamma$ is a
relatively hyperbolic group. We want to consider the action of
$\rho(\Gamma)$ on the flag manifold $\flags$, where $P$ is a parabolic
subgroup of $G$.

In this section, we give an overview of the definitions and notation
we will use to describe the dynamical behavior of the $\Gamma$-action
on $\flags$. We mostly follow the notation of \cite{ggkw2017anosov},
but we will also identify the connection to the language of
\cite{klp2017anosov}.

The exposition here is fairly brief, since most of this paper does not
use much of the technical theory of semisimple Lie groups and their
associated Riemannian symmetric spaces. In fact, in nearly every case,
our approach will be to use a representation of $G$ to reduce to the
case $G = \PGL(n, \R)$. The most important part of this section is
\ref{sec:flag_dynamics_background}, which identifies the connection
between \emph{$P$-divergence} (or equivalently
\emph{$\taumod$-regularity}) and \emph{contracting dynamics} in $G/P$.

Standard references for the general theory are \cite{eberlein},
\cite{helgason}, and \cite{knapp}. See also section 3 of
\cite{riestenberg2021quantified} for a careful discussion of the
theory as it relates to Anosov representations and the work of
Kapovich-Leeb-Porti.

\subsection{Parabolic subgroups}

Let $K$ be a maximal compact subgroup of the semisimple Lie group $G$,
and let $X$ be the Riemannian symmetric space $G/K$. A subgroup
$P \subset G$ is a \emph{parabolic subgroup} if it is the stabilizer
of a point in the visual boundary $\dee_\infty X$ of $X$. Two
parabolic subgroups $P, Q$ are \emph{opposite} if there is a
bi-infinite geodesic $c$ in $X$ so that $P$ is the stabilizer of
$c(\infty)$ and $Q$ is the stabilizer of $c(-\infty)$.

The compact homogeneous $G$-space $\flags$ is called a \emph{flag
  manifold}. If $P$ and $Q$ are parabolic subgroups, then we say that
two flags $\xi^+ \in G/P$ and $\xi^- \in G/Q$ are \emph{opposite} if
the stabilizers of $\xi^+$, $\xi^-$ are opposite parabolic
subgroups. (In particular a conjugate of $Q$ must be opposite to $P$).

\subsection{Root space decomposition}

Let $\mf{g}$ be the Lie algebra of $G$, and let $\mf{k}$ be the Lie
algebra of the maximal compact $K$. We can decompose $\mf{g}$ as
$\mf{k} \oplus \mf{p}$, and fix a maximal abelian subalgebra
$\mf{a} \subset \mf{p}$. The restriction of the Killing form $B$ to
$\mf{p}$ is positive definite, so any maximal abelian
$\mf{a} \subset \mf{p}$ is naturally endowed with a Euclidean
structure.

Each element of the abelian subalgebra $\mf{a}$ acts semisimply on
$\mf{g}$, with real eigenvalues. So we let $\Sigma \subset \mf{a}^*$
denote the set of \emph{roots} for this choice of $\mf{a}$, i.e. the
set of nonzero linear functionals $\alpha \in \mf{a}^*$ such that the
linear map $\mf{g} \to \mf{g}$ given by $X - \alpha(X)I$ has nonzero
kernel for every $X \in \mf{a}$. We have a \emph{restricted root space
  decomposition}
\[
  \mf{g} = \mf{g}_0 \oplus \bigoplus_{\alpha \in \Sigma}
  \mf{g}_\alpha,
\]
where $X \in \mf{a}$ acts on $\mf{g}_\alpha$ by multiplication by
$\alpha(X)$.

We choose a set of \emph{simple roots} $\Delta \subset \Sigma$ so that
each $\alpha \in \Sigma$ can be uniquely written as a linear
combination of elements of $\Delta$ with coefficients either all
nonnegative or all nonpositive. We let $\Sigma_+$ denote the
\emph{positive roots}, i.e. roots which are nonnegative linear
combinations of elements of $\Delta$.

The simple roots $\Delta$ determine a \emph{Euclidean Weyl chamber}
\[
  \mf{a}^+ = \{x \in \mf{a} : \alpha(x) \ge 0, \textrm{ for all }
  \alpha \in \Delta\}.
\]
The kernels of the roots $\alpha \in \Delta$ are the \emph{walls} of
the Euclidean Weyl chamber.

Choosing a maximal compact $K$, a maximal abelian
$\mf{a} \subset \mf{p}$, and a Euclidean Weyl chamber $\mf{a}^+$
determines a \emph{Cartan projection}
\[
  \mu:G \to \mf{a}^+,
\]
uniquely determined by the equation $g = k\exp(\mu(g))k'$, where
$k, k' \in K$ and $\mu(g) \in \mf{a}^+$.

\subsection{$P$-divergence}

Fix a subset $\theta$ of the simple roots $\Delta$. We define a
\emph{standard parabolic subgroup} $\stdpargroup$ to be the normalizer
of the Lie algebra
\[
  \bigoplus_{\alpha \in \Sigma^+_\theta} \mf{g}_\alpha,
\]
where $\Sigma^+_\theta$ is the set of positive roots which are
\emph{not} in the span of $\Delta - \theta$. The \emph{opposite
  subgroup} $\stdoppgroup$ is the normalizer of
\[
  \bigoplus_{\alpha \in \Sigma^+_\theta} \mf{g}_{-\alpha}.
\]

Every parabolic subgroup $P \subset G$ is conjugate to a unique
standard parabolic subgroup $\stdpargroup$, and every pair of opposite
parabolics $(P^+, P^-)$ is simultaneously conjugate to a unique pair
$(\stdpargroup, \stdoppgroup)$.

For a fixed $\theta \subset \Delta$, the group $\stdpargroup$ is the
stabilizer of the endpoint of a geodesic ray $\exp(t Z) \cdot p$,
where $p \in X$ is the image of the identity in $G/K$, and for any
$\alpha \in \Delta$, the element $Z \in \mf{a}^+$ satisfies
\[
  \alpha(Z) = 0 \iff \alpha \in \Delta - \theta.
\]

\begin{definition}
  Let $g_n$ be a sequence in $G$. The sequence $g_n$ is
  \emph{$\stdpargroup$-divergent} if for every $\alpha \in \theta$, we
  have
  \[
    \alpha(\mu(g_n)) \to \infty.
  \]

  That is, the Cartan projections of the sequence $g_n$ \emph{drift
    away from the walls} of $\mf{a}$ determined by the subset
  $\theta \subset \Delta$.

  For a general parabolic subgroup $P \subset G$, we say that $g_n$ is
  \emph{$P$-divergent} if $g_n$ is $\stdpargroup$-divergent for
  $\stdpargroup$ conjugate to $P$.
\end{definition}

\subsection{Affine charts}

\begin{definition}
  \label{defn:affine_chart}
  Let $\pargroup$, $\oppgroup$ be opposite parabolic subgroups in
  $G$. Given a flag $\xi \in \oppflags$, we define
  \[
    \Opp(\xi) = \{\eta \in \pflags : \xi \textrm{ is opposite to }
    \eta\}.
  \]
  We call a set of the form $\Opp(\xi)$ for some $\xi \in \oppflags$
  an \emph{affine chart} in $\pflags$.
\end{definition}

An affine chart is the unique open dense orbit of $\Stab_G(\xi)$ in
$\pflags$. When $G = \PGL(d, \R)$ and $\pargroup$ is the stabilizer of
a line $\ell \subset \R^d$, $\pflags$ is identified with $\P(\R^d)$
and this notion of affine chart agrees with the usual one in
$\P(\R^d)$.

\subsection{Dynamics in flag manifolds}
\label{sec:flag_dynamics_background}

There is a close connection between $P$-divergence in the group $G$
and the \emph{topological dynamics} of the action of $G$ on the
associated flag manifold $G/P$. Kapovich-Leeb-Porti frame this
connection in terms of a \emph{contraction property} for $P$-divergent
sequences.

\begin{definition}[\cite{klp2017anosov}, Definition 4.1]
  Let $g_n$ be a sequence of group elements in $G$. We say that $g_n$
  is \emph{$\pargroup$-contracting} if there exist $\xi \in \pflags$,
  $\xi_- \in \oppflags$ such that $g_n$ converges uniformly to $\xi$
  on compact subsets of $\Opp(\xi_-)$.

  The flag $\xi$ is the uniquely determined \emph{limit} of the
  sequence $g_n$.
\end{definition}

\begin{definition}
  \label{defn:limit_set}
  For an arbitrary sequence $g_n \in G$, a \emph{$\pargroup$-limit
    point} of $g_n$ in $\pflags$ is the limit point of some
  $\pargroup$-contracting subsequence of $g_n$.

  The \emph{$\pargroup$-limit set} of a group $\Gamma \subset G$ is the
  set of $\pargroup$-limit points of sequences in $\Gamma$.
\end{definition}

The importance of contracting sequences is captured by the following:
\begin{prop}[\cite{klp2017anosov}, Proposition 4.15]
  \label{prop:divergence_equivalent_contracting}
  A sequence $g_n \in G$ is $P^+$-divergent if and only if every
  subsequence of $g_n$ has a $P^+$-contracting subsequence.
\end{prop}

\Cref{prop:divergence_equivalent_contracting} implies in particular
that if $g_n \in G$ is $\pargroup$-divergent, then up to subsequence
there is an open subset $U \subset \pflags$ such that $g_n \cdot U$
converges to a singleton in $\pflags$. It turns out that this ``weak
contraction property'' is enough to characterize $P^+$-divergence.

\begin{restatable}{prop}{localContractingImpliesGlobalContracting}
  \label{prop:local_contracting_implies_global_contracting}
  Let $g_n$ be a sequence in $G$, and suppose that for some nonempty
  open subset $U \subset \pflags$, we have $g_n \cdot U \to \{\xi\}$
  for $\xi \in \pflags$. Then $g_n$ is $\pargroup$-divergent, and has
  a unique $\pargroup$-limit point $\xi \in \pflags$.
\end{restatable}

We provide a proof of this fact in \Cref{sec:contraction_appendix}.

\subsubsection{Dynamics of inverses of $\pargroup$-divergent
  sequences}

When $g_n$ is a $\pargroup$-divergent sequence, the inverse sequence
is $P^-$-divergent. Kapovich-Leeb-Porti show that this can be framed
in terms of the dynamical behavior of the inverse sequence.
\begin{lem}[\cite{klp2017anosov}, Lemma 4.19]
  \label{lem:klp_contraction_symmetry}
  For $g_n \in G$ and flags $\xi_- \in \oppflags, \xi_+ \in \pflags$,
  the following are equivalent:
  \begin{enumerate}
  \item $g_n$ is $P^+$-contracting and $g_n|_{\Opp(\xi_-)} \to \xi_+$
    uniformly on compacts.
  \item $g_n$ is $P^+$-divergent, $g_n$ has unique $P^+$-limit point
    $\xi_+$, and $g_n^{-1}$ has unique $P^-$-limit point $\xi_-$.
  \end{enumerate}
\end{lem}

\subsection{$\taumod$-regularity}
\label{sec:taumod_regular}

$P$-divergent sequences are equivalent to the \emph{$\taumod$-regular}
sequences discussed in the work of Kapovich-Leeb-Porti, where
$\taumod$ is the unique face corresponding to $P$ in a \emph{spherical
  model Weyl chamber}. We explain the connection here.

\begin{remark}
  The language of $\taumod$-regularity is not used anywhere else in
  this paper, so this part of the background is provided for
  convenience only and may be safely skipped.
\end{remark}

For any point $p \in X$, we let $\mf{p}$ be the uniquely determined
subspace of $\mf{g}$ such that $\mf{g} = \mf{k} \oplus \mf{p}$, where
$\mf{k}$ is the Lie algebra of the stabilizer of $p$ in $G$.

Let $z \in \dee_\infty X$. There is a point $p \in X$, a maximal
abelian subalgebra $\mf{a} \subset \mf{p}$, a Euclidean Weyl chamber
$\mf{a}^+ \subset \mf{a}$, and a unit-length $Z \in \mf{a}^+$ such
that $z$ is the endpoint of the geodesic ray
$c(t) = \exp(tZ) \cdot p$.

Up to the action of the stabilizer of $z$, the point $p$, the maximal
abelian subalgebra $\mf{a}$, the Euclidean Weyl chamber $\mf{a}^+$,
and the unit vector $Z \in \mf{a}^+$ are uniquely determined. In
addition, the stabilizer in $G$ of the triple $(p, \mf{a}, \mf{a}^+)$
acts trivially on $\mf{a}^+$.

This means that we can identify the space $\dee_\infty X / G$ with the
set of unit vectors in any Euclidean Weyl chamber $\mf{a}^+$. This set
has the structure of a \emph{spherical simplex}. We let $\sigmamod$
denote the \emph{model spherical Weyl chamber} $\dee_\infty X / G$.

We let $\pi:\dee_\infty X \to \sigmamod$ be the \emph{type map} to the
model spherical Weyl chamber. For fixed $z \in \dee_\infty X$, we let
$P_z$ denote the parabolic subgroup stabilizing $z$.

After choosing a maximal compact $K$, a maximal abelian
$\mf{a} \subset \mf{p}$, and a Euclidean Weyl chamber $\mf{a}^+$, the
data of a \emph{face} $\taumod$ of the spherical simplex $\sigmamod$
is the same as the data of a \emph{subset} of the simple roots of $G$:
the set of roots identifies a collection of walls of the Euclidean
Weyl chamber $\mf{a}^+$. The intersection of those walls with the unit
sphere in $\mf{a}$ is uniquely identified with a face of $\sigmamod$.

\begin{definition}
  Let $\taumod$ be a face of the model spherical Weyl chamber
  $\sigmamod$. We say that a sequence $g_n \in G$ is
  \emph{$\taumod$-regular} if $g_n$ is $P_z$-divergent for some
  $z \in \dee_\infty X$ such that $\pi(z)$ lies in the relative
  interior $\mathrm{int}(\taumod)$ of $\taumod$.
\end{definition}

For a fixed model face $\taumod \subset \sigmamod$, we let $P_\taumod$
denote any parabolic subgroup which is the stabilizer of a point
$z \in \pi^{-1}(\mathrm{int}(\taumod))$. All such parabolic subgroups
are conjugate, so as a $G$-space the flag manifold $G/P_\taumod$
depends only on the model face $\taumod$.


\section{EGF representations and relative Anosov representations}
\label{sec:egf_definition}
\label{sec:egf_representations}

In this section we cover basic properties of the central objects of
this paper: \emph{extended geometrically finite representations} from
a relatively hyperbolic group $\Gamma$ to a semisimple Lie group $G$
with no compact factor and trivial center. We also show that they
generalize a definition of relative Anosov representation
(\Cref{thm:bdry_dynamics_homeo_equiv_rel_asymp_embedded}), and prove
our Anosov relativization theorem
(\Cref{thm:changing_peripheral_structure}).

We refer also to Section 2 of the related paper \cite{WeismanExamples}
for an overview of the definition in the special case where
$G = \PGL(d, \R)$ or $\SL(d, \R)$ and the parabolic subgroup $P$ is
the stabilizer of a flag of type $(1, d-1)$ in $\R^d$.

\begin{definition}
  Let $P$ be a parabolic subgroup of $G$. We say that $P$ is
  \emph{symmetric} if $P = P^+$ is conjugate to a subgroup $P^-$
  opposite to $P$.

  When $P = P^+$ is symmetric, we can identify $\pflags$ with
  $\oppflags$, so that it makes sense to say that two flags $\xi_1,
  \xi_2 \in \flags$ are \emph{opposite}.
\end{definition}

For the rest of the section, we fix a relatively hyperbolic pair
$(\Gamma, \mc{H})$ and a symmetric parabolic subgroup $P$ in a
semisimple Lie group $G$.

\begin{definition}
  Let $A, B$ be two subsets of $\flags$. We say that $A$ and $B$ are
  \emph{opposite} if every $\xi \in A$ is opposite to every
  $\nu \in B$.
\end{definition}

\begin{definition}
  Let $\Lambda \subset \flags$. We say that a continuous surjective
  map $\phi:\Lambda \to \bgamh$ is \emph{antipodal} if for every pair
  of distinct points $z_1, z_2 \in \bgamh$, $\phi^{-1}(z_1)$ is
  opposite to $\phi^{-1}(z_2)$.
\end{definition}

We recall the main definition of the paper here:

\EGF*

\begin{remark}
  \label{rem:boundary_set_non_canonical}
  Unfortunately, the boundary set $\Lambda \subset \flags$ is
  \emph{not} necessarily uniquely determined by the representation
  $\rho$. In many contexts, we will be able to make a natural choice,
  but we do not give a procedure for doing so in general.
\end{remark}

\subsection{Discreteness and finite kernel}
\label{sec:discrete_finite_kernel}

When $\rho:\Gamma \to G$ is EGF, the action of $\rho(\Gamma)$ on the
boundary set $\Lambda$ is by definition an extension of the
topological dynamical system $(\Gamma, \bgamh)$. When $\Gamma$ is
non-elementary, convergence dynamics imply that the homomorphism
$\Gamma \to \Homeo(\bgamh)$ has finite kernel and discrete image. So
the map $\Gamma \to \Homeo(\Lambda)$ must also have discrete image and
finite kernel, and therefore so does the representation
$\rho:\Gamma \to G$. The case where $\Gamma$ is elementary can be
verified directly.

\subsection{Uniformity of repelling strata}
\label{sec:uniform_repelling_strata}

An earlier version of this paper gave a slightly weaker version of the
definition of an extended convergence group action. Specifically, in
place of condition \ref{item:strata_interior} in
\Cref{defn:extends_convergence_dynamics}, we gave the weaker
condition:
\begin{enumerate}[(C1*)]
\item \label{item:weak_strata_interior} For every $z \in \bgamh$, the
  set $C_z$ contains an open neighborhood of
  $\Lambda - \phi^{-1}(\{z\})$.
\end{enumerate}

It is not clear whether or not \ref{item:weak_strata_interior} implies
\ref{item:strata_interior} in the general context of extended
convergence group actions. However, the conditions are equivalent in
the context of discrete group actions on flag manifolds if we
additionally impose an assumption on $\rho$ called \emph{relative
  $P$-divergence}. To define this, let $\Pi \subset \bgamh$ be a
finite set containing exactly one point in each parabolic orbit in
$\bgamh$, and let $\pargps$ denote the set of $\Gamma$-stabilizers of
points in $\Pi$.

\begin{definition}
  Let $\rho$ be a representation of a relatively hyperbolic group
  $\Gamma$ with peripheral subgroups $\mc{H}$. The representation
  $\rho$ is \emph{relatively $P$-divergent} if, for every sequence of
  elements $\gamma_n \in \Gamma$ which is unbounded in the coned-off
  Cayley graph $\conedgraph$, the sequence $\rho(\gamma_n)$ is
  $P$-divergent.
\end{definition}

Below, we prove the following:

\begin{prop}
  \label{prop:weak_egf_to_egf}
  Let $\rho:\Gamma \to G$ be a representation, and suppose that there
  exists a compact $\rho$-invariant set $\Lambda \subset \flags$, an
  equivariant surjective antipodal map $\phi:\Lambda \to \bgamh$, and
  open sets $C_z \subset \flags$ for each $z \in \bgamh$ satisfying
  conditions \ref{item:weak_strata_interior} and
  \ref{item:extended_convergence}. If $\rho$ is relatively
  $P$-divergent, then $\rho$ is EGF.
\end{prop}

We will show later (see \Cref{cor:egf_rel_divergent}) that the
converse to \Cref{prop:weak_egf_to_egf} also holds, meaning that any
EGF representation is relatively
$P$-divergent. \Cref{prop:weak_egf_to_egf} itself will largely be a
consequence of the following more technical statement:


\begin{prop}
  \label{prop:weak_egf_to_egf_ii}
  Let $\rho:\Gamma \to G$ be a representation, and suppose that there
  exists a compact $\rho$-invariant set $\Lambda \subset \flags$, an
  equivariant surjective antipodal map $\phi:\Lambda \to \bgamh$, and
  open sets $C_z \subset \flags$ for each $z \in \bgamh$ satisfying:
  \begin{enumerate}
    \custitem[(C1*)] For every $z \in \bgamh$, the set $C_z$ contains
    an open neighborhood of $\Lambda - \phi^{-1}(\{z\})$.
    
    \custitem[(C2-par)]\label{item:parabolic_flag_convergence} For any
    parabolic point $p \in \bgamh$, any compact $K \subset C_p$, and
    any open set $U$ containing $\phi^{-1}(p)$, for all but finitely
    many $\gamma \in \Gamma_p = \Stab_\Gamma(p)$, we have
    $\rho(\gamma)K \subset U$.
  
    \custitem[(P-div)]\label{item:rel_p_div} The representation $\rho$
    is relatively $P$-divergent, and for every sequence $\gamma_n$
    tending to infinity in $\conedgraph$, the $P$-limit points of
    $\rho(\gamma_n)$ all lie in $\Lambda$.
  \end{enumerate}
  Then there is a subset $\hat{\Lambda} \subset \Lambda$ such that
  $\rho$ is EGF, with boundary extension $\phi|_{\hat{\Lambda}}$.
\end{prop}

\begin{proof}[Proof of \Cref{prop:weak_egf_to_egf}, assuming
  \Cref{prop:weak_egf_to_egf_ii}]
  Assume the hypotheses of \Cref{prop:weak_egf_to_egf} for a
  representation $\rho$ and map $\phi:\Lambda \to \bgamh$. Note that
  condition \ref{item:parabolic_flag_convergence} in
  \Cref{prop:weak_egf_to_egf_ii} is a special case of
  \ref{item:extended_convergence}, so we will be done if we can show
  that condition \ref{item:rel_p_div} holds for $\Lambda$ and
  $\rho$. Fix an arbitrary sequence $\gamma_n$ tending to infinity in
  $\conedgraph$, and let $\xi^+$ be a $P$-limit point of
  $\gamma_n$. We wish to show that $\xi^+ \in \Lambda$.

  By definition, we can find a flag $\xi^-$ and extract a subsequence
  so that $\rho(\gamma_n)$ converges to $\xi^+$, uniformly on compact
  subsets of $\Opp(\xi^-)$. Extract a further subsequence so that
  $\gamma_n^{\pm 1}$ converges to $z_\pm \in \bgamh$. Since
  $\Opp(\xi^-)$ is dense in $\flags$, it has nonempty intersection
  with $C_{z_-}$, so let $K$ be a nonempty compact subset in this
  intersection. By \ref{item:extended_convergence}, $\rho(\gamma_n)K$
  eventually lies in an arbitrary neighborhood of
  $\phi^{-1}(z_+) \subset \Lambda$. But we also know that
  $\rho(\gamma_n)K$ converges to the singleton $\{\xi^+\}$, hence
  $\xi^+ \in \Lambda$.
\end{proof}

We now turn to the proof of \Cref{prop:weak_egf_to_egf_ii}. Fix a
representation $\rho$, map $\phi:\Lambda \to \bgamh$, and family of
sets $C_z$ as in the statement of the proposition. For the proof of
the proposition, we need to construct a new surjective antipodal map
$\hat{\phi}:\hat{\Lambda} \to \bgamh$ and strata $\hat{C}_z$
satisfying conditions \ref{item:strata_interior} and
\ref{item:extended_convergence} in
\Cref{defn:extends_convergence_dynamics}.

We first construct the set $\hat{\Lambda}$. Using condition
\ref{item:rel_p_div}, we will shrink the set $\Lambda$ into a certain
nice form, while still preserving all of the assumptions in
\Cref{prop:weak_egf_to_egf_ii}. For each conical limit point
$z \in \bgamh$, define $\Psi_z \subset \flags$ to be the union of the
set of $P$-limit points of all sequences $\rho(\gamma_n)$, where
$\gamma_n \to z$ in the Bowditch compactification of $\Gamma$. (Note
that since $z$ is conical, any sequence $\gamma_n$ converging to $z$
must tend to infinity in $\conedgraph$, and therefore by
\ref{item:rel_p_div} $\rho(\gamma_n)$ is $P$-divergent and $\Psi_z$ is
nonempty.)

Then, define a subset $\hat{\Lambda} \subset \flags$ via:
\[
  \hat{\Lambda} = \left(\bigcup_{\substack{z \in \bgamh,\\z \textrm{
          parabolic}}}\phi^{-1}(z)\right) \cup
  \left(\bigcup_{\substack{z \in \bgamh,\\z \textrm{ conical}}}
    \Psi_z\right).
\]

\begin{lem}
  The set $\hat{\Lambda}$ is a compact $\rho$-invariant subset of
  $\Lambda$.
\end{lem}
\begin{proof}
  That $\hat{\Lambda} \subset \Lambda$ follows from assumption
  \ref{item:rel_p_div}. The $\rho$-invariance of $\hat{\Lambda}$
  follows from $\rho$-equivariance of $\phi$, as well as the
  $\Gamma$-invariance of the set of conical limit points in
  $\bgamh$. So, we just need to show that $\hat{\Lambda}$ is compact.

  Consider a sequence of flags $\xi_n \in \hat{\Lambda}$, and extract
  a subsequence so that $\xi_n \to \xi$. We wish to show
  $\xi \in \hat{\Lambda}$. Since $\hat{\Lambda} \subset \Lambda$,
  there is a sequence of points $z_n \in \bgamh$ such that
  $\xi_n \in \phi^{-1}(z_n)$. By compactness of $\Lambda$ and
  continuity of $\phi$, we know $\xi \in \Lambda$ and
  $z_n \to z = \phi(\xi)$. If $z$ is parabolic, then
  $\phi^{-1}(z) \subset \hat{\Lambda}$, hence $\xi \in \hat{\Lambda}$,
  so assume $z$ is conical.

  Let $\gamma_m$ be a sequence in $\Gamma$ converging conically to
  $z$; in particular, $\gamma_m$ must tend to infinity in
  $\conedgraph$, so applying assumption \ref{item:rel_p_div} means
  that, after extraction, there are flags $\xi^\pm \in \Lambda$ so
  that $\rho(\gamma_m)$ converges to the constant map $\xi^+$,
  uniformly on compacts in $\Opp(\xi^-)$. Note that since $\gamma_m$
  converges to the conical point $z$, we also have
  $\xi^+ \in \hat{\Lambda}$ by definition.

  Let $v = \phi(\xi^-)$. Recall that since $\gamma_m$ converges
  conically to $z$, there are distinct points $a, b \in \bgamh$ such
  that $\gamma_m^{-1}z \to a$ and $\gamma_m^{-1}y \to b$ for all
  $y \ne z$. We claim that we can find a sequence of indices $m_n$
  such that $m_n \to \infty$, but $\gamma_{m_n}^{-1}z_n$ does
  \emph{not} accumulate at $v$. There are two cases: either $v = a$,
  or $v \ne a$. In the case $v = a$, choose $m_n$ large enough so that
  the sequence $\gamma_{m_n}^{-1}z_n \to b$; in the case $v \ne a$,
  since $z_n \to z$, we can choose $m_n$ so that
  $\gamma_{m_n}^{-1}z_n \to a$.

  In either case, by continuity, equivariance, and antipodality of
  $\phi$, the sequence
  \[
    \rho(\gamma_{m_n}^{-1}\xi_n) \in \phi^{-1}(\gamma_{m_n}^{-1}z_n)
  \]
  lies in a compact subset of $\Opp(\xi^-)$. Therefore
  $\xi_n = \rho(\gamma_{m_n})\rho(\gamma_{m_n}^{-1})\xi_n$ must
  converge to $\xi^+$, hence $\xi = \xi^+$ and thus
  $\xi \in \hat{\Lambda}$.
\end{proof}

Restrict $\phi$ to obtain a new map
$\hat{\phi}:\hat{\Lambda} \to \bgamh$. This map is also equivariant,
surjective, and antipodal. Next, we define a new family $\hat{C}_z$ of
repelling strata. For each conical limit point $z \in \bgamh$, define
\[
  \hat{C}_z = \Opp(\hat{\phi}^{-1}(z)).
\]
Then, for each parabolic point $p \in \Pi$, let $\Gamma_p$ be the
stabilizer of $p$ in $\Gamma$, and define
\[
  \hat{C}_p = \Opp(\hat{\phi}^{-1}(p)) \cap \bigcup_{\gamma \in
    \Gamma_p}\rho(\gamma)C_p.
\]
Finally, for an arbitrary parabolic point $z \in \bgamh$, write
$z = gp$ for $p \in \Pi$, $g \in \Gamma$, and define
$\hat{C}_z = \rho(g)\hat{C}_p$.

\begin{lem}
  \label{lem:hats_weak_convergence}
  The map $\hat{\phi}$ and sets $\hat{C}_z$ satisfy conditions
  \ref{item:weak_strata_interior} and
  \ref{item:parabolic_flag_convergence}.
\end{lem}
\begin{proof}
  The fact that $\hat{\phi}$ is antipodal says precisely that
  $\hat{\Lambda} - \hat{\phi}^{-1}(z)$ lies in
  $\Opp(\hat{\phi}^{-1})(z)$, so each $\hat{C}_z$ is a union of open
  sets satisfying \ref{item:weak_strata_interior} and therefore also
  satisfies \ref{item:weak_strata_interior}. To see that
  \ref{item:parabolic_flag_convergence} holds, let $K$ be a compact
  subset of $\hat{C}_p$ for some $p \in \Pi$ parabolic, and let
  $\gamma_n$ be an unbounded sequence in the $\Gamma$-stabilizer of
  $p$. By compactness, $K$ is contained in a finite union of sets of
  the form $\rho(h)C_p$, and in fact $K$ can be written as a finite
  union $K_1, \ldots, K_m$ of compact pieces, each of which lies in
  some $\rho(h_i)C_p$. Each $h_i$ is fixed, so $\gamma_nh_i$ converges
  to $p$, and therefore
  \[
    \rho(\gamma_n)K_i = \rho(\gamma_n h_i) \rho(h_i^{-1})K_i
  \]
  eventually lies in an arbitrarily small neighborhhod of
  $\phi^{-1}(p) = \hat{\phi}^{-1}(p)$, and therefore the same is true
  for $K$.
\end{proof}

Note that we have not verified that the set $\hat{\Lambda}$ satisfies
property \ref{item:rel_p_div}. It turns out that this is true, but the
proof will come much later (see
\Cref{lem:automaton_divergent_limits}). However, we can still use the
fact that $\rho$ is relatively $P$-divergent throughout the arguments
below.

We strengthen the second part of the above lemma to:
\begin{prop}
  \label{prop:weak_to_convergence}
  The map $\hat{\phi}:\hat{\Lambda} \to \bgamh$ and sets $\hat{C}_z$
  satisfy condition \ref{item:extended_convergence} in
  \Cref{defn:extends_convergence_dynamics}.
\end{prop}
\begin{proof}
  Consider a sequence $\gamma_n \in \Gamma$ such that
  $\gamma_n^{\pm 1} \to z_\pm \in \bgamh$, and let
  $K \subset \hat{C}_{z_-}$ be compact. Fix an arbitrary subsequence
  of $\gamma_n$. We will show that there is a further subsequence so
  that $\rho(\gamma_n)K$ lies in an arbitrarily small neighborhood of
  $\hat{\phi}^{-1}(z_+)$. After extraction we can assume that either
  $\gamma_n$ tends to infinity in $\conedgraph$, or else $\gamma_n$ is
  bounded in $\conedgraph$.

  In the first case, the sequences $\rho(\gamma_n^{\pm 1})$ are both
  $P$-divergent by assumption \ref{item:rel_p_div}. Using
  \Cref{lem:klp_contraction_symmetry}, extract a subsequence so that
  $\rho(\gamma_n^{\pm 1})$ converges to a flag $\xi_{\pm}$, uniformly
  on compact subsets of $\Opp(\xi_\mp)$; again by the assumption
  \ref{item:rel_p_div}, $\xi_\pm$ both lie in $\Lambda$. We may assume
  that $\Gamma$ is non-elementary, so there are at least two points
  $y, y' \in \bgamh$ not equal to $\phi(\xi_-)$. Then, since
  $\phi:\Lambda \to \bgamh$ is antipodal, we have
  $\phi^{-1}(\{y, y'\}) \subset \Opp(\xi_-)$. Then, for any flags
  $\eta \in \phi^{-1}(y), \eta' \in \phi^{-1}(y')$, both of
  $\rho(\gamma_n)\eta, \rho(\gamma_n)\eta'$ converge to $\xi_+$, and
  by equvariance and continuity of $\phi$, it follows that
  $\gamma_ny, \gamma_ny'$ both converge to $\phi(\xi_+)$. Since
  $\gamma_n \to z_+$ we must have $\xi_+ \in \phi^{-1}(z_+)$. If $z_+$
  is conical, then by construction we have
  $\xi_+ \in \Psi_{z_+} = \hat{\phi}^{-1}(z_+)$. Otherwise, $z_+$ is
  parabolic and $\hat{\phi}^{-1}(z_+) = \phi^{-1}(z_+)$.

  Applying an identical argument to the sequence $\gamma_n^{-1}$ shows
  that $z_- = \hat{\phi}(\xi_-)$, and thus, by construction,
  $\hat{C}_{z_-} \subseteq \Opp(\xi_-)$. In particular, since the
  compact set $K$ lies in $\hat{C}_{z_-} \subseteq \Opp(\xi_-)$, the
  set $\rho(\gamma_n)K$ eventually lies in an arbitrarily small
  neighborhood of $\hat{\phi}^{-1}(z_+)$, as required for this case.

  We now turn to the case where $\gamma_n$ is bounded in
  $\conedgraph$. After extracting a subsequence, we can fix some
  $k \in \N$ and write each $\gamma_n$ as an alternating product
  \[
    g_1h_1^{(n)} \ldots g_kh_k^{(n)}g_{k+1},
  \]
  where, for each $1 \le i \le k$, $g_i$ is a fixed element in
  $\Gamma$, and $h_i^{(n)}$ lies in a fixed parabolic subgroup
  $H_i \in \mc{P}$. Let $p_i \in \Pi$ be the parabolic point fixed by
  $H_i$. We can assume that $k$ has been chosen minimally, which
  ensures that $g_{i+1}p_{i+1} \ne p_i$ for every $i$.
  
  We claim that $\gamma_n$ converges to $z_+ = g_1p_1$,
  $\gamma_n^{-1}$ converges to $z_- = g_{k+1}^{-1} p_k$, and for any
  compact $K \subset \hat{C}_{z_-}$ and open $U$ containing
  $\hat{\phi}^{-1}(z_+)$, for large $n$, we have
  $\rho(\gamma_n)K \subset U$; we will prove this by inducting on $k$.

  When $k = 1$, then $p = p_1 = p_k$, and $\gamma_n = g_1h_ng_2$ for
  $h_n \in \Gamma_p$ and $g_1, g_2 \in \Gamma$ fixed. The distance
  between $h_ng_2$ and $h_n$ is bounded in any word metric on
  $\Gamma$, so $h_ng_2$ converges to $p$ in $\overline{\Gamma}$ and
  $g_1h_ng_2$ converges to $g_1p = z_+$. Then, if $K$ is a compact
  subset of $\hat{C}_{z_-} = \hat{C}_{g_2^{-1}p}$, the set
  $\rho(g_2)K$ is a compact subset of $\hat{C}_p \subset C_p$. Thus by
  \Cref{lem:hats_weak_convergence} we see that $\rho(h_ng_2)K$
  eventually lies in an arbitrarily small neighborhood of
  $\hat{\phi}^{-1}(p)$, so $\rho(g_1h_ng_2)K$ eventually lies in an
  arbitrary neighborhood of
  $\hat{\phi}^{-1}(g_1p) = \hat{\phi}^{-1}(z_+)$.

  When $k > 1$, consider the sequence
  \[
    \gamma_n' = g_2h_2^{(n)} \cdots g_kh_k^{(n)}g_{k+1}.
  \]
  Inductively we can assume that $\gamma_n' \to g_2p_2$,
  $(\gamma_n')^{-1}$ converges to $z_-$, and that for any compact
  $K \subset \hat{C}_{z_-}$, the set $\rho(\gamma_n')K$ lies in an
  arbitrarily small neighborhood of $\hat{\phi}^{-1}(g_2p_2)$. Then
  since $p_1 \ne g_2p_2$, for large enough $n$, since the strata
  $\hat{C}_z$ satisfy property \ref{item:weak_strata_interior}, the
  set $\rho(\gamma_n')K$ is eventually a compact subset of
  $\hat{C}_{p_1}$. Applying \Cref{lem:hats_weak_convergence} and
  equivariance again, we see that for large enough $n$,
  \[
    \rho(\gamma_n)K = \rho(g_1h_1^{(n)}) \rho(\gamma_n')K \subset U.
  \]

\end{proof}

The proposition below completes the proof of
\Cref{prop:weak_egf_to_egf_ii}, hence of \Cref{prop:weak_egf_to_egf}.

\begin{prop}
  The map $\hat{\phi}:\hat{\Lambda} \to \bgamh$ and sets $\hat{C}_z$
  satisfy condition \ref{item:strata_interior} in
  \Cref{defn:extends_convergence_dynamics}.
\end{prop}
\begin{proof}
  Fix $K \subset \bgamh$ compact. It suffices to show that for every
  $w \in \bgamh - K$, the set $\bigcap_{z \in K} \hat{C}_z$
  contains a neighborhood of $\hat{\phi}^{-1}(w)$. Suppose for a
  contradiction that this is not the case. We have already seen that
  the strata $\hat{C}_z$ satisfy condition
  \ref{item:weak_strata_interior}, so the only possibility is
  that there is a sequence of points $z_n \in K$ so that for every
  neighborhood $U$ of $\hat{\phi}^{-1}(w)$, the set
  $\hat{C}_n = \hat{C}_{z_n}$ eventually fails to contain $U$.

  Up to a subsequence, every $z_n$ is either conical or
  parabolic. Suppose first that every $z_n$ is conical, in which case
  $\hat{C}_n = \Opp(\hat{\phi}^{-1}(z_n))$ for every $n$. Thus, there
  is a sequence of flags $\xi_n \in \hat{\phi}^{-1}(z_n)$, and a
  sequence of flags $\eta_n \in \flags$, so that $\eta_n$ accumulates
  in $\hat{\phi}^{-1}(w)$, and $\xi_n, \eta_n$ are not opposite. After
  extraction, $\xi_n$ converges to $\xi \in \hat{\phi}^{-1}(K)$ and
  $\eta_n$ converges to $\eta \in \hat{\phi}^{-1}(w)$. The flags $\xi$
  and $\eta$ are not opposite, which contradicts antipodality of
  $\hat{\phi}$.

  So, now consider the case where each $z_n$ is parabolic. Let $Y$ be
  a Gromov model for $(\Gamma, \mc{H})$, and fix a bi-infinite
  geodesic $c_n:(-\infty, \infty) \to Y$ with its forward endpoint at
  $z_n$ and backward endpoint at $w$. This geodesic must eventually
  enter a horoball $\mc{B}$ in $Y$ centered at $z_n$, so let $y_n$ in
  $Y$ be the earliest point where this occurs. Since the action of
  $\Gamma$ on the complement of the horoballs is cocompact, there is
  some $\gamma_n \in \Gamma$ so that the distance from
  $\gamma_n^{-1}y_n$ to some fixed basepoint $y_0 \in Y$ is bounded,
  independent of $n$. Thus the horoball $\gamma_n^{-1}\mc{B}$ passes
  within uniform distance of $y_0$. There are only finitely many
  horoballs in $Y$ which do this, so up to subsequence
  $\gamma_n^{-1}z_n = p$ for a fixed parabolic point $p \in \bgamh$.

  We may assume that $z_n$ converges to some $z \in K$. Now, if the
  sequence $\gamma_n$ is bounded in $\Gamma$, then eventually
  $z_n = z$ for all $n$, and therefore an open neighborhood of
  $\hat{\phi}^{-1}(w)$ is contained in every $\hat{C}_n = \hat{C}_z$
  by condition \ref{item:weak_strata_interior}. So, assume that
  $\gamma_n$ tends to infinity in $\Gamma$. We claim that in this
  case, $\gamma_n$ converges to $z$. To see this, observe that since
  $z_n$ does not accumulate at $w$, the bi-infinite geodesic $c_n$
  passes within a uniform distance of a basepoint $y_0$ in $Y$. This
  means that the point $\gamma_n$ must actually lie uniformly close
  (in $Y$) to a point $r_n(t_n)$, where $r_n$ is a geodesic ray in $Y$
  from $y_0$ to $z_n$, and $t_n \to \infty$. Since $z_n \to z$ it
  follows that $r_n(t_n)$, hence $\gamma_n$, tends to $z$.

  Now, since our chosen strata satisfy condition
  \ref{item:weak_strata_interior}, and $w \ne z$, there is a
  neighborhood $U$ of $\hat{\phi}^{-1}(w)$ whose closure lies in
  $\hat{C}_z$. Extract a subsequence so that $\gamma_n^{-1} \to u$ for
  some $u \in \bgamh$. From \Cref{prop:weak_to_convergence}, the
  repelling strata satisfy condition \ref{item:extended_convergence},
  so $\rho(\gamma_n^{-1})U$ eventually lies in an arbitrarily small
  neighborhood of $\hat{\phi}^{-1}(u)$. However, we also know that
  $u \ne p$: there is a bi-infinite geodesic $\gamma_n^{-1}c_n$ in $Y$
  from $p$ to $\gamma_n^{-1}w$ passing through the identity, and since
  $w \ne z$, we have $\gamma_n^{-1}w \to u$ and therefore there is
  also a geodesic in $Y$ joining $p$ to $u$.

  By \ref{item:weak_strata_interior} again, the interior of
  $\hat{C}_p$ contains $\hat{\phi}^{-1}(u)$, and therefore eventually
  $\hat{C}_p$ contains $\rho(\gamma_n^{-1})U$. Then by equivariance,
  $\hat{C}_n = \rho(\gamma_n)\hat{C}_p$ contains $U$ for all
  sufficiently large $n$. This contradicts our original assumption.
\end{proof}

Note that by following the proof of \Cref{prop:weak_egf_to_egf}, we
see that one can always modify the repelling strata $C_z$ of an EGF
representation so that they are contained in the sets
$\Opp(\phi^{-1}(z))$ for the boundary extension $\phi$. We record this
fact in the proposition below.
\begin{prop}
  \label{prop:point_cell_transversality}
  Let $\rho:\Gamma \to G$ be an EGF representation. If
  $\phi:\Lambda \to \bgamh$ is a boundary extension for $\rho$, and
  $\{C_z\}_{z \in \bgamh}$ are a family of repelling strata, then so
  are the sets $\hat{C}_z = C_z \cap \Opp(\phi^{-1}(z))$.
\end{prop}
\begin{proof}
  By compactness of $\Lambda$, and the fact that antipodality is an
  open condition, the sets $\Opp(\phi^{-1}(z))$ satisfy condition
  \ref{item:strata_interior}. Since the sets $C_z$ also satisfy this
  condition by assumption, it is apparent that the intersections
  $C_z \cap \Opp(\phi^{-1}(z))$ satisfy it as well. Condition
  \ref{item:extended_convergence} is automatic from the fact that
  $\hat{C}_z \subseteq C_z$.
\end{proof}

\subsection{The ``conical/peripheral'' characterization of EGF
  representations}

It is often possible to prove properties of relatively hyperbolic
groups by first showing that the property holds for conical
subsequences in the group, and then showing that the property holds
inside of peripheral subgroups. There is a characterization of the EGF
property along these lines, which is frequently useful for
constructing examples of EGF representations (see
\cite{WeismanExamples}). As before, $(\Gamma, \mc{H})$ is a fixed
relatively hyperbolic pair, and $P \subset G$ is a symmetric parabolic
subgroup of $G$.

\begin{restatable}{prop}{conicalPeripheralEGF}
  \label{prop:conical_peripheral_implies_egf}
  Let $\rho:\Gamma \to G$ be a representation. Suppose that there is a
  continuous surjective $\rho$-equivariant antipodal map
  $\phi:\Lambda \to \bgamh$ and open subsets $C_z \subset \flags$ for
  each $z \in \bgamh$ satisfying property \ref{item:strata_interior}
  in \Cref{defn:extends_convergence_dynamics}.

  Then $\rho$ is an EGF representation if and only if the following
  conditions hold:
  \begin{enumerate}    
    \custitem[(P-conical)]\label{item:conical_extended_con} For any
    sequence $\gamma_n \in \Gamma$ limiting conically to some point in
    $\bgamh$, $\rho(\gamma_n^{\pm 1})$ is $P$-divergent and every
    $P$-limit point of $\rho(\gamma_n^{\pm 1})$ lies in $\Lambda$.

    \custitem[(C2-par)]\label{item:peripheral_extended_con}For any
    parabolic point $p \in \bgamh$, any compact $K \subset C_p$, and
    any open set $U$ containing $\phi^{-1}(p)$, for all but finitely
    many $\gamma \in \Gamma_p = \Stab_\Gamma(p)$, we have
    $\rho(\gamma)K \subset U$.
  \end{enumerate}
\end{restatable}

The proof of \Cref{prop:conical_peripheral_implies_egf} requires the
technical machinery of \emph{relative quasigeodesic automata}, so we
defer it to \Cref{sec:weak_egf}. At the end of \Cref{sec:weak_egf}, we
also prove another similar ``conical/peripheral'' characterization of
EGF representations which may be of interest.

\subsection{Properties of $\Lambda$}

\begin{prop}
  \label{prop:boundary_contains_limset}
  Let $(\Gamma, \mc{H})$ be a relatively hyperbolic pair, and let
  $\rho:\Gamma \to G$ be a representation which is EGF with respect to
  a symmetric parabolic $P$, with boundary extension
  $\phi:\Lambda \to \bgamh$. Then $\Lambda$ contains the $P$-limit set
  of $\rho(\Gamma)$.
\end{prop}
\begin{proof}
  Let $\xi \in \flags$ be a flag in the $P$-limit set of
  $\rho(\Gamma)$. Then there is a $P$-contracting sequence
  $\rho(\gamma_n)$ for $\gamma_n \in \Gamma$ and a flag
  $\xi_- \in G/P$ such that $\rho(\gamma_n)\eta$ converges to $\xi$
  for any $\eta$ in $\Opp(\xi_-)$. Up to subsequence
  $\gamma_n^{\pm 1}$ converges to $z_\pm \in \bgamh$, so for any flag
  $\eta \in C_{z_-}$, the sequence $\rho(\gamma_n)\eta$ subconverges
  to a point in $\phi^{-1}(z_+)$. But since $\Opp(\xi_-)$ is open and
  dense, for some $\eta \in C_{z_-}$ we have
  $\rho(\gamma_n)\eta \to \xi$ and hence $\xi \in \phi^{-1}(z_+)$.
\end{proof}

In particular, \Cref{prop:boundary_contains_limset} implies that the
EGF boundary set $\Lambda \subset \flags$ of an EGF representation
$\rho:\Gamma \to G$ must always contain the \emph{$P$-proximal limit
  set} of $\rho(\Gamma)$. (Recall that $g \in G$ is
\emph{$P$-proximal} if it has a unique attracting fixed point in
$\flags$; the \emph{$P$-proximal limit set} of a subgroup of $G$ is
the closure of the set of attracting fixed points of $P$-proximal
elements).

We will see that most of the power of EGF representations lies in the
fact that their associated boundary extensions
$\phi:\Lambda \to \bgamh$ do not have to be homeomorphisms (so
the Bowditch boundary of $\Gamma$ does not need to be equivariantly
embedded in any flag manifold). However, it turns out that it is
always possible to choose the boundary extension $\phi$ so that it has
a well-defined inverse on \emph{conical limit points} in
$\bgamh$. In fact, we can even get a somewhat precise description
of all the fibers of $\phi$. Concretely, we have the following:
\begin{prop}
  \label{prop:conical_pts_singletons}
  Let $\rho:\Gamma \to G$ be an EGF representation, with boundary
  extension $\phi:\Lambda \to \bgamh$. There is a
  $\rho(\Gamma)$-invariant closed subset $\Lambda' \subset \flags$ and
  a $\rho$-equivariant map $\phi':\Lambda' \to \bgamh$
  such that:
  \begin{enumerate}
  \item $\phi':\Lambda' \to \bgamh$ is also a boundary
    extension for $\rho$,
  \item for every $z \in \conbdry(\Gamma, \mc{H})$, $\phi'^{-1}(z)$ is
    a singleton, and
  \item for every $p \in \parbdry(\Gamma, \mc{H})$, $\phi'^{-1}(p)$ is
    the closure of the set of all accumulation points of orbits
    $\rho(\gamma_n) \cdot x$ for $\gamma_n$ a sequence of pairwise
    distinct elements in $\Gamma_p$ and $x \in C_p$.
  \end{enumerate}
\end{prop}

We will prove \Cref{prop:conical_pts_singletons} at the end of
\Cref{sec:relative_stability}, where it will follow as a consequence
of the proof of the relative stability theorem for EGF representations
(\Cref{thm:cusp_stable_stability})---see
\Cref{rem:conical_pts_singletons}.

We will rely on both \Cref{prop:conical_peripheral_implies_egf} and
\Cref{prop:conical_pts_singletons} to prove some of the remaining
results in this section, which are not needed anywhere else in this
paper.

\subsection{Relatively Anosov representations}

EGF representations give a strict generalization of the relative
Anosov representations mentioned in the introduction. We give a
precise definition here.

\begin{definition}[{\cite[Definition 7.1]{kl2018relativizing} or
    \cite[Definition 1.1]{zz1}; see also \cite[Proposition
    4.4]{zz1}}]
  \label{defn:relative_anosov}
  Let $\Gamma$ be a subgroup of $G$ and suppose $(\Gamma, \mc{H})$ is
  a relatively hyperbolic pair. Let $P \subset G$ be a symmetric
  parabolic subgroup.

  The subgroup $\Gamma$ is \emph{relatively $P$-Anosov} if it is
  $P$-divergent, and there is a $\Gamma$-equivariant antipodal
  embedding $\bgamh \to \flags$ whose image $\Lambda$ is the $P$-limit
  set of $\Gamma$.
\end{definition}

Here, we say an embedding $\psi:\bgamh \to \flags$ is \emph{antipodal}
if for every distinct $\xi_1, \xi_2$ in $\bgamh$, $\psi(\xi_1)$ and
$\psi(\xi_2)$ are opposite flags. 

\begin{remark} \
  \begin{enumerate}[(a)]
  \item When $\Gamma$ is a \emph{hyperbolic} group (and the collection
    of peripheral subgroups $\mc{H}$ is empty), then the Bowditch
    boundary $\bgamh$ is identified with the Gromov boundary
    $\dee \Gamma$. In this case, \Cref{defn:relative_anosov} coincides
    with the usual definition of an Anosov representation.

  \item In general, it is possible to define (relatively) $P$-Anosov
    representations for a \emph{non-symmetric} parabolic subgroup
    $P$. However, there is no loss of generality in assuming that $P$
    is symmetric: a representation $\rho:\Gamma \to G$ is $P$-Anosov
    if and only if it is $P'$-Anosov for a symmetric parabolic
    subgroup $P' \subset G$ depending only on $P$.
  \end{enumerate}
\end{remark}

In their paper \cite{kl2018relativizing}, Kapovich--Leeb gave several
possible definitions for a ``relatively Anosov representation.''
\Cref{defn:relative_anosov} is equivalent to (essentially) their
weakest definition, that of a \emph{relatively asymptotically
  embedded} representation; their stronger definitions are those of a
\emph{relatively Morse} representation and a \emph{uniformly regular
  relatively asymptotically embedded representation}. In \cite{zz1},
Zhu--Zimmer showed that when $P = P_k$ is the stabilizer of a
$(k, d - k)$-flag in $\SL(d, \K)$ ($\K = \R$ or $\C$), then relatively
asymptotically embedded representations are automatically uniformly
regular, but the following example shows that this is not true in
general.
\begin{example}[See also {\cite[Example 5.6]{kl2018relativizing}}]
  \label{ex:nonuniform_example}
  Let $S$ be a finite-type non-compact orientable surface with
  negative Euler characteristic, and equip $S$ with a pair of
  hyperbolic structures $X_1, X_2$ so that in $X_1$, an end of $S$ is
  a cusp, and in $X_2$, the same end of $S$ is an infinite-area
  ``funnel.'' For $i = 1,2$, let $\rho_i:\pi_1S \to \PSL(2, \R)$ be
  the holonomy representation for $X_i$, and let
  \[
    \rho:\pi_1S \to \PSL(2, \R) \times \PSL(2, \R)
  \]
  be the product $\rho = (\rho_1, \rho_2)$.
  
  If $P \subset \PSL(2, \R)$ is the stabilizer of a point in
  $\dee \H^2$, then $P \times \PSL(2, \R)$ is a symmetric parabolic
  subgroup in $\PSL(2, \R) \times \PSL(2, \R)$. Then, letting $\mc{H}$
  be the collection of cusp subgroups of $\pi_1S$ for the hyperbolic
  structure $X_1$, the representation $\rho$ is relatively
  $(P \times \PSL(2, \R))$-Anosov relative to $\mc{H}$ in the sense of
  \Cref{defn:relative_anosov}. However, this representation is not
  \emph{uniformly $\taumod$-regular} in the sense of
  Kapovich--Leeb--Porti \cite{klp2017anosov} if $\taumod$ is the
  ``model face'' corresponding to $P \times \PSL(2, \R)$.

  Note that this example also shows that in general semisimple Lie
  groups, the peripheral subgroups for a relatively Anosov
  representation need not have weakly unipotent image.
\end{example}

The proposition below (a restatement of
\Cref{thm:bdry_dynamics_homeo_equiv_rel_asymp_embedded}) gives the
relationship between relatively Anosov representations and EGF
representations.
\begin{prop}
  \label{prop:rel_anosov_egf}
  Let $(\Gamma, \mc{H})$ be a relatively hyperbolic pair, and let
  $P \subset G$ be a symmetric parabolic subgroup. A representation
  $\rho:\Gamma \to G$ is relatively $P$-Anosov in the sense of
  \Cref{defn:relative_anosov} if and only if $\rho$ is EGF with
  respect to $P$, and has an injective boundary extension
  $\phi:\Lambda \to \bgamh$.
\end{prop}

As a first step towards the proof, we show:
\begin{prop}
  \label{prop:homeomorphism_implies_asymptotically_embedded}
  Let $\rho:\Gamma \to G$ be an EGF representation with respect to
  $P$, and suppose that the boundary extension
  $\phi:\Lambda \to \bgamh$ is a homeomorphism. Then:
  \begin{enumerate}
  \item $\rho(\Gamma)$ is $P$-divergent, and $\Lambda$ is the
    $P$-limit set of $\rho(\Gamma)$.
    
  \item The sets $C_z$ for $z \in \bgamh$ can be taken to be
    \[
      \Opp(\phi^{-1}(z)) = \{\nu \in \flags : \nu \textrm{ is opposite to
      } \phi^{-1}(z)\}.
    \]
  \end{enumerate}
\end{prop}
\begin{proof}
  (1). Let $\gamma_n$ be any infinite sequence of elements in
  $\Gamma$. After extracting a subsequence, we have
  $\gamma_n^{\pm1} \to z_\pm$, and since $\phi$ is a homeomorphism,
  $\rho(\gamma_n)$ converges to the point $\phi^{-1}(z_+)$ uniformly
  on compacts in the open set $C_{z_-}$. Then
  \Cref{prop:local_contracting_implies_global_contracting} implies
  that $\rho(\gamma_n)$ is $P$-divergent, with unique $P$-limit point
  $\phi^{-1}(z_+) \in \Lambda$.

  (2). Condition \ref{item:strata_interior} follows directly from
  antipodality of $\phi$ and compactness of $\Lambda$, so we just need
  to see that these sets satisfy condition
  \ref{item:extended_convergence}. Let $\gamma_n$ be an infinite
  sequence in $\Gamma$ with $\gamma_n^{\pm 1} \to z_{\pm}$ for
  $z_\pm \in \bgamh$.

  We know that for open subsets $U_\pm \subset \flags$, we have
  $\rho(\gamma_n) \cdot U_+ \to \phi^{-1}(z_+)$ and
  $\rho(\gamma_n^{-1}) U_- \to \phi^{-1}(z_-)$, uniformly on
  compacts. \Cref{prop:local_contracting_implies_global_contracting}
  implies that $\rho(\gamma_n)$ and $\rho(\gamma_n^{-1})$ are both
  $P$-divergent with unique $P$-limit points $\phi^{-1}(z_+)$,
  $\phi^{-1}(z_-)$. So in fact by \Cref{lem:klp_contraction_symmetry}
  $\rho(\gamma_n)$ converges to $\phi^{-1}(z_+)$ uniformly on compacts
  in $\Opp(\phi^{-1}(z_-))$.
\end{proof}

\begin{proof}[Proof of \Cref{prop:rel_anosov_egf}]
  \Cref{prop:homeomorphism_implies_asymptotically_embedded} ensures
  that if $\rho$ is an EGF representation, and the boundary extension
  $\phi$ is a homeomorphism, then $\rho$ is $P$-divergent and
  $\phi^{-1}$ is an antipodal embedding whose image is the $P$-limit
  set.

  On the other hand, if $\rho$ is relatively $P$-Anosov, with boundary
  embedding $\psi:\bgamh \to \Lambda$, for each $z \in \bgamh$, we can
  take
  \[
    C_z = \Opp(\psi(z)).
  \]
  Antipodality means that $C_z$ contains an open neighborhood of
  $\Lambda - \psi(z)$, and $P$-divergence and
  \Cref{lem:klp_contraction_symmetry} imply that $\rho(\Gamma)$ has
  the appropriate convergence dynamics.
\end{proof}

\subsection{Relativization}
\label{sec:relativization}

We now turn to the situation where we have an EGF representation of a
\emph{hyperbolic} group $\Gamma$ with a \emph{nonempty} collection of
peripheral subgroups. That is, for some invariant set
$\Lambda \subset \flags$, we have an EGF boundary extension
$\phi:\Lambda \to \bgamh$, where $\bgamh$
is the Bowditch boundary of $\Gamma$ with peripheral structure
$\mc{H}$.

We want to prove \Cref{thm:changing_peripheral_structure}, which says
that in this situation, if $\rho$ restricts to a $P$-Anosov
representation on each $H \in \mc{H}$, then $\rho$ is a $P$-Anosov
representation of $\Gamma$. For the rest of this section, we assume
that $\Gamma$ is a hyperbolic group, and $\mc{H}$ is a collection of
subgroups of $\Gamma$ so that the pair $(\Gamma, \mc{H})$ is
relatively hyperbolic. We let $\rho:\Gamma \to G$ be an EGF
representation for the pair $(\Gamma, \mc{H})$ with respect to a
symmetric parabolic subgroup $P \subset G$, and we assume that for
each $H \in \mc{H}$, $\rho|_H:H \to G$ is $P$-Anosov, with Anosov
limit map $\psi_H:\dee H \to \flags$.

The main step in the proof is to observe that it is always possible to
choose the boundary extension $\phi:\Lambda \to \bgamh$
so that $\Lambda$ is equivariantly homeomorphic to the Gromov boundary
of $\Gamma$ (which we here denote $\dee \Gamma$).

Whenever $\Gamma$ is a hyperbolic group and $\mc{H}$ is a collection
of subgroups so that $(\Gamma, \mc{H})$ is a relatively hyperbolic
pair, there is an explicit description of the Bowditch boundary
$\bgamh$ in terms of the Gromov boundary $\dee \Gamma$
of $\Gamma$---see \cite{gerasimov2012floyd}, \cite{gp2013quasi}, or
\cite{tran2013relations}. Specifically, we can say:
\begin{prop}
  \label{prop:bb_quotient_description}
  There is an equivariant surjective continuous map
  $\phi_\Gamma:\dee \Gamma \to \bgamh$ such that for each conical
  limit point $z$ in $\bgamh$, $\phi_\Gamma^{-1}(z)$ is a singleton,
  and for each parabolic point $p \in \bgamh$ with
  $H = \Stab_\Gamma(p)$, $\phi_\Gamma^{-1}(p)$ is an embedded copy of
  $\dee H$ in $\dee \Gamma$.
\end{prop}

In our situation, we can see that the boundary extension
$\phi:\Lambda \to \bgamh$ satisfies similar properties.




\begin{lem}
  \label{lem:parabolic_preimages_anosov_images}
  There is a closed $\rho(\Gamma)$-invariant subset
  $\Lambda' \subset \flags$ and an EGF boundary extension
  $\phi':\Lambda' \to \bgamh$ such that:
  \begin{enumerate}
  \item For each conical limit point $z \in \bgamh$,
    $\phi'^{-1}(z)$ is a singleton.
  \item For each parabolic point $p \in \bgamh$, with
    $H = \Stab_\Gamma(p)$, we have $\phi'^{-1}(p) = \psi_H(\dee H)$.
  \end{enumerate}
\end{lem}

\begin{proof}
  We choose $\Lambda'$ as in \Cref{prop:conical_pts_singletons}. The
  only thing we need to check is that for $H = \Stab_\Gamma(p)$, the
  set $\psi_H(\dee H)$ is exactly the closure of the set of
  accumulation points of $\rho(H)$-orbits in $C_p$. But since we may
  assume $C_p$ is contained in $\Opp(\psi_H(\dee H))$, this follows
  immediately from the fact that $\rho(H)$ is $P$-divergent and the
  closed set $\psi_H(\dee H)$ is the $P$-limit set of $\rho(H)$.   



\end{proof}

Next we need a lemma which will allow us to characterize the Gromov
boundary of $\Gamma$ as an extension of the Bowditch boundary
$\bgamh$. First recall that if $\Gamma$ acts as a
convergence group on a space $Z$, the \emph{limit set} of $\Gamma$ is
the set of points $z \in Z$ such that for some $y \in Z$ and some
sequence $\gamma_n \in \Gamma$, we have
\[
  \gamma_n |_{Z - \{y\}} \to z
\]
uniformly on compacts.

\begin{lem}
  \label{lem:sequence_convergence_agrees}
  Let $\Gamma$ act on compact metrizable spaces $X$ and $Y$, and let
  $\phi_X:X \to \bgamh$, $\phi_Y:Y \to \bgamh$ be continuous
  equivariant surjective maps such that for every conical limit point
  $z \in \bgamh$, $\phi_X^{-1}(z)$ and $\phi_Y^{-1}(z)$ are both
  singletons, and for every parabolic point $p \in \bgamh$,
  $H = \Stab_\Gamma(p)$ acts as a convergence group on $X$ and $Y$,
  with limit sets $\phi_X^{-1}(p)$, $\phi_Y^{-1}(p)$ equivariantly
  homeomorphic to $\dee H$.

  Then for any sequences $z_n, z_n' \in \conbdry(\Gamma, \mc{H})$, we
  have
  \[
    \lim_{n \to \infty} \phi_X^{-1}(z_n) = \lim_{n \to \infty}\phi_X^{-1}(z_n')
  \]
  if and only if
  \[
    \lim_{n \to \infty} \phi_Y^{-1}(z_n) = \lim_{n \to
      \infty}\phi_Y^{-1}(z_n').
  \]
\end{lem}
\begin{proof}
  We proceed by contradiction, and suppose that for a pair of
  sequences $z_n, z_n' \in \conbdry(\Gamma, \mc{H})$, we have
  \[
    \lim_{n \to \infty} \phi_X^{-1}(z_n) = \lim_{n \to
      \infty}\phi_X^{-1}(z_n') = x,
  \]
  but
  \[
    \lim_{n \to \infty}\phi_Y^{-1}(z_n) \ne \lim_{n \to
    \infty}\phi_Y^{-1}(z_n').
  \]
  After taking a subsequence we may assume $z_n$ converges to
  $z \in \bgamh$, and that $y_n = \phi_Y^{-1}(z_n)$
  converges to $y$ and $y_n' = \phi_Y^{-1}(z_n')$ converges to $y'$
  for $y \ne y'$. By continuity, we have
  \[
    \phi_Y(y) = \phi_Y(y') = \phi_X(x) = z.
  \]

  Since $\phi_X$ and $\phi_Y$ are bijective on
  $\phi_X^{-1}(\conbdry(\Gamma, \mc{H}))$ and
  $\phi_Y^{-1}(\conbdry(\Gamma, \mc{H}))$ respectively, we must have
  $z = p$ for a parabolic point $p \in \parbdry(\Gamma, \mc{H})$. Let
  $H = \Stab_\Gamma(p)$.

  Since $p$ is a bounded parabolic point, we can find sequences of
  group elements $h_n, h_n' \in H$ so that for a fixed compact subset
  $K \subset \bgamh - \{p\}$, we have
  \begin{equation}
    \label{eq:bounded_parabolic}
    h_nz_n \in K, \quad h_n'z_n' \in K. 
  \end{equation}
  This implies that no subsequence of $h_n y_n$ or $h_n'y_n'$
  converges to a point in $\phi_Y^{-1}(p)$.

  Then, since $H$ acts as a convergence group on $Y$ with limit set
  $\phi_Y^{-1}(p)$, up to subsequence there are points
  $u, u' \in \phi_Y^{-1}(p)$ so that $h_n$ converges to a point in
  $\phi_Y^{-1}(p)$ uniformly on compacts in $Y - \{u\}$, and $h_n'$
  converges to a point in $\phi_Y^{-1}(p)$ uniformly on compacts in
  $Y - \{u'\}$. So, we must have $u = y$ and $u' = y'$.

  This means that the sequences $h_n^{-1}$ and $h_n'^{-1}$ have
  distinct limits in the compactification
  $\overline{H} = H \sqcup \dee H$. So, there are distinct points
  $v, v' \in \phi_X^{-1}(p)$ so that (again up to subsequence) $h_n$
  converges to a point in $\phi_X^{-1}(p)$ uniformly on compacts in
  $X - \{v\}$, and $h_n'$ converges to a point in $\phi_X^{-1}(p)$
  uniformly on compacts in $X - \{v'\}$. Without loss of generality,
  we can assume $x \ne v$.

  But then $\phi_X^{-1}(z_n)$ lies in a compact subset of $X - \{v\}$,
  so $h_n\phi_X^{-1}(z_n)$ converges to a point in $\phi_X^{-1}(p)$
  and $h_nz_n$ converges to $p$. But this contradicts
  (\ref{eq:bounded_parabolic}) above.
\end{proof}

\begin{prop}
  \label{prop:gromov_boundary_homeo}
  If the set $\Lambda$ satisfies the conclusions of
  \Cref{lem:parabolic_preimages_anosov_images}, then $\Lambda$ is
  equivariantly homeomorphic to the Gromov boundary of $\Gamma$.
\end{prop}
\begin{proof}  
  Let $\phi_\Gamma : \dee \Gamma \to \bgamh$ denote the
  quotient map identifying the limit set of each $H \in \mc{H}$ to the
  parabolic point $p$ with $H = \Stab_\Gamma(p)$. For each conical
  limit point $z \in \bgamh$, the fiber
  $\phi_\Gamma^{-1}(z)$ is a singleton. So, there is an equivariant
  bijection $f$ from $\phi_\Gamma^{-1}(\conbdry(\Gamma, \mc{H}))$ to
  $\phi^{-1}(\conbdry(\Gamma, \mc{H}))$.

  Moreover, since $\phi_\Gamma^{-1}(\conbdry(\Gamma, \mc{H}))$ is
  $\Gamma$-invariant, and the action of $\Gamma$ on its Gromov
  boundary $\dee \Gamma$ is minimal,
  $\phi_\Gamma^{-1}(\conbdry(\Gamma, \mc{H}))$ is dense in
  $\dee \Gamma$. We claim that $f$ extends to a continuous injective
  map $\dee \Gamma \to \Lambda$ by defining $f(x) = \lim f(x_n)$ for
  any sequence $x_n \to x$.

  To see this, we can apply \Cref{lem:sequence_convergence_agrees},
  taking $\dee \Gamma = X$ and $\Lambda = Y$. We know that $\Gamma$
  always acts on its own Gromov boundary as a convergence group (so in
  particular each $H \in \mc{H}$ acts on $\dee \Gamma$ as a
  convergence group with limit set $\dee H$). And, since $\rho$
  restricts to a $P$-Anosov representation on each $H \in \mc{H}$, for
  any infinite sequence $h_n \in H$, up to subsequence there are
  $u, u_- \in \psi_H(\dee H)$ so that $\rho(h_n)$ converges to $u$
  uniformly on compacts in $\Opp(u_-)$. Antipodality of $\phi$ implies
  that $\rho(h_n)$ converges to $u$ uniformly on compacts in
  $\Lambda - \psi_H(\dee H)$. The other hypotheses of
  \Cref{lem:sequence_convergence_agrees} follow from
  \Cref{prop:bb_quotient_description} and
  \Cref{lem:parabolic_preimages_anosov_images}.

  We still need to check that $f$ is actually surjective. We know that
  $f$ restricts to a bijection on
  $\phi_\Gamma^{-1}(\conbdry(\Gamma, \mc{H}))$, and that $f$ takes
  $\phi_\Gamma^{-1}(p)$ to $\phi^{-1}(p)$ for each parabolic point $p$
  in $\bgamh$. So we just need to check that for every $H \in \mc{H}$,
  $f$ restricts to a surjective map $\dee H \to \psi_H(\dee H)$. If
  $H$ is non-elementary, this must be the case because the action of
  $H$ on $\psi_H(\dee H)$ is minimal and $f$ maps $\dee H$ into
  $\psi_H(\dee H)$ as an invariant closed subset. Otherwise, $H$ is
  virtually cyclic and $\dee H$, $\psi_H(\dee H)$ both contain exactly
  two points. Then injectivity of $f$ implies surjectivity.

  So we conclude that there is a continuous bijection
  $f:\dee \Gamma \to \Lambda$, and since $\dee \Gamma$ is compact and
  $\Lambda$ is metrizable, $f$ is a homeomorphism.
\end{proof}

We let $f:\dee \Gamma \to \Lambda$ denote the equivariant
homeomorphism from \Cref{prop:gromov_boundary_homeo}.
\begin{prop}
  \label{prop:egf_anosov_extends_convergence}
  The equivariant homeomorphism $f^{-1}:\Lambda \to \dee \Gamma$ is
  antipodal. Moreover, together with the strata
  $C_z = \Opp(f^{-1}(z))$ for all $z \in \dee \Gamma$, it satisfies
  conditions \ref{item:strata_interior},
  \ref{item:conical_extended_con}, and
  \ref{item:parabolic_flag_convergence} in
  \Cref{prop:conical_peripheral_implies_egf}.
\end{prop}
\begin{proof}
  First we show antipodality. Fix distinct points
  $z, z' \in \dee \Gamma$. If $\phi_\Gamma(z), \phi_\Gamma(z')$ are
  distinct, then $f(z), f(z')$ are opposite by antipodality of
  $\phi$. Otherwise, $z, z'$ correspond to distinct points in $\dee H$
  for some $H \in \mc{H}$, and $f(z), f(z')$ are opposite by
  antipodality of $\psi_H$.
  
  Now we turn to the remaining conditions. Since $\Gamma$ is
  hyperbolic, condition \ref{item:parabolic_flag_convergence} is
  vacuous, and condition \ref{item:strata_interior} follows from
  antipodality of $f$, compactness of $\Lambda$, and the fact that
  antipodality is an open condition in $G/P$. So we just need to
  consider condition \ref{item:conical_extended_con}. That is, we need
  to show that if $\gamma_n \in \Gamma$ is a conical limit sequence
  with $\gamma_n^{\pm 1} \to z_{\pm}$ for $z_\pm \in \dee \Gamma$,
  then every $P$-limit point of $\rho(\gamma_n^{\pm 1})$ lies in
  $\Lambda$.

  We consider two cases:
  \begin{aside}{Case 1: $\phi \circ f(z_+)$ is a parabolic point $p$
      in $\bgamh$}
    In this case, $\gamma_n$ lies along a quasigeodesic ray in
    $\Gamma$ limiting to some $z_+ \in \dee H$, with
    $H = \Stab_\Gamma(p)$. This means that for a bounded sequence
    $b_n \in \Gamma$, we have $\gamma_nb_n \in H$. Since $\rho$
    restricts to a $P$-Anosov representation on $H$, this means that
    $\rho(\gamma_nb_n)$ is $P$-divergent and every $P$-limit point of
    $\rho(\gamma_nb_n)$ lies in $\psi_H(z_+)$. For the same reason,
    every $P$-limit point of $\rho(b_n^{-1}\gamma_n^{-1})$ lies in
    $\psi_H(z_+)$.

    Up to subsequence $b_n$ is a constant $b$. We can use
    \Cref{prop:divergence_equivalent_contracting} to see that
    $\rho(\gamma_n)$ has the same $P$-limit set as
    $\rho(\gamma_nb)$. This $P$-limit set lies in $\Lambda$. And,
    every $P$-limit point of $\rho(\gamma_n^{-1})$ is a $b$-translate
    of a $P$-limit point of $\rho(b^{-1}\gamma_n^{-1})$. This
    $P$-limit set also lies in $\Lambda$.
  \end{aside}
  \begin{aside}{Case 2: $\phi \circ f(z_+)$ is a conical limit point
      in $\bgamh$}
    In this case, a subsequence of $\gamma_n$ is a conical limit
    sequence for the action of $\Gamma$ on $\bgamh$, and the desired
    result follows from the ``only if'' part of
    \Cref{prop:conical_peripheral_implies_egf}.
  \end{aside}
\end{proof}

\begin{proof}[Proof of \Cref{thm:changing_peripheral_structure}]
  Let $\Gamma$ be hyperbolic, let $\mc{H}$ be a collection of
  subgroups such that $(\Gamma, \mc{H})$ is a relatively hyperbolic
  pair, and let $\rho:\Gamma \to G$ be an EGF representation with
  respect to $P$, for the peripheral structure $\mc{H}$.

  Suppose that $\rho$ restricts to a $P$-Anosov representation on each
  $H \in \mc{H}$. Propositions \ref{prop:boundary_contains_limset} and
  \ref{prop:egf_anosov_extends_convergence} imply that $\rho$ is
  \emph{also} an EGF representation of $\Gamma$ for its \emph{empty}
  peripheral structure, whose boundary extension can be chosen to be a
  homeomorphism. Then
  \Cref{thm:bdry_dynamics_homeo_equiv_rel_asymp_embedded} says that
  $\rho$ is relatively $P$-Anosov (again for the empty peripheral
  structure on $\Gamma$). This ensures that $\rho$ is actually
  (non-relatively) $P$-Anosov; see e.g. \cite[Theorem
  1.1]{klp2017anosov}.
\end{proof}


\section{Relative quasigeodesic automata}
\label{sec:relative_automata}
In the next three sections, we develop the technical tools needed to
prove the main results of the paper: namely, a \emph{relative
  quasigeodesic automaton} for a relatively hyperbolic group $\Gamma$
acting on a flag manifold $\flags$, and a \emph{system of open sets}
in $\flags$ which is in some sense \emph{compatible} with both the
relative quasigeodesic automaton and the action of $\Gamma$ on
$\flags$.

The basic idea is motivated by the computational theory of hyperbolic
groups. Given a hyperbolic group $\Gamma$ with finite generating set
$S$, it is always possible to find a finite directed graph $\mc{G}$,
with edges labeled by elements of $S$, so that directed paths on
$\mc{G}$ starting at a fixed vertex $v_{\mr{id}} \in \mc{G}$ are in
one-to-one correspondence with geodesic words in $\Gamma$. The graph
$\mc{G}$ is called a \emph{geodesic automaton} for $\Gamma$.

Geodesic automata are really a manifestation of the local-to-global
principle for geodesics in hyperbolic metric spaces: the fact that the
automaton exists means that it is possible to recognize a geodesic
path in a hyperbolic group just by looking at bounded-length subpaths.

In this section of the paper, we consider a \emph{relative} version of
a geodesic automaton. This is a finite directed graph $\mc{G}$ which
encodes the behavior of quasigeodesics in the coned-off Cayley graph
of a relatively hyperbolic group $\Gamma$. Eventually, our goal is to
build such an automaton by looking at the dynamics of the action of
$\Gamma$ on its Bowditch boundary $\bgamh$. The main result of this section is
\Cref{prop:relative_geodesic_automaton_from_bb}, which says that we
can construct such a \emph{relative quasigeodesic automaton} for a
relatively hyperbolic pair $(\Gamma, \mc{H})$ using an open covering
of the Bowditch boundary $\bgamh$ which satisfies
certain technical conditions.

In this section of the paper and the next, we will work in the general
context of a relatively hyperbolic group $\Gamma$ acting by
homeomorphisms on a connected compact metrizable space $M$, before
returning to the case where $M$ is a flag manifold $\flags$ for the
rest of the paper.

Throughout the rest of this section, we fix a \emph{non-elementary}
relatively hyperbolic pair $(\Gamma, \mc{H})$, and let
$\Pi \subset \parbdry(\Gamma, \mc{H})$ be a finite set, containing
exactly one point from each $\Gamma$-orbit in
$\parbdry(\Gamma, \mc{H})$. We also fix a finite generating set $S$
for $\Gamma$, which allows us to refer to the coned-off Cayley graph
$\conedgraph$ (\Cref{defn:conedgraph}).

\begin{definition}
  A \emph{$\Gamma$-graph} is a finite directed graph $\mc{G}$ where
  each vertex $v$ is labelled with a subset $T_v \subset \Gamma$,
  which is either:
  \begin{itemize}
  \item A singleton $\{\gamma\}$, with $\gamma \ne \mr{id}$, or
  \item A cofinite subset of a coset $g\Gamma_p$ for some $p \in \Pi$,
    $g \in \Gamma$.
  \end{itemize}

  A sequence $\{\alpha_n\} \subset \Gamma$ is a \emph{$\mc{G}$-path}
  if $\alpha_n \in T_{v_n}$ for a vertex path $\{v_n\}$ in
  $\mc{G}$.
\end{definition}

\begin{remark}
  We will often refer to ``the'' vertex path $\{v_n\}$ corresponding
  to a $\mc{G}$-path $\{\alpha_n\}$, although we will never actually
  verify that such a vertex path is uniquely determined by the
  sequence of group elements $\{\alpha_n\}$ in $\Gamma$.
\end{remark}

A vertex of a $\Gamma$-graph which is labeled by a cofinite subset of
a (necessarily unique) coset $g\Gamma_p$ is a \emph{parabolic
  vertex}. If $v$ is a parabolic vertex, we let $p_v = g \cdot p$
denote the corresponding parabolic point in
$\parbdry(\Gamma, \mc{H})$.

\begin{remark}
  It will be convenient to allow parabolic vertices to be labeled by
  \emph{cofinite} subsets of peripheral cosets (instead of just the
  entire coset) when we construct $\Gamma$-graphs using the
  convergence dynamics of the $\Gamma$-action on $\bgamh$.
\end{remark}

\begin{definition}
  Let $z \in \bgamh$. We say that a $\mc{G}$-path $\{\alpha_n\}$
  \emph{limits to $z$} if either:
  \begin{itemize}
  \item $z \in \conbdry(\Gamma, \mc{H})$, $\{\alpha_n\}$ is infinite,
    and the sequence
    \[
      \{\gamma_n = \alpha_1 \cdots \alpha_n\}_{n=1}^\infty
    \]
    limits to $z$ in the compactification
    $\overline{\Gamma} = \Gamma \sqcup \bgamh$, or
  \item $z \in \parbdry(\Gamma, \mc{H})$, $\{\alpha_n\}$ is a finite
    $\mc{G}$-path whose corresponding vertex path $\{v_n\}$ ends at a
    parabolic vertex $v_N$, and
    \[
      z = \alpha_1 \cdots \alpha_{N-1} p_{v_N}.
    \]
  \end{itemize}
\end{definition}

\begin{definition}
  Let $\mc{G}$ be a $\Gamma$-graph. The \emph{endpoint} of a finite
  $\mc{G}$-path $\{\alpha_n\}_{n=1}^N$ is
  \[
    \alpha_1 \cdots \alpha_N.
  \]
\end{definition}

\begin{definition}
  \label{defn:rel_automaton}
  A $\Gamma$-graph $\mc{G}$ is a \emph{relative quasigeodesic
    automaton} if:
  \begin{enumerate}
  \item\label{item:rel_automaton_geodesic} There is a constant $D > 0$
    and a Gromov model $Y$ for $(\Gamma, \mc{H})$ with basepoint $y_0$
    so that for any infinite $\mc{G}$-path $\{\alpha_n\}$, the
    sequence
    \[
      \alpha_1 \cdots \alpha_ny_0
    \]
    lies within a $D$-neighborhood of a geodesic ray $c$ in $Y$.
  \item For every $z \in \bgamh$, there exists a $\mc{G}$-path
    limiting to $z$.
  \end{enumerate}
\end{definition}

One way to think of a relative quasigeodesic automaton is that it
gives us a system for finding quasigeodesic representatives of every
element in the group. More concretely, we have the following:

\begin{lem}
  \label{lem:relative_automaton_quasidense_image}
  Let $\mc{G}$ be a relative quasigeodesic automaton. There is a
  constant $R > 0$ so that set of endpoints of $\mc{G}$-paths is
  $R$-dense in $\Gamma$.
\end{lem}
\begin{proof}
  If $\Gamma$ is hyperbolic and $\mc{H}$ is empty, then this is a
  consequence of the Morse lemma and the fact that the union of the
  images of all infinite geodesic rays based at the identity in
  $\Gamma$ is coarsely dense in $\Gamma$ (see
  \cite{bogopolski97infinite}).

  If $\mc{H}$ is nonempty, there is some $R > 0$ so that the union of
  all of the cosets $g \cdot \Gamma_p$ for $p \in \Pi$ is $R$-dense in
  $\Gamma$. So it suffices to show that for each $p \in \Pi$, there is
  some $R > 0$ so that all but $R$ elements in any coset
  $g \cdot \Gamma_p$ are the endpoints of a $\mc{G}$-path.

  For any such coset $g \cdot \Gamma_p$, we can find a finite
  $\mc{G}$-path $\{\alpha_n\}_{n=1}^{N-1}$ limiting to the vertex
  $g \cdot p$. That is,
  \[
    g \cdot p = \alpha_1 \cdots \alpha_{N-1} p_{v_N}.
  \]
  By definition $p_{v_N} = g' \cdot p$ with $T_{v_N}$ a cofinite
  subset of the coset $g' \Gamma_p$ That is,
  \[
    g \cdot \Gamma_p = \alpha_1 \cdots \alpha_{N-1} g' \Gamma_p,
  \]
  so for all but finitely many $\gamma \in g \cdot \Gamma_p$
  (depending only on the size of the complement of $T_{v_N}$ in
  $g' \cdot \Gamma_p$), we can find $\alpha_N \in T_{v_N}$ with
  \[
    \alpha_1 \cdots \alpha_N = \gamma.
  \]
\end{proof}

\begin{remark}
  In general, we do \emph{not} require the set of elements in $\Gamma$
  labelling the vertices of a relative quasigeodesic automaton
  $\mc{G}$ to generate the group $\Gamma$ (although the proposition
  above implies that they at least generate a finite-index subgroup).
\end{remark}

\subsection{Compatible systems of open sets}

A relative quasigeodesic automaton always exists for any relatively
hyperbolic group (although we will not prove this fact in full
generality). We will give a way to construct a relative quasigeodesic
automaton using the convergence group action of a group acting on its
Bowditch boundary.

\begin{definition}
  Suppose that $\Gamma$ acts on a metrizable space $M$ by
  homeomorphisms, and let $\mc{G}$ be a $\Gamma$-graph. A
  \emph{$\mc{G}$-compatible system of open sets} for the action of
  $\Gamma$ on $M$ is an assignment of an open subset $U_v \subset M$
  to each vertex $v$ of $\mc{G}$ such that for each edge $e = (v,w)$
  in $\mc{G}$, for some $\eps > 0$, we have
  $\overline{N_M(U_w, \eps)} \ne M$ and
  \begin{equation}
    \label{eq:compatibility}
    \alpha \cdot N_M(U_w, \eps)  \subset U_v
  \end{equation}
  for all $\alpha \in T_v$.
\end{definition}

\begin{remark}
  If $\mc{G}$ has no parabolic vertices (so each set $T_v$ contains a
  single group element $\alpha_v \in \Gamma$), then
  (\ref{eq:compatibility}) is equivalent to requiring
  $\alpha_v \cdot \overline{U_w} \subset U_v$ for every edge $(v, w)$
  in $\mc{G}$. When $\mc{G}$ has parabolic vertices (so $T_v$ may be
  infinite), (\ref{eq:compatibility}) may be a stronger condition.
\end{remark}

\begin{prop}
  \label{prop:bowditch_graphs_geodesic}
  Let $\mc{G}$ be a $\Gamma$-graph, and let
  $\{U_v : v \textrm{ vertex of \mc{G}}\}$ be a $\mc{G}$-compatible
  system of subsets of $\bgamh$ for the action of
  $\Gamma$ on $\bgamh$.

  There is a Gromov model $Y$ with basepoint $y_0$ satisfying the
  following. Let $\{\alpha_n\}$ be an infinite $\mc{G}$-path,
  corresponding to a vertex path $\{v_n\}$, and suppose the sequence
  $\{\gamma_n = \alpha_1 \cdots \alpha_n\}$ is divergent in
  $\Gamma$. Then for any point $z$ in the intersection
  \[
    U_\infty = \bigcap_{n=1}^\infty \alpha_1 \cdots \alpha_n U_{v_{n+1}},
  \]
  there is some $R > 0$ and a geodesic ray $c$ in $Y$ from $y_0$ to
  $z$ such that the sequence $\gamma_ny_0$ lies in an $R$-neighborhood
  of $c$. In particular, $z$ is a conical limit point.
\end{prop}
\begin{proof}
  Fix a point $z \in U_\infty$, and write $z = z_+$ and
  $U_n = U_{v_n}$. We first claim that there is a uniform $\eps > 0$
  and a point $z_- \in \bgamh$ such that
  \begin{equation}
    \label{eq:endpoints_not_close}
    d(\gamma_n^{-1} z_+, \gamma_n^{-1} z_-) \ge \eps
  \end{equation}
  for all $n \ge 0$.

  To prove the claim, choose a uniform $\eps > 0$ so that for every
  vertex $v$ in $\mc{G}$, we have
  $\overline{N(U_v, \eps)} \ne \bgamh$, and for every edge $(v, w)$ in
  $\mc{G}$ and every $\alpha \in T_v$, we have
  $\alpha \cdot N(U_w, \eps) \subset U_v$. Then we choose some
  $z_- \in \bgamh - \overline{N(U_1, \eps)}$.

  By the $\mc{G}$-compatibility condition, we know that for any $n$,
  $\gamma_nU_{n+1} \subset \ldots \subset \gamma_1U_2 \subset U_1$, so
  we know that $d(z_+, z_-) \ge \eps$.
  
  Then, for any $n \ge 1$, we have
  \[
    \gamma_n^{-1} z_+ \in U_{n+1}.
  \]
  Moreover since $\gamma_n N(U_{n+1}, \eps) \subset U_1$, we also have
  \[
    \gamma_n^{-1} z_- \in \bgamh - N(U_{n+1}, \eps).
  \]
  So for all $n$ we have
  $d(\gamma_n^{-1}z_+, \gamma_n^{-1}z_-) \ge \eps$, which establishes
  that \eqref{eq:endpoints_not_close} holds for all $n$.

  Now, consider a bi-infinite geodesic $c$ in a Gromov model $Y$ for
  $\Gamma$ joining $z_+$ and $z_-$. The sequence of geodesics
  $\gamma_n^{-1} \cdot c$ has endpoints in $\dee Y = \bgamh$ lying
  distance at least $\eps$ apart, so each geodesic in the sequence
  passes within a uniformly bounded neighborhood of a fixed basepoint
  $y_0 \in Y$. Equivalently every $\gamma_n \cdot y_0$ lies within a
  fixed neighbood of the geodesic $c$.

  Since $\gamma_n$ is divergent, $\gamma_n y_0$ can only accumulate at
  either $z_+$ or $z_-$. But in fact $\gamma_n y_0$ can only
  accumulate at $z_+$---for in the construction of $c$ above, we could
  have chosen any $z_-$ in the nonempty open set
  $\bgamh - \overline{N}(U_1, \eps)$, and since $\bgamh$ is perfect
  there is at least one such $z_-' \ne z_-$. This implies that the set
  $\{\gamma_ny_0\}$ actually lies within a neighborhood of a ray with
  ideal endpoint $z$, and therefore within a neighborhood of a ray
  from $y_0$ to $z$.
\end{proof}

\begin{definition}
  Let $\mc{G}$ be a $\Gamma$-graph. An infinite $\mc{G}$-path
  $\{\alpha_n\}$ is \emph{divergent} if the sequence
  $\{\gamma_n = \alpha_1 \cdots \alpha_n\}$ leaves every bounded
  subset of $\Gamma$.

  We say that a $\Gamma$-graph $\mc{G}$ is divergent if \emph{every}
  infinite $\mc{G}$-path is divergent.
\end{definition}

Whenever $\{U_v\}$ is a $\mc{G}$-compatible system of open sets for a
$\Gamma$-graph $\mc{G}$, one can think of a $\mc{G}$-path
$\{\alpha_n\}_{n=1}^N$ as giving a \emph{symbolic coding} of a point
in the intersection
\[
  \bigcap_{n=1}^N\alpha_1 \cdots \alpha_n U_{n+1}.
\]
In general, the coding for a given point in this intersection may not
be unique. However, the next result says that (at least for divergent
$\mc{G}$-paths) the coding is \emph{coarsely} unique.
\begin{lem}
  \label{lem:codings_hausdorff_close}
  Let $\mc{G}$ be a $\Gamma$-graph, and let
  $\{U_v : v \textrm{ vertex of \mc{G}}\}$ be a $\mc{G}$-compatible
  system of subsets of $\bgamh$ for the action of $\Gamma$ on
  $\bgamh$. Suppose that $\{\alpha_n\}$ is a divergent $\mc{G}$-path
  limiting to $z \in \bgamh$, and $\{\beta_n\}$ is a divergent
  $\mc{G}$-path limiting to $g_0z$ for some $g_0 \in \Gamma$. Then the
  Hausdorff distance in $\cay(\Gamma)$ between the sets
  \[
    \{g_0\alpha_1 \cdots \alpha_n\}_{n \in \N}, \qquad \{\beta_1
    \cdots \beta_m\}_{m \in \N}
  \]
  is finite.
\end{lem}
Actually, it is possible to say more: the Hausdorff distance between
the sets in question is bounded by a constant which depends only on
the group element $g \in \Gamma$ (see e.g. \cite[Lemma
4.15]{MMW2}). The statement above is sufficient for our purposes,
however.

Before giving the proof of the above lemma, we need a basic result in
hyperbolic geometry.
\begin{lem}
  \label{lem:hyp_geom_balls}
  Let $Y$ be a $\delta$-hyperbolic metric space, and let $B \subset Y$
  be a horoball with $Y - B$ nonempty. Suppose that $x, y \in Y$ lie
  within distance $R > 0$ of $B$, and that $z \in Y$ lies on a
  geodesic $[x,y] \subset Y$ joining $x$ to $y$. Then we have
  \[
    d(z, Y - B) \ge d(z, \{x,y\}) - R - 4\delta.
  \]
\end{lem}
\begin{proof}
  It suffices to prove the claim when $B$ is instead a metric ball in
  $Y$, since horoballs are local Hausdorff limits of these. So, let
  $B$ be a ball in $Y$ with center $o$ and radius $r$, and assume that
  $x,y$ lie in an $R$-neighborhood of $B$. Let $w \in Y - B$, and
  consider a pair of geodesic triangles with vertices $o, x, y$ and
  $w, x, y$, each with an edge equal to a fixed geodesic
  $[x,y]$. Since $Y$ is $\delta$-hyperbolic, for any $z \in [x,y]$,
  there are points $z' \in [o, x] \cup [o,y]$ and
  $z'' \in [w, x] \cup [w,y]$ with $d(z, z') < \delta$ and
  $d(z, z'') < \delta$. Without loss of generality $z' \in [o,x]$, so
  that we have
  \[
    d(z, x) < \delta + d(z', x) = d(o, x) - d(o, z') + \delta \le r +
    R - d(o, z') + \delta.
  \]
  On the other hand we also have
  \[
    r \le d(o, w) < d(o, z') + d(z'', w) + 2\delta,
  \]
  hence
  \[
    d(z, w) > d(z'', w) - \delta \ge r - d(o, z') - 3\delta > d(z, x) -
    R - 4\delta.
  \]
\end{proof}
\begin{proof}[Proof of \Cref{lem:codings_hausdorff_close}]
  Let $\{g_n\}_{n \in \N \cup \{0\}}$ and $\{h_m\}_{m \in \N}$
  respectively denote the sets
  \[
    \{g\alpha_1 \cdots \alpha_n\}_{n \in \N \cup \{0\}}, \qquad
    \{\beta_1 \cdots \beta_m\}_{m \in \N}.
  \]
  Fix a Gromov model $Y$ with basepoint $y_0$ for $(\Gamma, \mc{H})$,
  and let $\mc{B}$ be an invariant system of horoballs in $Y$. From
  \Cref{prop:bowditch_graphs_geodesic}, we know that there is a
  geodesic ray $c$ in $Y$ from $y_0$ to $gz$ and a constant $R > 0$ so
  that both sets $\{g_ny_0\}_{n \in \N}$ and $\{h_my_0\}_{m \in \N}$
  lie within distance $R$ of $c$. We claim that the Hausdorff distance
  in $Y$ between these two sets is actually bounded.

  We will show that there is a uniform constant $D > 0$ so that if $a$
  is any point belonging to the set $\{h_my_0\}_{m \in \N}$, then $a$
  lies within distance $D$ of some point in the set
  $\{g_ny_0\}_{n \in \N \cup \{0\}}$; the other argument is nearly
  symmetric. So, fix such an $a$, and let $a'$ be a point on $c$ lying
  within distance $R$ of $a$. Then, for each $n \ge 1$, let
  $y_n = g_ny_0$, and let $y_n'$ be a point on $c$ within distance $R$
  of $y_n$. We also define $y_0' = y_0$.

  Since $\{g_n\}$ is divergent and the action of $\Gamma$ on $Y$ is
  proper, the $y_n$ (hence the $y_n'$) give an unbounded subset of
  $Y$. It follows that there is some $n \ge 0$ such that $a'$ lies on
  the closed sub-geodesic of $c$ joining the points $y_n', y_{n+1}'$.

  By definition we know that $g_{n+1} = g_n\alpha_{n+1}$. Now, since
  $\alpha_{n+1}$ belongs to one of the finitely many label sets for
  the $\Gamma$-graph $\mc{G}$, it belongs to one of finitely many
  cosets of the form $g\Gamma_p$ for $g \in \Gamma$ and $p \in
  \Pi$. Since each $\Gamma_p$ preserves a horoball in the
  $\Gamma$-invariant system $\mc{B}$, there is a constant $R' > 0$
  (independent of $n$) and a horoball $B \in \mc{B}$ so that $y_0$ and
  $\alpha_{n+1}y_0$ both lie within distance $R'$ of $B$.

  Consequently the points $y_n'$ and $y_{n+1}'$ both lie within
  distance $R' + R$ of the horoball $g_nB \in \mc{B}$. Thus, by
  \Cref{lem:hyp_geom_balls}, we have
  \[
    d(a', Y - g_nB) \ge d(a', \{y_n', y_{n+1}'\}) - R' - R - 4\delta.
  \]
  Now, we may assume that our basepoint $y_0$ does not lie in any
  horoball in $\mc{B}$, which means that
  $d(a', Y - g_nB) \le d(a', a) \le R$. Therefore we have
  \[
    d(a', \{y_n', y_{n+1}'\}) \le R' + 2R + 4\delta,
  \]
  hence
  \[
    d(a, \{y_n, y_{n+1}\}) \le R' + 4R + 4\delta.
  \]
  This proves that the Hausdorff distance (in $Y$) between the sets
  $\{g_ny_0\}_{n \in \N \cup \{0\}}$ and $\{h_ny_0\}_{n \in \N}$ is
  bounded. The analogous claim for the corresponding subsets of
  $\cay(\Gamma)$ then follows from the fact that $\Gamma$ acts
  properly on $Y$.
\end{proof}

Now that we have established (coarse) uniqueness for $\Gamma$-graph
codings, we turn to the question of existence. The proposition below
says that, given an appropriate pair of open coverings of the Bowditch
boundary $\bgamh$ compatible with a $\Gamma$-graph $\mc{G}$, one can
always construct a $\mc{G}$-path limiting to a given point
$z \in \bgamh$.
\begin{prop}
  \label{prop:relative_geodesic_automaton_from_bb}
  Let $\mc{G}$ be a \emph{divergent} $\Gamma$-graph. Suppose that for
  each vertex $a \in \mc{G}$, there exist open subsets $V_a, W_a$ of
  $\bgamh$ such that the following conditions hold:
  \begin{enumerate}
  \item\label{item:g_compatible_system} The sets $\{W_a\}$ give a
    $\mc{G}$-compatible system of sets for the action of $\Gamma$ on
    $\bgamh$.
  \item\label{item:closure_proper} For all vertices $a$, we have
    $V_a \subset W_a$ and $\overline{W_a} \ne \bgamh$.
  \item The sets $V_a$ give an open covering of $\bgamh$.
  \item \label{item:conical_edge_exists} For every $z \in \bgamh$
    and every non-parabolic vertex $a$ such that $z \in V_a$, there is
    an edge $(a, b)$ in $\mc{G}$ such that
    $\alpha_a^{-1} \cdot z \in V_b$ for $\{\alpha_a\} = T_a$.
  \item \label{item:parabolic_edge_exists} For every $z \in \bgamh$
    and every parabolic vertex $a$ such that $z \in V_a - \{p_a\}$,
    there is an edge $(a, b)$ in $\mc{G}$ and $\alpha \in T_a$ such
    that $\alpha^{-1} \cdot z \in V_b$.
  \end{enumerate}
  Then $\mc{G}$ is a relative quasigeodesic automaton for $\Gamma$.
\end{prop}

\begin{figure}[h]
  \centering
  \def\svgwidth{8cm}
  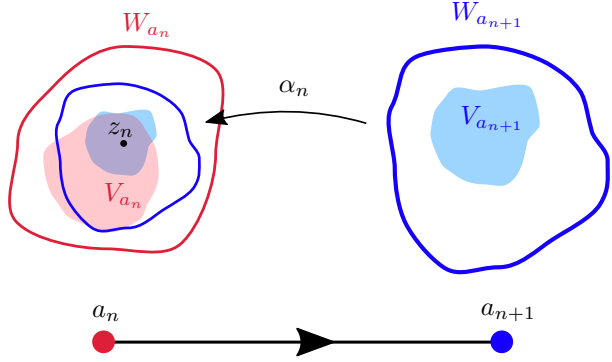
  \caption{Illustration for the proof of
    \Cref{prop:relative_geodesic_automaton_from_bb}. The group element
    $\alpha_n$ nests an $\eps$-neighborhood of $W_{a_{n+1}}$ inside of
    $W_{a_n}$ whenever $\alpha_n \cdot V_{a_{n+1}}$ intersects
    $V_{a_n}$.}
  \label{fig:g_compatible_inclusion}
\end{figure}
\begin{figure}[h]
  \centering
  \def\svgwidth{.8\textwidth}
  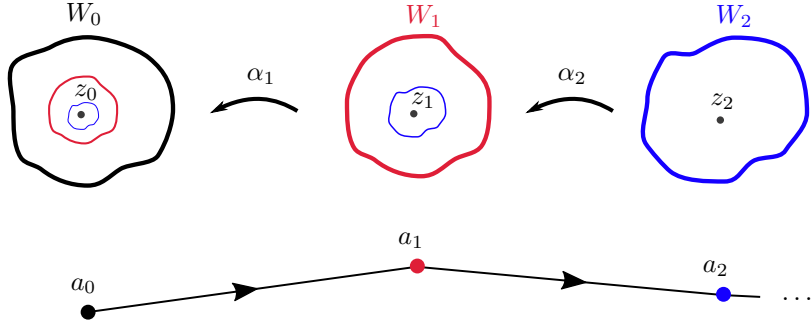
  \caption{By iterating the nesting procedure backwards, we produce an
    infinite $\mc{G}$-path and a sequence of subsets intersecting in
    the initial point $z = z_0$.}
  \label{fig:nested_subsets}
\end{figure}
\begin{proof}
  \Cref{prop:bowditch_graphs_geodesic} establishes the first part of
  the definition of a relative quasi-geodesic automaton
  (\Cref{defn:rel_automaton}), so we just need to show the second
  part, stating that every $z \in \bgamh$ is the limit of a
  $\mc{G}$-path.

  The idea behind the proof is to use the fact that the sets $V_a$
  cover $\bgamh$ to show that we can keep ``expanding''
  a neighborhood of $z$ in $\bgamh$ to construct a path
  in $\mc{G}$ limiting to $z$. The $\{V_a\}$ covering tells us how to
  find the next edge in the path, and the $\{W_a\}$ cover gives us the
  $\mc{G}$-compatible system we need to show that the path is a
  geodesic. 

  We let $A$ denote the vertex set of $\mc{G}$. When $a \in A$ is not
  a parabolic vertex, we write $T_a = \{\gamma_a\}$.

  \begin{aside}{Case 1: $z$ is a conical limit point}
    Fix $a \in A$ so that $z \in V_a$. We take $z_0 = z$, $a_0 = a$,
    and define sequences
    $\{z_n\}_{n=0}^\infty \subset \conbdry \Gamma$,
    $\{a_n\}_{n=0}^\infty \subset A$, and
    $\{\alpha_n\}_{n=1}^\infty \subset \Gamma$ as follows:
    \begin{itemize}
    \item If $a_n$ is not a parabolic vertex, then we choose
      $\alpha_{n+1} = \gamma_{a_n}$. Let
      $z_{n+1} = \alpha_{n+1}^{-1} \cdot z_n$. Since conical limit
      points are invariant under the action of $\Gamma$, $z_{n+1}$ is
      a conical limit point. By condition
      \ref{item:conical_edge_exists}, there is a vertex $a_{n+1}$
      satisfying $z_{n+1} \in V_{a_{n+1}}$ with $(a_n, a_{n+1})$ an
      edge in $\mc{G}$.
    \item If $a_n$ is a parabolic vertex, then since $z_n$ is a
      conical limit point, $z_n \ne p$ for $p = p_{a_n}$. Then
      condition \ref{item:parabolic_edge_exists} implies that there
      exists some $\alpha_{n+1} \in T_{a_n}$ so that
      $\alpha_{n+1}^{-1} \cdot z_n \in V_{a_{n+1}}$ for an edge
      $(a_n, a_{n+1})$ in $\mc{G}$. Again,
      $z_{n+1} = \alpha_{n+1}^{-1} \cdot z_n$ must be a conical limit
      point since $\conbdry\Gamma$ is $\Gamma$-invariant.
    \end{itemize}

    The sequence $\{\alpha_n\}$ necessarily gives a $\mc{G}$-path. By
    assumption the sequence
    \[
      \gamma_n = \alpha_1 \cdots \alpha_n
    \]
    is divergent. And by construction $z = \gamma_nz_n$ lies in
    $\gamma_n W_{a_n}$ for all $n$. So,
    \Cref{prop:bowditch_graphs_geodesic} implies that $\gamma_{n}$ is
    a conical limit sequence, limiting to $z$. See
    \Cref{fig:nested_subsets}.
  \end{aside}
  \begin{aside}{Case 2: $z$ is a parabolic point}
    As before fix $a \in A$ so that $z \in V_a$, and take $z_0 = z$,
    $a_0 = a$. We inductively define sequences $z_n$, $a_n$,
    $\alpha_n$ as before, but we claim that for some finite $N$, $a_N$
    is a parabolic vertex with $z_N = p_{a_N}$. For if not, we can
    build an infinite $\mc{G}$-path (as in the previous case) limiting
    to $z$. But then, \Cref{prop:bowditch_graphs_geodesic} would imply
    that $z$ is actually a conical limit point. So, we must have
    \[
      z = \gamma_N p_{a_N} = \alpha_1 \cdots \alpha_N p_{a_N}
    \]
    as required.
  \end{aside}
\end{proof}

\begin{remark}
  In a typical application of
  \Cref{prop:relative_geodesic_automaton_from_bb}, it will not be
  possible to choose the open coverings $\{V_a\}$ and $\{W_a\}$ so
  that they exactly agree. In particular we expect this to be
  impossible whenever $\bgamh$ is connected. To get an idea why,
  assume that $V_a = W_a$ for all vertices $a$; then conditions
  \ref{item:g_compatible_system} and \ref{item:conical_edge_exists} in
  the proposition together imply that whenever $a$ is a non-parabolic
  vertex with $T_a = \{\alpha_a\}$, the set $\alpha_a^{-1}W_a$ can be
  written as a finite union of sets of the form $\overline{W_b}$. Thus
  in this case $W_a$ must be a proper clopen subset of $\bgamh$.
\end{remark}

To conclude this section, we make one more observation about systems
of $\mc{G}$-compatible sets as in
\Cref{prop:relative_geodesic_automaton_from_bb}.

\begin{lem}
  \label{lem:close_points_share_prefixes}
  Let $\Gamma$ be a relatively hyperbolic group, let $\mc{G}$ be a
  $\Gamma$-graph, and let $\{V_a\}$, $\{W_a\}$ be an assignment of
  open subsets of $\bgamh$ to vertices of $\mc{G}$ satisfying the
  hypotheses of \Cref{prop:relative_geodesic_automaton_from_bb}.

  Fix $z \in \conbdry \Gamma$ and $N \in \N$. There exists
  $\delta > 0$ so that if $d(z, z') < \delta$, then there are
  $\mc{G}$-paths $\{\alpha_n\}, \{\beta_n\}$ limiting to $z, z'$
  respectively, with $\alpha_i = \beta_i$ for all $i < N$.
\end{lem}
\begin{proof}
  Let $\{\alpha_n\}$ be a $\mc{G}$-path limiting to $z$ coming from
  the construction in \Cref{prop:relative_geodesic_automaton_from_bb},
  passing through vertices $v_n$. We choose $\delta > 0$ small enough
  so that if $d(z, z') < \delta$, then $z'$ lies in the set
  \[
    \alpha_1 \cdots \alpha_N V_{v_{N+1}}.
  \]
  Then for every $i < N$, we have
  \[
    \alpha_i^{-1} \alpha_{i-1}^{-1} \cdots \alpha_1^{-1} z' \in
    V_{v_{i+1}}.
  \]
  As in \Cref{prop:relative_geodesic_automaton_from_bb}, we can then
  extend $\{\alpha_n\}_{n=1}^{N-1}$ to a $\mc{G}$-path limiting to
  $z'$.
\end{proof}


\section{Extended convergence dynamics}
\label{sec:extended_convergence}
Let $\Gamma$ be a relatively hyperbolic group acting on a connected
compact metrizable space $M$. In this section, we will show that if
the action of $\Gamma$ on $M$ extends the convergence action of
$\Gamma$ (\Cref{defn:extends_convergence_dynamics}), then we can
construct a relative quasigeodesic automaton $\mc{G}$ and a
$\mc{G}$-compatible system of open subsets of $M$ which are in some
sense reasonably well-behaved with respect to the group action.

To give the precise statement, let $\Lambda \subset M$ be a closed
$\Gamma$-invariant subset, and let $\phi:\Lambda \to \bgamh$ be an
equivariant, surjective, and continuous map. Assume that there is an
assigment $z \mapsto C_z$ of points in $\bgamh$ to open subsets of
$M$, so that the map $\phi$ and sets $C_z$ satisfy property
\ref{item:strata_interior} in \Cref{defn:convergence_group}, as well
as both of the following:

\begin{enumerate}
  \custitem[(C2-con)]\label{item:conical_extended_convergence} For any
  sequence $\gamma_n \in \Gamma$ limiting conically to $z$, with
  $\gamma_n^{-1} \to z_-$, any open set $U$ containing $\phi^{-1}(z)$,
  and any compact $K \subset C_{z_-}$, we have
  $\gamma_n \cdot K \subset U$ for all sufficiently large $n$.
  \custitem[(C2-par)]\label{item:peripheral_extended_convergence} For
  any parabolic point $p$, any compact $K \subset C_p$, and any open
  set $U$ containing $\phi^{-1}(p)$, for all but finitely many
  $\gamma \in \Gamma_p$, we have $\gamma \cdot K \subset U$.
\end{enumerate}
Note that \ref{item:conical_extended_convergence} and
\ref{item:peripheral_extended_convergence} are both special cases of
condition \ref{item:extended_convergence} in
\Cref{defn:extends_convergence_dynamics}, so any map $\phi$ and strata
$C_z$ satisfying \Cref{defn:extends_convergence_dynamics} fulfill
these requirements. For the rest of this section, however, we
\emph{only} assume that \ref{item:strata_interior},
\ref{item:conical_extended_convergence} and
\ref{item:peripheral_extended_convergence} hold for $\phi$ and
$C_z$. In this context, we show:
\begin{prop}
  \label{prop:convergence_dynamics_implies_automaton}
  For any $\eps > 0$, there is a relative quasigeodesic automaton
  $\mc{G}$ for $\Gamma$, a $\mc{G}$-compatible system of open sets
  $\{U_v\}$ for the action of $\Gamma$ on $M$, and a
  $\mc{G}$-compatible system of open sets $\{W_v\}$ for the action of
  $\Gamma$ on $\bgamh$ such that:
  \begin{enumerate}
  \item For every vertex $v$, there is some $z \in W_v$ so that
    \[
      \phi^{-1}(W_v) \subset U_v \subset N_M(\phi^{-1}(z), \eps).
    \]
  \item For every $p \in \Pi$, there is a parabolic vertex $a$ with
    $p_a = p$. Moreover, for every parabolic vertex $w$ with
    $p_w = g \cdot p$, $(a, b)$ is an edge of $\mc{G}$ if and only if
    $(w, b)$ is an edge of $\mc{G}$.
  \item If $q = g \cdot p$ for $p \in \Pi$, $a$ is a parabolic vertex
    with $p_a = q$, and $(a,b)$ is an edge of $\mc{G}$, then
    $q \in W_a$ and $U_b \subset C_p$.
  \end{enumerate}
\end{prop}

\begin{remark}
  By equivariance of $\phi$, for each $p \in \parbdry \Gamma$, we can
  replace $C_p$ with $\Gamma_p C_p$ and assume that $C_p$ is
  $\Gamma_p$-invariant, and that if $q,p \in \parbdry\Gamma$ satisfy
  $q = gp$ for some $g \in \Gamma$, then $C_q = gC_p$.
\end{remark}

The proof of \Cref{prop:convergence_dynamics_implies_automaton}
involves some technicalities, so we first outline the general
approach:
\begin{enumerate}
\item For each $z \in \bgamh$, we construct a pair $V_z$, $W_z$
  of small open neighborhoods of $z$ and a subset $T_z \subset \Gamma$
  so that for each $\alpha \in T_z$, $\alpha^{-1}$ is ``expanding''
  about some point in $V_z$. When $z$ is a conical limit point, then
  we can choose a single element $\alpha_z \in \Gamma$ which expands
  about every point in $V_z$. When $z$ is a parabolic point, we may
  use a different element of $\Gamma$ to ``expand'' about each
  $u \in V_z - \{z\}$.

  We choose $V_z$, $W_z$, and $T_z$ so that if $\alpha^{-1}$ is
  ``expanding'' about $u \in V_z$, and $\alpha^{-1}u \in V_y$, then
  $\alpha^{-1}W_z \supset W_y$. See \Cref{fig:expansion_inclusions}.

  \begin{figure}[h]
  \centering
  \def\svgwidth{.6\textwidth}
  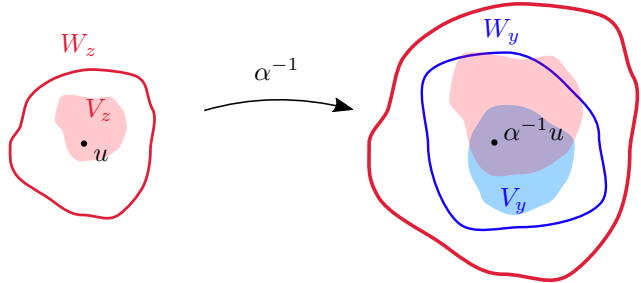
  \caption{The group element $\alpha^{-1}$ is ``expanding'' about
    $u \in V_z$. We will construct $V_z$, $W_z$ and $V_y, W_y$ so that
    if $\alpha^{-1} u$ lies in $V_y$, then $\alpha^{-1} W_z$ contains
    $W_y$. Equivalently, we get the containment
    $\alpha W_y \subset W_z$ illustrated earlier in
    \Cref{fig:g_compatible_inclusion}.}
  \label{fig:expansion_inclusions}
\end{figure}

\item Using compactness of $\bgamh$, we pick a finite set of
  points $a \in \bgamh$ so that the sets $\{V_a\}$ give an open
  covering of $\bgamh$. These points in $\bgamh$ are
  identified with the vertices of a $\Gamma$-graph $\mc{G}$. We define
  the edges of $\mc{G}$ in such a way so that if, for some
  $\alpha \in T_a$, $\alpha^{-1}$ expands about $u \in V_a$ and
  $\alpha^{-1}u \in V_b$, then there is an edge from $a$ to $b$. This
  ensures that $\{W_a\}$ is a $\mc{G}$-compatible system of open
  subsets of $\bgamh$.

\item Simultaneously, we construct a $\mc{G}$-compatible system
  $\{U_a\}$ of open sets in $M$ by taking $U_a$ to be a small
  neighborhood of $\phi^{-1}(a)$. The idea is to use the extended
  convergence dynamics to ensure that if, for some $\alpha \in T_z$,
  $\alpha^{-1}$ ``expands'' about some $u \in V_z$ and the point
  $\alpha^{-1} u$ lies in $V_y$, then $\alpha^{-1}U_z$ contains
  $U_y$. See \Cref{fig:conical_convergence_setup} below.

\item Finally, we use \Cref{prop:relative_geodesic_automaton_from_bb}
  to prove that $\mc{G}$ is actually a relative quasigeodesic
  automaton. The open sets $V_a, W_a$ are constructed exactly to
  satisfy the conditions of the proposition, so the main thing to
  check in this step is that the graph $\mc{G}$ is actually divergent
  (using the action of $\Gamma$ on $M$).
\end{enumerate}

Throughout the rest of the section, we will work with fixed metrics on
both $\bgamh$ and $M$. Critically, none of our ``expansion'' arguments
will depend sensitively on the precise choice of metric. That is, in
the sketch above, when we say that some group element
$\alpha \in \Gamma$ ``expands'' on a small open subset $U$ of a metric
space $X$, we just mean that $\alpha U$ is quantifiably ``bigger''
than $U$, and \emph{not} that for any $x, y \in U$, we have
$d(\alpha \cdot x, \alpha \cdot y) \ge C \cdot d(x, y)$ for some
expansion constant $C$. \Cref{lem:conical_limit_neighborhoods} and
\Cref{lem:large_parabolic_neighborhood} below describe precisely what
we mean by ``bigger.'' The general idea is captured by the following
example.

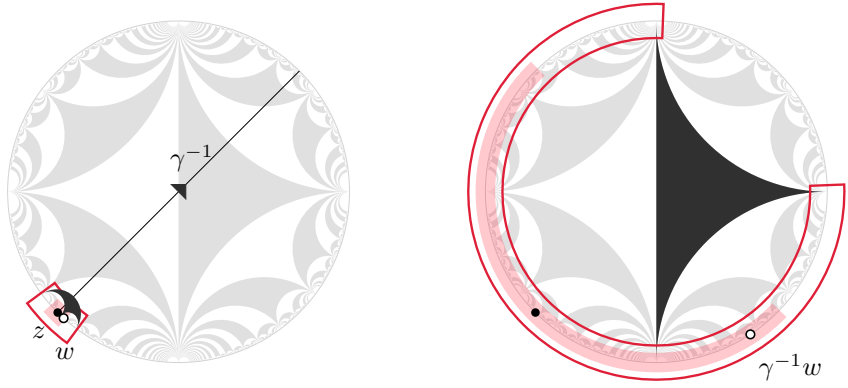
\begin{figure}[h]
  \centering
  {
    \def\svgwidth{6cm}
\begingroup%
  \makeatletter%
  \providecommand\color[2][]{%
    \errmessage{(Inkscape) Color is used for the text in Inkscape, but the package 'color.sty' is not loaded}%
    \renewcommand\color[2][]{}%
  }%
  \providecommand\transparent[1]{%
    \errmessage{(Inkscape) Transparency is used (non-zero) for the text in Inkscape, but the package 'transparent.sty' is not loaded}%
    \renewcommand\transparent[1]{}%
  }%
  \providecommand\rotatebox[2]{#2}%
  \newcommand*\fsize{\dimexpr\f@size pt\relax}%
  \newcommand*\lineheight[1]{\fontsize{\fsize}{#1\fsize}\selectfont}%
  \ifx\svgwidth\undefined%
    \setlength{\unitlength}{576bp}%
    \ifx\svgscale\undefined%
      \relax%
    \else%
      \setlength{\unitlength}{\unitlength * \real{\svgscale}}%
    \fi%
  \else%
    \setlength{\unitlength}{\svgwidth}%
  \fi%
  \global\let\svgwidth\undefined%
  \global\let\svgscale\undefined%
  \makeatother%
  \begin{picture}(1,1)%
    \lineheight{1}%
    \setlength\tabcolsep{0pt}%
    \put(0,0){\includegraphics[width=\unitlength,page=1]{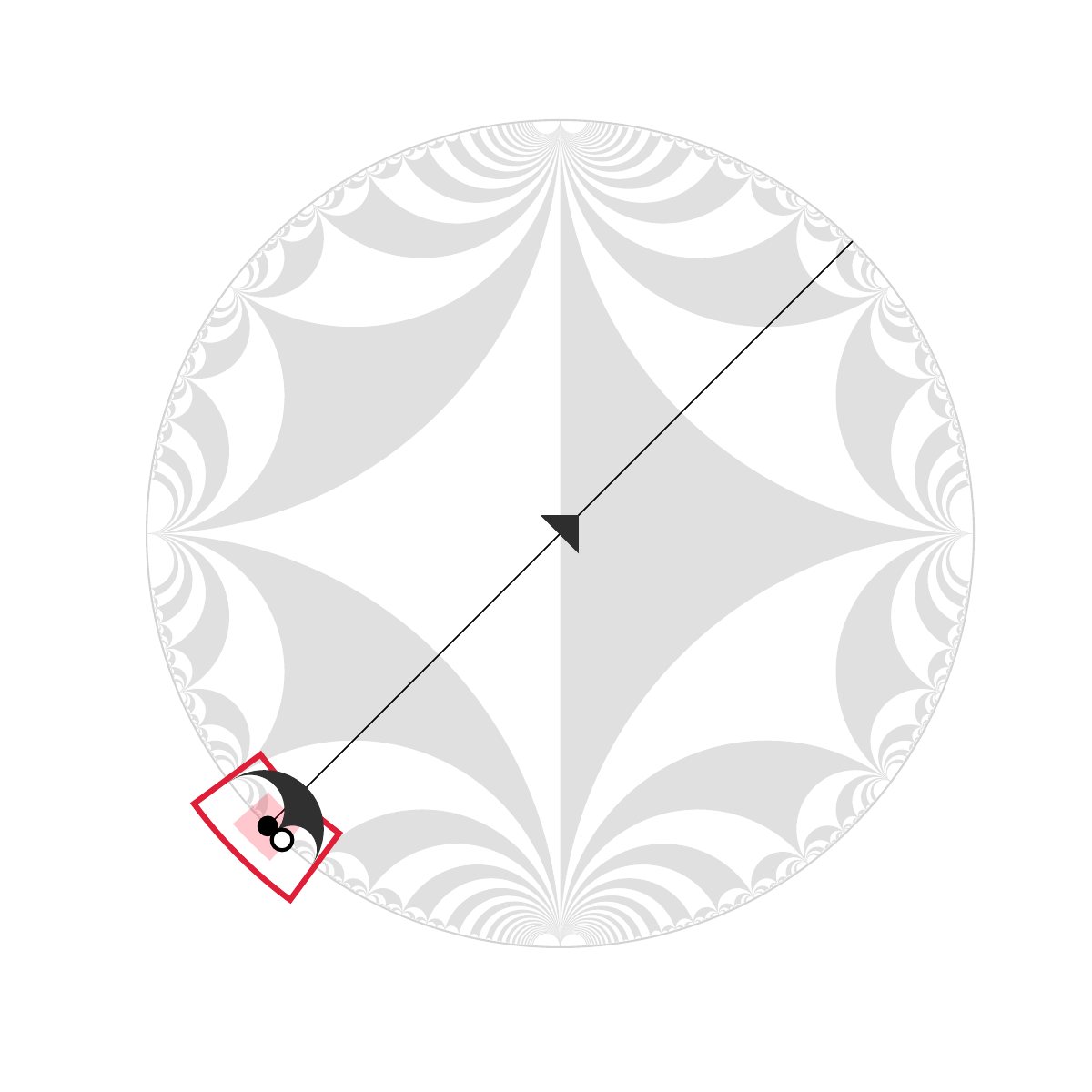}}%
    \put(0.49402894,0.56809759){\makebox(0,0)[lt]{\lineheight{1.25}\smash{\begin{tabular}[t]{l}$\gamma^{-1}$\end{tabular}}}}%
    \put(0.18833522,0.18665335){\makebox(0,0)[lt]{\lineheight{1.25}\smash{\begin{tabular}[t]{l}$z$\end{tabular}}}}%
    \put(0.23936929,0.14002376){\makebox(0,0)[lt]{\lineheight{1.25}\smash{\begin{tabular}[t]{l}$w$\end{tabular}}}}%
  \end{picture}%
\endgroup%

  }{
    \def\svgwidth{6cm}
\begingroup%
  \makeatletter%
  \providecommand\color[2][]{%
    \errmessage{(Inkscape) Color is used for the text in Inkscape, but the package 'color.sty' is not loaded}%
    \renewcommand\color[2][]{}%
  }%
  \providecommand\transparent[1]{%
    \errmessage{(Inkscape) Transparency is used (non-zero) for the text in Inkscape, but the package 'transparent.sty' is not loaded}%
    \renewcommand\transparent[1]{}%
  }%
  \providecommand\rotatebox[2]{#2}%
  \newcommand*\fsize{\dimexpr\f@size pt\relax}%
  \newcommand*\lineheight[1]{\fontsize{\fsize}{#1\fsize}\selectfont}%
  \ifx\svgwidth\undefined%
    \setlength{\unitlength}{576bp}%
    \ifx\svgscale\undefined%
      \relax%
    \else%
      \setlength{\unitlength}{\unitlength * \real{\svgscale}}%
    \fi%
  \else%
    \setlength{\unitlength}{\svgwidth}%
  \fi%
  \global\let\svgwidth\undefined%
  \global\let\svgscale\undefined%
  \makeatother%
  \begin{picture}(1,1)%
    \lineheight{1}%
    \setlength\tabcolsep{0pt}%
    \put(0,0){\includegraphics[width=\unitlength,page=1]{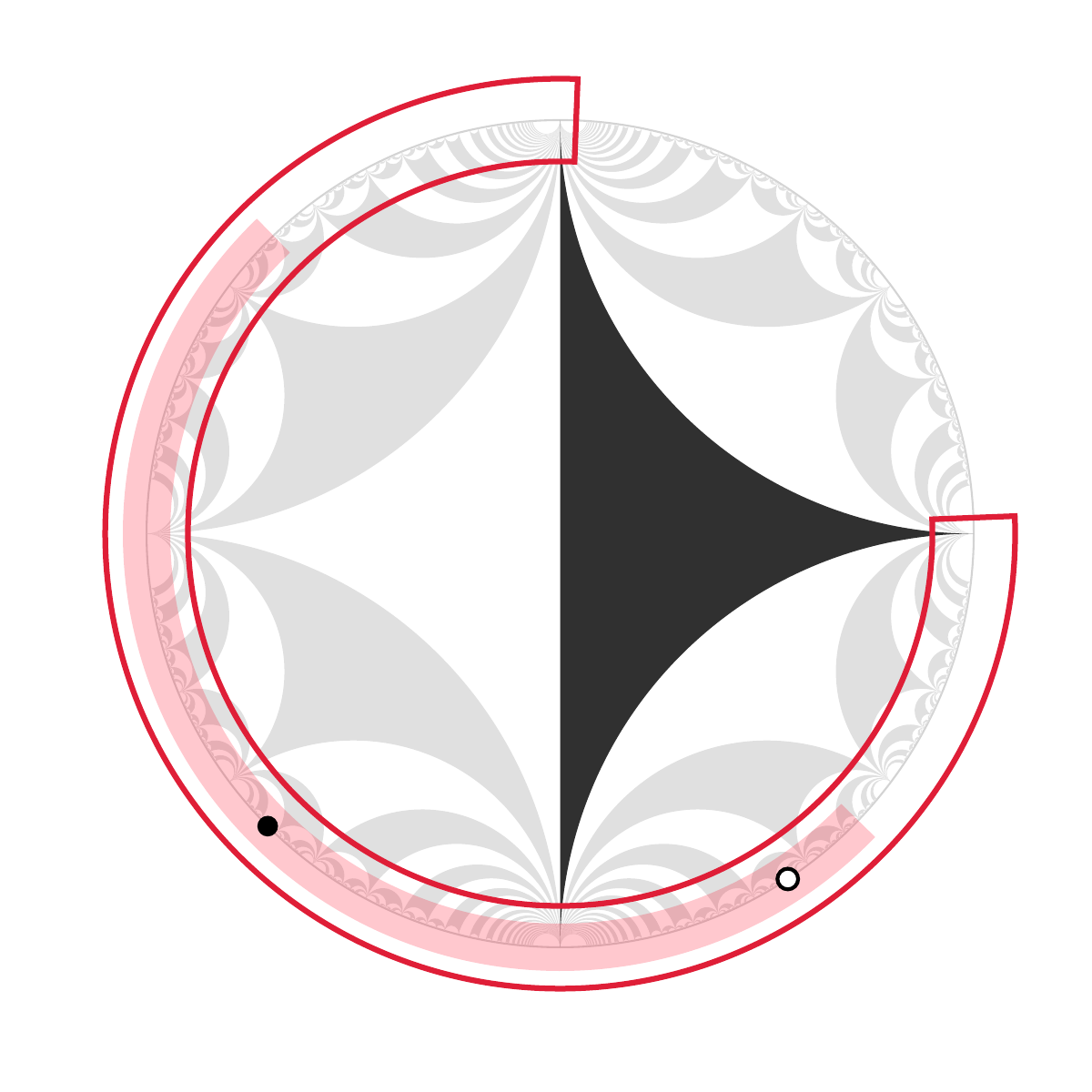}}%
    \put(0.74035078,0.10415077){\makebox(0,0)[lt]{\lineheight{1.25}\smash{\begin{tabular}[t]{l}$\gamma^{-1} w$\end{tabular}}}}%
  \end{picture}%
\endgroup%

  }
  \caption{For any point $w$ in a sufficiently small neighborhood
    $V_z$ (pink) of $z$, the expanded neighborhood $\gamma^{-1}W_z$
    (red) contains a uniform neighborhood of $\gamma^{-1} w$.}
  \label{fig:sl2z_loxodromic_expansion}
\end{figure}

\begin{figure}[h]
  \centering
  {
    \def\svgwidth{6cm}
\begingroup%
  \makeatletter%
  \providecommand\color[2][]{%
    \errmessage{(Inkscape) Color is used for the text in Inkscape, but the package 'color.sty' is not loaded}%
    \renewcommand\color[2][]{}%
  }%
  \providecommand\transparent[1]{%
    \errmessage{(Inkscape) Transparency is used (non-zero) for the text in Inkscape, but the package 'transparent.sty' is not loaded}%
    \renewcommand\transparent[1]{}%
  }%
  \providecommand\rotatebox[2]{#2}%
  \newcommand*\fsize{\dimexpr\f@size pt\relax}%
  \newcommand*\lineheight[1]{\fontsize{\fsize}{#1\fsize}\selectfont}%
  \ifx\svgwidth\undefined%
    \setlength{\unitlength}{576bp}%
    \ifx\svgscale\undefined%
      \relax%
    \else%
      \setlength{\unitlength}{\unitlength * \real{\svgscale}}%
    \fi%
  \else%
    \setlength{\unitlength}{\svgwidth}%
  \fi%
  \global\let\svgwidth\undefined%
  \global\let\svgscale\undefined%
  \makeatother%
  \begin{picture}(1,1)%
    \lineheight{1}%
    \setlength\tabcolsep{0pt}%
    \put(0,0){\includegraphics[width=\unitlength,page=1]{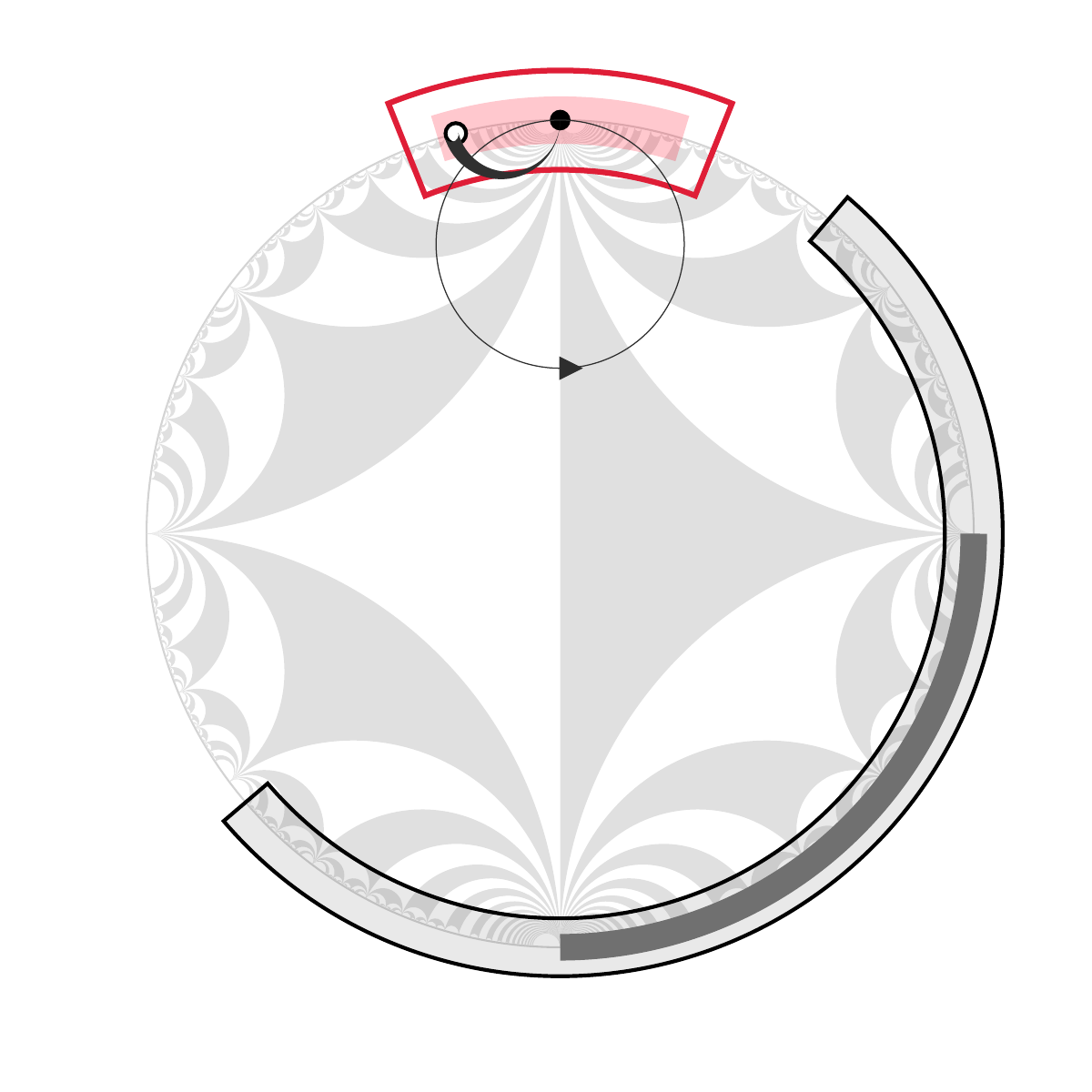}}%
    \put(0.47509135,0.60597285){\makebox(0,0)[lt]{\lineheight{1.25}\smash{\begin{tabular}[t]{l}$\gamma^{-1}$\end{tabular}}}}%
    \put(0.50160403,0.95821271){\makebox(0,0)[lt]{\lineheight{1.25}\smash{\begin{tabular}[t]{l}$q$\end{tabular}}}}%
    \put(0.39934082,0.93548756){\makebox(0,0)[lt]{\lineheight{1.25}\smash{\begin{tabular}[t]{l}$w$\end{tabular}}}}%
    \put(0.81596862,0.20828268){\makebox(0,0)[lt]{\lineheight{1.25}\smash{\begin{tabular}[t]{l}$K_q$\end{tabular}}}}%
  \end{picture}%
\endgroup%

  }{
    \def\svgwidth{6cm}
\begingroup%
  \makeatletter%
  \providecommand\color[2][]{%
    \errmessage{(Inkscape) Color is used for the text in Inkscape, but the package 'color.sty' is not loaded}%
    \renewcommand\color[2][]{}%
  }%
  \providecommand\transparent[1]{%
    \errmessage{(Inkscape) Transparency is used (non-zero) for the text in Inkscape, but the package 'transparent.sty' is not loaded}%
    \renewcommand\transparent[1]{}%
  }%
  \providecommand\rotatebox[2]{#2}%
  \newcommand*\fsize{\dimexpr\f@size pt\relax}%
  \newcommand*\lineheight[1]{\fontsize{\fsize}{#1\fsize}\selectfont}%
  \ifx\svgwidth\undefined%
    \setlength{\unitlength}{576bp}%
    \ifx\svgscale\undefined%
      \relax%
    \else%
      \setlength{\unitlength}{\unitlength * \real{\svgscale}}%
    \fi%
  \else%
    \setlength{\unitlength}{\svgwidth}%
  \fi%
  \global\let\svgwidth\undefined%
  \global\let\svgscale\undefined%
  \makeatother%
  \begin{picture}(1,1)%
    \lineheight{1}%
    \setlength\tabcolsep{0pt}%
    \put(0,0){\includegraphics[width=\unitlength,page=1]{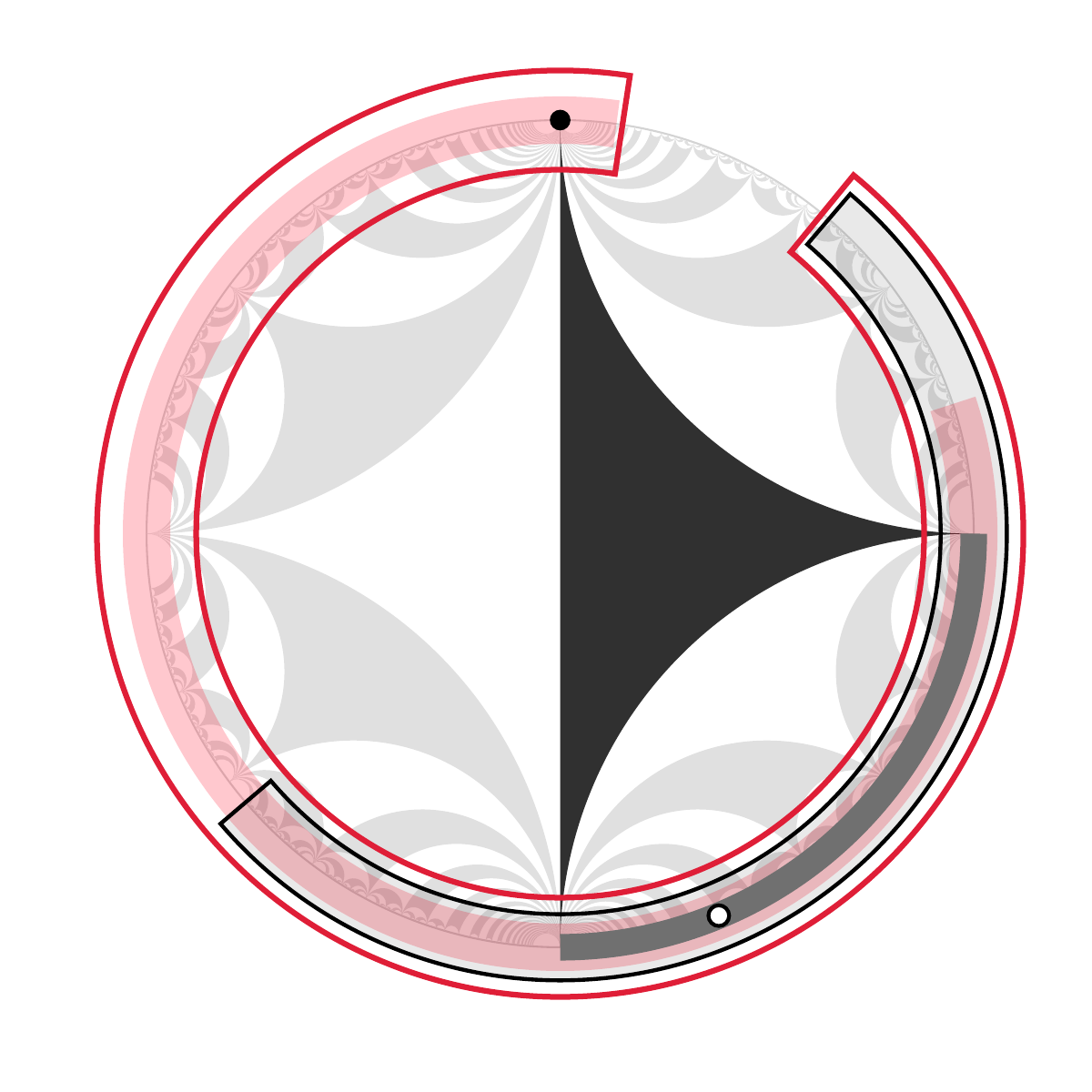}}%
    \put(0.6644676,0.0757193){\makebox(0,0)[lt]{\lineheight{1.25}\smash{\begin{tabular}[t]{l}$\gamma^{-1} w$\end{tabular}}}}%
  \end{picture}%
\endgroup%

  }
  \caption{For any point $w \ne q$ in a neighborhood $V_q$ (pink) of
    the parabolic point $q$, we find some $\gamma \in \Gamma_q$ so
    that $\gamma^{-1} w$ lies in $K_q$ (dark gray), a fundamental
    domain for the action of $\Gamma_q$ on $\dee \H^2 - \{q\}$. The
    expanded neighborhood $\gamma^{-1} W_q$ (red) contains a uniform
    neighborhood of $K_q$, so $\gamma^{-1}W_q$ contains a uniform
    neighborhood of $\gamma^{-1} w$.}
  \label{fig:sl2z_parabolic_expansion}
\end{figure}
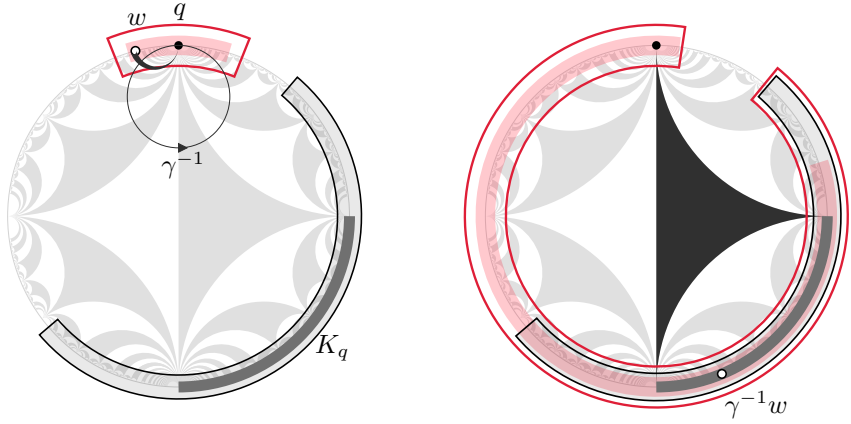

\begin{example}
  We consider the group $\PGL(2, \Z)$. While $\PGL(2, \Z)$ is
  virtually a free group (and therefore word-hyperbolic), it is also
  relatively hyperbolic, relative to the collection $\mc{H}$ of
  conjugates of the parabolic subgroup
  $\left\{\begin{pmatrix} \pm 1 & t\\0 & 1
    \end{pmatrix} : t \in \Z \right\}$.

  Since $\PGL(2, \Z)$ acts with finite covolume on the hyperbolic
  plane $\H^2$, the Bowditch boundary of the pair
  $(\PGL(2, \Z), \mc{H})$ is equivariantly identified with
  $\dee \H^2$, the visual boundary of $\H^2$. Given a non-parabolic
  point $w \in \dee \H^2$, we can find an element of $\PGL(2, \Z)$
  which ``expands'' a neighborhood of $w$. There are two distinct
  possibilities:
  \begin{enumerate}
  \item Suppose $w$ is in a small neighborhood $V_z$ of a conical
    limit point $z \in \dee \H^2$. Then choose some loxodromic element
    $\gamma \in \PGL(2, \Z)$ whose attracting fixed point is close to
    $z$. Then, if $W_z$ is a slightly larger neighborhood of $z$,
    $\gamma^{-1} \cdot W_z$ is large enough to contain a uniformly
    large neighborhood of $\gamma^{-1} \cdot w$. See
    \Cref{fig:sl2z_loxodromic_expansion}.
  \item On the other hand, suppose $w$ is in a small neighborhood
    $V_q$ of a parabolic fixed point $q \in \dee \H^2$, but $w \ne
    q$. We can find some element
    $\gamma \in \Gamma_q = \Stab_\Gamma(q)$ so that $\gamma^{-1}$
    takes $w$ into a fundamental domain for the action of $\Gamma_q$
    on $\dee \H^2 - \{q\}$. Then, if $W_q$ is a slightly larger
    neighborhood of $q$, $\gamma^{-1} \cdot W_q$ is again large enough
    to contain a uniformly large neighborhood of
    $\gamma^{-1} \cdot w$. See \Cref{fig:sl2z_parabolic_expansion}.
  \end{enumerate}
  There is a slight issue with this approach: in the second case above
  (when $w$ is close to a parabolic point $q$), it is actually not
  quite good enough to ``expand'' a neighborhood of $w$ by using
  $\Gamma_q$ to push $w$ into a fundamental domain for $\Gamma_q$ on
  $\dee \H^2 - \{q\}$. The reason is that there might not be a way a
  choose these fundamental domains so that they all lie uniformly far
  away from $q$, as $q$ varies over all of the (infinitely many)
  parabolic points in $\bgamh$. We resolve this issue by instead
  choosing $\gamma$ to lie in a \emph{coset} $g\Gamma_p$, where
  $q = gp$ for some $p \in \Pi$. Then $\gamma^{-1} \cdot w$ lies in a
  fundamental domain for $\Gamma_p$ on $\dee \H^2 - \{p\}$, which
  allows us to get uniform control on the size of the expanded
  neighborhood $\gamma^{-1} W_q$.
\end{example}

The two technical lemmas below (\Cref{lem:conical_limit_neighborhoods}
and \Cref{lem:large_parabolic_neighborhood}) essentially say that one
can set up this kind of expansion \emph{simultaneously} on the
Bowditch boundary of our relatively hyperbolic group $\Gamma$
\emph{and} in a neighborhood of the $\Gamma$-invariant set
$\Lambda \subset M$. The precise formulation of the expansion
condition found in these two lemmas is best motivated by the proof of
\Cref{prop:edge_inclusions} below, which shows that the ``expanding''
open sets we construct give rise to a $\mc{G}$-compatible system of
open sets on a $\Gamma$-graph $\mc{G}$.

We first need the following consequence of condition
\ref{item:strata_interior} on strata:
\begin{lem}
  \label{lem:eps_neighborhoods_separated}
  Given $D > 0$, there exists $\eps > 0$ so that for any
  $a, b \in \bgamh$ with $d(a,b) > D$, the $\eps$-neighborhood of
  $\phi^{-1}(a)$ in $M$ is contained in $C_b$.
\end{lem}
\begin{proof}
  Fix $\delta < D/4$. For each closed ball $B$ of radius $\delta$ in
  $\bgamh$, the set of points $z \in \bgamh$ satisfying
  $d(z, B) \ge D/2$ is a compact subset $K$ of $\bgamh$, so by
  condition \ref{item:strata_interior} there is some $\eps > 0$ so
  that an $\eps$-neighborhood of $\phi^{-1}(B)$ is contained in
  $\bigcap_{z \in K} C_z$. In particular for any $y \in B$, if
  $a \in \bgamh$ satisfies $d(y, a) > D$, then
  $d(a, B) > D - \mathrm{diam}(B) \ge D/2$ and so the
  $\eps$-neighborhood of $\phi^{-1}(y)$ is contained in $C_a$. Since
  $\bgamh$ is compact we can cover it with finitely many of these
  balls, and take the minimum required $\eps$.
\end{proof}

\begin{remark}
  This uniformity lemma is where we need condition
  \ref{item:strata_interior}, rather than the weaker condition
  \ref{item:weak_strata_interior} used in an earlier version of the
  paper.
\end{remark}

\begin{lem}
  \label{lem:conical_limit_neighborhoods}
  There exists $\eps_{\mr{con}} > 0, \delta_{\mr{con}} > 0$ satisfying
  the following: for any $\eps > 0$, $\delta > 0$ with
  $\eps < \eps_{\mr{con}}$, $\delta < \delta_{\mr{con}}$, and every
  conical limit point $z$, we can find:
  \begin{itemize}
  \item A group element $\gamma_z \in \Gamma$
  \item Open subsets $W_z, V_z \subset \bgamh$ with
    $z \in V_z \subset W_z$
  \end{itemize}
  such that:
  \begin{enumerate}
  \item $\mr{diam}(W_z) < \delta$,
  \item \label{item:conical_nbhds_nested} In $\bgamh$,
    we have
    \[
      N_{\dee\Gamma}(\gamma_z^{-1} V_z, \delta) \subset \gamma_z^{-1} W_z.
    \]
  \item \label{item:conical_eps_nbhds_nested} In $M$ we have
    \[
      N_M(\gamma_z^{-1}\phi^{-1}(W_z), 2\eps) \subset
      \gamma_z^{-1}N_M(\phi^{-1}(z), \eps).
    \]
  \end{enumerate}
\end{lem}

\begin{remark}
  Conditions (1) and (2) together imply that for any
  $y,z \in \conbdry\Gamma$, if $\gamma_z^{-1}V_z$ intersects $V_y$,
  then $\gamma_zW_y \subset W_z$. Later, we will see that condition
  (3) implies that if $\gamma_z^{-1}V_z$ intersects $V_y$, then also
  $\gamma_z N_M(\phi^{-1}(y), 2\eps) \subset N_M(\phi^{-1}(z), \eps)$
  (giving us the inclusion indicated by
  \Cref{fig:expansion_inclusions}).
\end{remark}

\begin{figure}[h]
  \centering
  \def\svgwidth{.8\textwidth}
  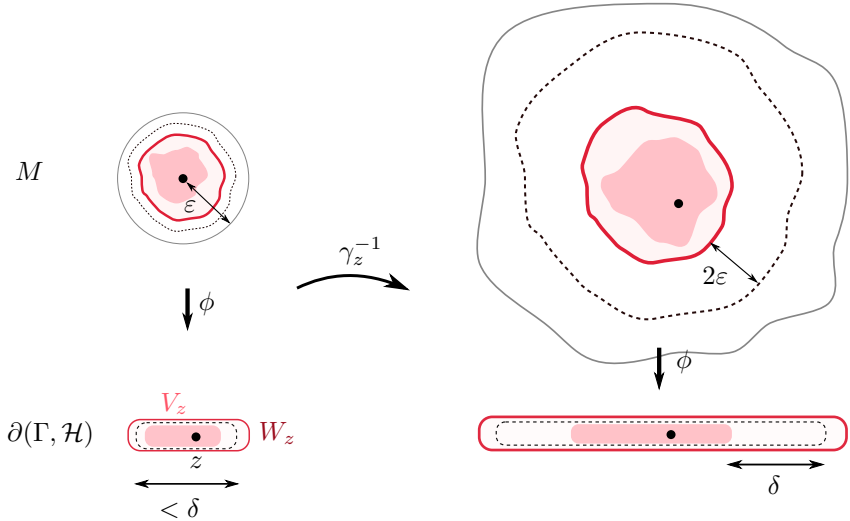
  \caption{The group element $\gamma_z^{-1}$ is ``expanding'' about
    $V_z \subset \bgamh$: while $W_z$ has diameter $< \delta$,
    $\gamma_z^{-1}W_z$ contains a $\delta$-neighborhood of
    $\gamma_z^{-1}V_z$. At the same time, $\gamma_z^{-1}$ enlarges an
    $\eps$-neighborhood of $\phi^{-1}(z)$ in $M$, so that it contains
    a $2\eps$-neighborhood of $\gamma_z^{-1}\phi^{-1}(W_z)$.}
  \label{fig:conical_convergence_setup}
\end{figure}
\begin{proof}
  For a conical limit point $z$, we choose a sequence $\gamma_n$ so
  that for distinct $a, b \in \bgamh$, we have $\gamma_n^{-1} z \to a$
  and $\gamma_n ^{-1} w \to b$ for any $w \ne z$. That is, $\gamma_n$
  limits conically to $z$ in $\overline{\Gamma}$, and $\gamma_n^{-1}$
  converges (not necessarily conically) to $b$. Since the
  $\Gamma$-action on distinct pairs in $\bgamh$ is cocompact
  (\Cref{prop:pairs_cocompact}), we may assume that $d(a, b) > D$ for
  a uniform constant $D > 0$.

  We choose $\eps_{\mr{con}} > 0$ from
  \Cref{lem:eps_neighborhoods_separated} so that if
  $a, b \in \bgamh$ satisfy $d(a, b) > D/2$, then a
  $2\eps_{\mr{con}}$-neighborhood of $\phi^{-1}(a)$ is contained in
  $C_b$. Let $\eps > 0$ satisfy $\eps < \eps_{\mr{con}}$, and let
  $\delta$ satisfy $\delta < \delta_{\mr{con}} := D / 4$.

  By the triangle inequality, we have $d(c, b) > D/2$ for all
  $c \in B_{\dee\Gamma}(a, 2\delta)$, so the closed
  $2\eps$-neighborhood of $\phi^{-1}(B_{\dee\Gamma}(a, 2\delta))$ is
  contained in $C_b$. This means that we can choose $n$ large enough
  so that
  \[
    \gamma_n \cdot N_M(\phi^{-1}(B_{\dee \Gamma}(a, 2\delta)), 2\eps)
  \]
  is contained in $N_M(\phi^{-1}(z), \eps)$ and
  \[
    \gamma_n \cdot B_{\dee\Gamma}(a, 2\delta)
  \]
  is contained in $B_{\dee\Gamma}(z, \delta / 2)$. We let
  $\gamma_z = \gamma_n$ for this large $n$, and take
  \[
    W_z = \gamma_z \cdot B_{\dee \Gamma}(a, 2\delta)
  \]
  and
  \[
    V_z = \gamma_z \cdot B_{\dee\Gamma}(a, \delta).
  \]
\end{proof}

The next lemma is a version of \Cref{lem:conical_limit_neighborhoods}
for parabolic points. As before, we want to show that for a point $q$
in the Bowditch boundary, we can find a neighborhood $W_q$ of $q$ in
$\bgamh$ with uniformly bounded diameter $\delta$, and group elements
$\gamma \in \Gamma$ so that $\gamma^{-1}$ enlarges $W_q$ enough to
contain a $2\delta$-neighborhood of $\gamma^{-1}z$, for some $z$ close
to $q$. Simultaneously we want to choose $\gamma$ so that
$\gamma^{-1}$ enlarges an $\eps$-neighborhood of $\phi^{-1}(q)$ in a
similar manner. This case is more complicated, because we need to
allow $\gamma$ to depend on the point $z \in W_q$: if $q$ is a
parabolic point in $\bgamh$, then in general there is \emph{not} a
single group element in $\Gamma$ which expands distances in a
neighborhood of $q$.

\begin{lem}
  \label{lem:large_parabolic_neighborhood}
  For each point $p \in \Pi$, there exist constants $\eps_p > 0$,
  $\delta_p > 0$ such that for any $q = g \cdot p \in \Gamma \cdot p$,
  any $\eps < \eps_p$, and any $\delta < \delta_p$, we can find:
  \begin{itemize}
  \item A cofinite subset $T_q$ of the coset $g\Gamma_p$,
  \item Open neighborhoods $V_q, W_q$ of $\bgamh$, with
    $q \in V_q \subset W_q$,
  \item Open neighborhoods $\hat{V}_q, \hat{W}_q$ of $\bgamh$
    with $\hat{V}_q \subset \hat{W}_q$
  \end{itemize}
  such that:
  \begin{enumerate}
  \item $\mr{diam}(W_q) < \delta$, and
    $\phi^{-1}(W_q) \subset N_M(\phi^{-1}(q), \eps)$.
  \item \label{item:parabolic_nbhds_nested} in $\bgamh$, we have
    \[
      N_{\dee\Gamma}(\hat{V}_q, \delta) \subset \hat{W}_q.
    \]
  \item \label{item:parabolic_backward_edge} For every
    $z \in V_q - \{q\}$, there exists $\gamma \in T_q$ with
    $\gamma^{-1} \cdot z \in \hat{V}_q$.
  \item \label{item:parabolic_eps_nbhds_nested} For every
    $\gamma \in T_q$, we have
    \[
      N_M(\phi^{-1}(\hat{W}_q), 2\eps) \subset \gamma^{-1}
      N_M(\phi^{-1}(q), \eps)
    \]
    and
    \[
      \hat{W}_q \subset \gamma^{-1} W_q.
    \]
  \item \label{item:parabolic_in_attracting}
    $N_M(\phi^{-1}(\hat{W}_q), 2\eps)$ is contained in $C_p$ and
    $g\Gamma_p \cdot \hat{V}_q$ contains $\bgamh - \{q\}$.
  \end{enumerate}
\end{lem}

\begin{figure}[h]
  \centering
  \def\svgwidth{3.5in}
  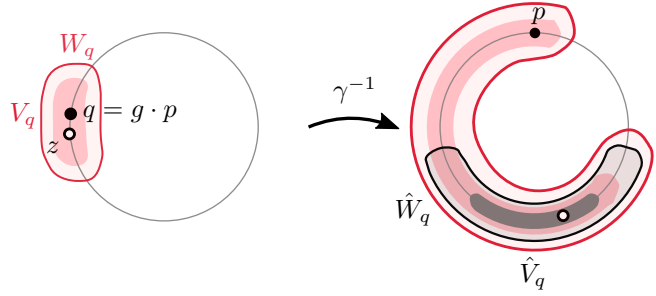

  \caption{The behavior of sets in $\bgamh$ described by
    \Cref{lem:large_parabolic_neighborhood}. Given $z \in V_q$, we
    pick an element $\gamma \in g\Gamma_p$ so that a uniformly large
    neighborhood of $\gamma^{-1}z$ is contained in
    $\gamma^{-1}W_q$. We cannot pick $\gamma^{-1}$ to expand the
    metric everywhere close to $q$---some points in $V_q$ get
    contracted close to $p$.}
  \label{fig:parabolic_expansion}
\end{figure}
\begin{remark}
  If $z \in V_q - \{q\}$ and $\gamma^{-1} z \in \hat{V}_q$ for some
  $\gamma \in T_q$, we think of $\gamma^{-1}$ as ``expanding'' about
  $z$. Conditions (1) and (2) imply that if $\gamma^{-1}z \in V_y$ for
  some $\gamma \in T_q$, then $\hat{W}_q$ contains $W_y$, and by
  condition (4), $\gamma^{-1}W_q$ contains $W_y$. Here $V_y, W_y$ are
  the sets from either \Cref{lem:conical_limit_neighborhoods} or
  \Cref{lem:large_parabolic_neighborhood}.
\end{remark}
\begin{proof}
  Pick a compact set $K \subset \bgamh - \{p\}$ so that
  $\Gamma_p \cdot K$ covers $\bgamh - \{p\}$. Choose
  $\delta_p$ small enough so that the closure of
  $N_{\dee\Gamma}(K, 2\delta_p)$ does not contain $p$. Then, for any
  $\delta < \delta_p$, we can pick
  \[
    \hat{V}_q = N_{\dee\Gamma}(K, \delta), \quad \hat{W}_q = N_{\dee
      \Gamma}(K, 2\delta).
  \]

  We can choose $\eps_p$ sufficiently small so that a
  $2\eps_p$-neighborhood of $\phi^{-1}(N_{\dee\Gamma}(K, 2\delta_p))$
  is contained in $C_p$. Now, fix $\eps < \eps_p$. We claim that there
  exists a cofinite subset $T_q \subset g \cdot \Gamma_p$ so that for
  any $\gamma \in T_q$, we have
  \begin{eqnarray}
    \label{eq:bdry_coset_containment}
    \gamma \cdot \hat{W}_q \subset& B_{\dee\Gamma}(q, \delta / 2)\\
    \label{eq:coset_containment}
    \gamma \cdot N_M(\phi^{-1}(\hat{W}_q), 2\eps) \subset& N_M(\phi^{-1}(q),
                                                           \eps)
  \end{eqnarray}
  To see that this claim holds, it suffices to verify that for any
  infinite sequence $\gamma_n$ of distinct group elements in
  $g \Gamma_p$, (\ref{eq:bdry_coset_containment}) and
  (\ref{eq:coset_containment}) both hold for all sufficiently large
  $n$.

  We write $\gamma_n = g \cdot \gamma_n'$ for
  $\gamma_n' \in \Gamma_p$. Then $\gamma_n'$ converges uniformly to
  $p$ on compact subsets of $\bgamh - \{p\}$, so $\gamma_n$
  converges uniformly to $q$ on compact subsets of
  $\bgamh - \{p\}$, implying that
  (\ref{eq:bdry_coset_containment}) eventually holds. And by our
  assumptions, we know that
  \[
    \gamma_n' \cdot N_M(\phi^{-1}(\hat{W}_q), 2\eps) \subset g^{-1}
    \cdot N_M(\phi^{-1}(q), \eps)
  \]
  for sufficiently large $n$, implying that
  (\ref{eq:coset_containment}) also eventually holds.

  So we can take $W_q$ to be the set
  \[
    \{q\} \cup \bigcup_{\gamma \in T_q} \gamma \cdot \hat{W}_q,
  \]
  and $V_q$ to be the set
  \[
    \{q\} \cup \bigcup_{\gamma \in T_q} \gamma \cdot \hat{V}_q.
  \]

  To see that $W_q$ and $V_q$ are open we just need to verify that
  they each contain a neighborhood of $q$. Since $\hat{V}_q$ and
  $\hat{W}_q$ each contain $K$, and $\Gamma_p \cdot K$ covers
  $\bgamh - \{p\}$, $V_q$ and $W_q$ each contain the set
  \[
    \bgamh - \bigcup_{\gamma \in g\Gamma_p - T_q} \gamma K.
  \]
  Since $T_q$ is cofinite in $g\Gamma_p$ this is an open set
  containing $q$.
\end{proof}

\subsection{Construction of the relative automaton}

We will construct the relative automaton $\mc{G}$ satisfying the
conditions of \Cref{prop:convergence_dynamics_implies_automaton} by
choosing a suitable open covering of $\bgamh$, and then using
compactness to take a finite subcover. The subsets of this subcover
will be the vertices of $\mc{G}$.

We choose constants $\eps > 0$, $\delta > 0$ so that
$\eps < \eps_{\mr{con}}$, $\delta < \delta_{\mr{con}}$ (where
$\eps_{\mr{con}}$, $\delta_{\mr{con}}$ are the constants coming from
\Cref{lem:conical_limit_neighborhoods}) and $\eps < \eps_p$,
$\delta < \delta_p$ for each $p \in \Pi$ (where $\eps_p, \delta_p$ are
the constants coming from \Cref{lem:large_parabolic_neighborhood}).

Then:
\begin{itemize}
\item For each $z \in \conbdry\Gamma$, we define $W_z$, $V_z$,
  $\gamma_z$ as in \Cref{lem:conical_limit_neighborhoods}, with
  parameters $\eps$, $\delta$.
\item For each $q \in \parbdry \Gamma$, we define
  $V_q, W_q, \hat{V}_q, \hat{W}_q$, and $T_q$ as in
  \Cref{lem:large_parabolic_neighborhood}, again with parameters
  $\eps$, $\delta$.
\end{itemize}
The collections of sets $\{V_z : z \in \conbdry \Gamma\}$ and
$\{V_q : q \in \parbdry\Gamma\}$ together give an open covering of
the Bowditch boundary $\bgamh$. So we choose a finite subcover
$\mc{V}$, which we can write as
\[
  \mc{V} = \{V_a : a \in A\}
\]
where $A$ is a finite subset of $\bgamh$. We can in
particular ensure that $A$ contains the finite set $\Pi$.

We identify the vertices of our $\Gamma$-graph $\mc{G}$ with $A$. For
each $a \in A$, the set $T_a$ is either $\{\gamma_a\}$ (if $a$ is a
conical limit point) or $T_q$ (if $a = q$ for a parabolic point
$q$). Then, for each $a \in A$, we define the open sets $U_a$ by
\[
  U_a = N_M(\phi^{-1}(a), \eps).
\]

The edges of the $\Gamma$-graph $\mc{G}$ are defined as follows:
\begin{itemize}
\item For $a, b \in A$ with $a \in \conbdry\Gamma$, there is an edge
  from $a$ to $b$ if $(\gamma_a^{-1} \cdot V_a) \cap V_b$ is nonempty.
\item If $a, b \in A$ with $a \in \parbdry\Gamma$, there is an edge
  from $a$ to $b$ if $\hat{V}_a \cap V_b$ is nonempty.
\end{itemize}

Note that every vertex $a$ of $\mc{G}$ has at least one outgoing edge:
if $a$ is conical, then $\gamma_a^{-1}V_a$ is nonempty, and since
$\mc{V}$ is a covering of $\bgamh$, the set $\gamma_a^{-1}V_a$ has
nonempty intersection with at least one $V_b \in \mc{V}$. The exact
same reasoning applies if $a$ is parabolic, since $\hat{V}_a$ is
nonempty. Moreover, for any parabolic point $a$, the set $\hat{V}_a$
depends only on the orbit of $a$ in $\bgamh$, so $\mc{G}$ must satisfy
condition (2) in \Cref{prop:convergence_dynamics_implies_automaton}.
\begin{prop}
  \label{lem:small_nbhds_in_lg_nbhds}
  For each $a \in A$, we have
  \[
    \phi^{-1}(W_a) \subset U_a.
  \]
\end{prop}
\begin{proof}
  When $a$ is not a parabolic vertex, Part (3) of
  \Cref{lem:conical_limit_neighborhoods} implies:
  \[
    \phi^{-1}(W_a) = \gamma_a\gamma_a^{-1}\phi^{-1}(W_a) \subset
    \gamma_aN_M(\gamma_a^{-1}\phi^{-1}(W_a), 2\eps) \subset
    N_M(\phi^{-1}(a), \eps) = U_a.
  \]
  When $a$ is a parabolic vertex, then the claim follows directly from
  Part (1) of \Cref{lem:large_parabolic_neighborhood}.
\end{proof}

Next we verify:

\begin{prop}
  \label{prop:edge_inclusions}
  The collection of sets $\{W_v\}$ and $\{U_v\}$ are both
  $\mc{G}$-compatible systems of open sets for the $\Gamma$-graph
  $\mc{G}$.
\end{prop}
\begin{proof}
  First fix an edge $(a, b)$ with $a \in \conbdry \Gamma$. Since
  $(\gamma_a^{-1}V_a) \cap V_b$ is nonempty, part
  \ref{item:conical_nbhds_nested} of
  \Cref{lem:conical_limit_neighborhoods} implies that
  $\gamma_a^{-1} \cdot W_a$ contains the $\delta$-neighborhood of some
  point $z \in V_b$. Since $\mr{diam}(W_b) < \delta$ and
  $V_b \subset W_b$, we can find a small $\eps' > 0$ so that
  $\gamma_aN_{\dee\Gamma}(W_b, \eps') \subset W_a$.

  In particular, $\gamma_a^{-1} \cdot W_a$ contains $b$, which means
  that $N_M(\gamma_a^{-1}\phi^{-1}(W_a), 2\eps)$ contains
  $N_M(\phi^{-1}(b), 2\eps)$, which contains $N_M(U_b, \eps)$. Then,
  part \ref{item:conical_eps_nbhds_nested} of
  \Cref{lem:conical_limit_neighborhoods} implies that
  $\gamma_a \cdot N_M(U_b, \eps)$ is contained in
  $N_M(\phi^{-1}(a), \eps) = U_a$.

  Next fix an edge $(q, b)$ with $q \in \parbdry \Gamma$. From part
  \ref{item:parabolic_nbhds_nested} of
  \Cref{lem:large_parabolic_neighborhood}, we know that $\hat{W}_q$
  contains the $\delta$-neighborhood of a point
  $z \in \hat{V}_q \cap V_b$. Since $\mr{diam}(W_b) < \delta$ and
  $V_b \subset W_b$, this means that $\hat{W}_q$ contains an
  $\eps'$-neighborhood of $W_b$ for some small $\eps' > 0$. So part
  \ref{item:parabolic_eps_nbhds_nested} of
  \Cref{lem:large_parabolic_neighborhood} implies that for any
  $\gamma \in T_q$, we have $\gamma \cdot N(W_b, \eps') \subset W_q$.

  In particular $\hat{W}_q$ contains $b$, so
  $N_M(\phi^{-1}(\hat{W}_q), 2\eps)$ contains
  $N_M(\phi^{-1}(b), 2\eps)$, which contains $N_M(U_b, \eps)$. Then,
  part \ref{item:parabolic_eps_nbhds_nested} of
  \Cref{lem:large_parabolic_neighborhood} implies that
  \[
    \gamma N_M(U_b, \eps) \subset N_M(\phi^{-1}(q), \eps) = U_q
  \]
  for any $\gamma \in T_q$.
\end{proof}

We observe:
\begin{prop}
  The $\mc{G}$-compatible systems of open subsets $\{U_v\}$ and
  $\{W_v\}$ satisfy conditions (1) - (3) in
  \Cref{prop:convergence_dynamics_implies_automaton}.
\end{prop}
\begin{proof}
  Part (1) follows directly from \Cref{lem:small_nbhds_in_lg_nbhds},
  and the fact that we defined each $U_a$ to be the
  $\eps$-neighborhood of $\phi^{-1}(a)$. Part (2) is true by the
  construction of the open covering $\mc{V}$ and the graph
  $\mc{G}$. Part (3) is true by construction and part
  \ref{item:parabolic_in_attracting} of
  \Cref{lem:large_parabolic_neighborhood}.
\end{proof}

To finish the proof of
\Cref{prop:convergence_dynamics_implies_automaton}, we now just need
to show:
\begin{prop}
  \label{prop:graph_is_automaton}
  The $\Gamma$-graph $\mc{G}$ is a relative quasigeodesic automaton
  for the pair $(\Gamma, \mc{H})$.
\end{prop}
\begin{proof}
  We apply \Cref{prop:relative_geodesic_automaton_from_bb}, using the
  $\mc{G}$-compatible system $\{W_a\}$ and the sets $\{V_a\}$ we
  defined in the construction of $\mc{G}$.

  The first three conditions of
  \Cref{prop:relative_geodesic_automaton_from_bb} are satisfied by
  construction. To see that conditions \ref{item:conical_edge_exists}
  and \ref{item:parabolic_edge_exists} hold, first observe that if
  $z \in V_a$ for a non-parabolic vertex $a$, then
  $\gamma_a^{-1} \cdot z$ lies in some $V_b$ and $(a,b)$ is an edge in
  $\mc{G}$. And if $z \in V_a - \{p_a\}$ for a parabolic vertex $a$,
  then part (\ref{item:parabolic_backward_edge}) of
  \Cref{lem:large_parabolic_neighborhood} says that there is some
  $\gamma \in T_a$ such that $\gamma^{-1} \cdot z \in \hat{V}_a$. If
  $V_b$ contains $\gamma^{-1} \cdot z$, the edge $(a,b)$ must be in
  $\mc{G}$.

  It only remains to check that $\mc{G}$ is a divergent
  $\Gamma$-graph. Let $\{\alpha_n\}$ be an infinite $\mc{G}$-path,
  following a vertex path $\{v_n\}$. The $\mc{G}$-compatibility
  condition implies that $\gamma_n \overline{U}_{v_{n+1}}$ is a subset
  of $\gamma_{n-1} U_{v_n}$ for every $n$. Since $M$ is connected and
  each $U_v$ is a proper subset of $M$, this inclusion must be
  proper. This implies that in the sequence $\gamma_n$, no element can
  appear more than $\#A$ times and therefore $\gamma_n$ is divergent.
\end{proof}

\begin{remark}
  This last step is the only part of the proof of
  \Cref{prop:convergence_dynamics_implies_automaton} which uses the
  connectedness of $M$. Note that the proof does \emph{not} assume
  that the Bowditch boundary $\bgamh$ itself is connected, and this
  does not follow from connectedness of $M$ since only the subset
  $\Lambda \subset M$ surjects onto $\bgamh$.

  Thus, the above arguments almost (but do not quite) show the general
  statement that \emph{any relatively hyperbolic group has a relative
    quasigeodesic automaton}. This statement was later proved in
  \cite{MMW2}, using arguments very similar to the above. It is likely
  possible to modify those arguments to also prove
  \Cref{prop:convergence_dynamics_implies_automaton} without assuming
  that $M$ is connected. However, such a proof would involve
  introducing additional technicalities in the construction of the
  sets $V_a$ and $W_a$, so we elect not to do so here.
\end{remark}

We conclude this section by observing that one can slightly refine the
construction in \Cref{prop:convergence_dynamics_implies_automaton} to
obtain some stronger conditions on the resulting automaton.
\begin{prop}
  \label{prop:compact_bowditch_nbhd}
  Fix a compact subset $Z$ of the Bowditch boundary $\bgamh$. Then,
  for any open set $U \subset M$ containing $\phi^{-1}(Z)$, there is a
  relative quasigeodesic automaton $\mc{G}$ and a pair of
  $\mc{G}$-compatible systems of open sets $\{U_a\}$, $\{W_a\}$ as in
  \Cref{prop:convergence_dynamics_implies_automaton}, additionally
  satisfying the following: any $z \in Z$ is the limit of a
  $\mc{G}$-path $\{\alpha_n\}$ (with corresponding vertex path
  $\{v_n\}$) such that $U_{v_1} \subset U$.
\end{prop}
\begin{proof}
  We choose $\eps > 0$ so that $U$ contains $N_M(\phi^{-1}(Z),
  \eps)$. We then construct our relative quasigeodesic automaton
  $\mc{G}$ as in the proof of
  \Cref{prop:convergence_dynamics_implies_automaton}, but we also
  choose a finite subset $A_Z \subset Z$ so that the sets $V_a$ for
  $a \in A_Z$ give a finite open covering of $Z$. We can ensure that
  the vertex set $A$ of $\mc{G}$ contains $A_Z$.
  
  Then, for any $z \in Z$, by the construction in
  \Cref{prop:relative_geodesic_automaton_from_bb}, we can find a
  $\mc{G}$-path limiting to $z$ whose first vertex is some
  $a \in A_Z$. The corresponding open set for this vertex is
  $U_a = N_M(\phi^{-1}(a), \eps) \subset U$.
\end{proof}

\begin{prop}
  \label{prop:parabolic_outward_edge_nbhds}
  For each parabolic point $p \in \Pi$, let $K_p$ be a compact subset
  of $\bgamh - \{p\}$ such that $\Gamma_p \cdot K_p = \bgamh -
  \{p\}$. Then the relative quasigeodesic automaton in
  \Cref{prop:convergence_dynamics_implies_automaton} can be chosen to
  satisfy the following:

  For every parabolic vertex $w$ with $p_w = p \in \Pi$, and
  every $z \in K_p$, there is a $\mc{G}$-path limiting to $z$ whose
  first vertex $u$ is connected to $w$ by an edge $(w, u)$.
\end{prop}
\begin{proof}
  The proof of \Cref{lem:large_parabolic_neighborhood} shows that in
  our construction of the relative automaton, we can ensure that each
  set $\hat{V}_p$ contains $K_p$. So if $z \in K_p$, then by
  definition of the automaton, $z$ lies in $V_u$ with $w$ connected to
  $u$ by a directed edge. Then, following the proof of
  \Cref{prop:relative_geodesic_automaton_from_bb}, we can find a
  $\mc{G}$-path limiting to $z$ whose first vertex is $u$.
\end{proof}


\section{Contracting paths in flag manifolds}
\label{sec:contracting_paths}
\label{sec:flag_convergence}

Let $\Gamma \subset G$ be a discrete relatively hyperbolic group, and
let $\mc{G}$ be a $\Gamma$-graph. Fix a pair of opposite parabolic
subgroups $\pargroup$, $\oppgroup$. Our goal in this section is to
show that under certain conditions, if $\{U_v\}$ is a
$\mc{G}$-compatible system of open subsets of $\pflags$ for the action
of $\Gamma$ on $\pflags$, then the sequence of group elements lying
along an infinite $\mc{G}$-path is $\pargroup$-divergent.

\subsection{Contracting paths in $\Gamma$-graphs}

\begin{definition}
  Let $\Gamma$ be a discrete subgroup of $G$, let $\mc{G}$ be a
  $\Gamma$-graph, and let $\{U_v\}_{v \in V(\mc{G})}$ be a
  $\mc{G}$-compatible system of open subsets of $\pflags$. We say that
  a $\mc{G}$-path $\{\alpha_n\}_{n \in \N}$ is \emph{contracting} if
  the decreasing intersection
  \begin{equation}
    \label{eq:contracting_intersection}
    \bigcap_{n=1}^\infty \alpha_1 \cdots \alpha_n \cdot U_{v_{n+1}}
  \end{equation}
  is a singleton in $\pflags$.
\end{definition}

\begin{definition}
  We say that an open set $\Omega \subset \pflags$ is a \emph{proper
    domain} if the closure of $\Omega$ lies in an affine chart
  $\Opp(\xi) \subset \pflags$ for some $\xi \in \oppflags$.
\end{definition}

Here is the main result in this section:
\begin{prop}
  \label{prop:compatible_system_in_flags_contracting}
  Let $\mc{G}$ be a $\Gamma$-graph for $(\Gamma, \mc{H})$, and let
  $\{U_v\}_{v \in V(\mc{G})}$ be $\mc{G}$-compatible system of open
  subsets of $\pflags$.

  If the set $U_v$ is a proper domain for each vertex $v$ of the
  automaton, then every infinite $\mc{G}$-path is contracting.
\end{prop}

\subsection{A metric property for proper domains in flag manifolds}

To prove \Cref{prop:compatible_system_in_flags_contracting}, we
consider a metric $C_\Omega$ defined by Zimmer \cite{zimmer2018proper}
on any proper domain $\Omega \subset \pflags$. $C_\Omega$ is defined
so that it is invariant under the action of $G$ on $\pflags$: for any
$x,y$ in some proper domain $\Omega \subset \pflags$, and any
$g \in G$, we have
\begin{equation}
  \label{eq:g_invariance}
  C_\Omega(x,y) = C_{g\Omega}(gx, gy).
\end{equation}

In general, $C_\Omega$ is not a complete metric. However, $C_\Omega$
induces the standard topology on $\Omega$ as an open subset of
$\pflags$. We will show that for a $\mc{G}$-path $\{\alpha_n\}$, the
diameter of
\[
  \alpha_1 \cdots \alpha_n U_{v_{n+1}}
\]
with respect to $C_{U_{v_1}}$ tends to zero as $n \to \infty$.

Zimmer's construction of $C_\Omega$ depends on an irreducible
representation $\zeta:G \to \PGL(V)$ for some real vector space
$V$. This is provided by a theorem of
Gu\'eritaud-Guichard-Kassel-Wienhard.
\begin{theorem}[\cite{ggkw2017anosov}, see also
  \cite{zimmer2018proper}, Theorem 4.6]
  \label{thm:flag_irrep_embedding}
  There exists a real vector space $V$, an irreducible representation
  $\zeta:G \to \PGL(V)$, a line $\ell \subset V$, and a hyperplane
  $H \subset V$ such that:
  \begin{enumerate}
  \item $\ell + H = V$.
  \item The stabilizer of $\ell$ in $G$ is $\pargroup$ and the
    stabilizer of $H$ in $G$ is $\oppgroup$.
  \item $g\pargroup g^{-1}$ and $h\oppgroup h^{-1}$ are opposite if
    and only if $\zeta(g)\ell$ and $\zeta(h)H$ are transverse.
  \end{enumerate}
\end{theorem}

The representation $\zeta$ determines a pair of embeddings
$\iota:\pflags \to \P(V)$ and $\iota^*:\oppflags \to \P(V^*)$ by
\[
  \iota(g\pargroup) = \zeta(g)\ell, \qquad \iota^*(g\oppgroup) =
  \zeta(g)H.
\]
In this section, we will identify $\P(V^*)$ with the space of
projective hyperplanes in $\P(V)$, by identifying the projectivization
of a functional $w \in V^*$ with the projectivization of its kernel.

\begin{definition}
  Let $\Omega$ be an open subset of $\pflags$. The \emph{dual domain}
  $\Omega^* \subset \oppflags$ is
\[
  \Omega^* = \{\nu \in \oppflags : \nu \textrm{ is opposite to }\xi
  \textrm{ for every } \xi \in \overline{\Omega}\}.
\]
Note that $\Omega^*$ is always open, and it is nonempty precisely when
$\Omega$ is a proper domain.
\end{definition}

\begin{definition}
  Let $w_1, w_2 \in \P(V^*)$, and let $z_1, z_2 \in \P(V)$. The
  \emph{cross-ratio} $[w_1, w_2 ; z_1, z_2]$ is defined by
  \[
    \frac{\tilde{w}_1(\tilde{z}_2)\tilde{w}_2(\tilde{z}_1)}
    {\tilde{w}_1(\tilde{z}_1)\tilde{w}_2(\tilde{z}_2)},
  \]
  where $\tilde{w}_i$, $\tilde{z}_i$ are respectively lifts of $w_i$
  and $z_i$ in $V^*$ and $V$.
\end{definition}

\begin{remark}
  When $V$ is two-dimensional, we can identify the projective line
  $\P(V^*)$ with $\P(V)$ by identifying each $[w] \in \P(V^*)$ with
  $[\ker(w)] \in \P(V)$. In that case, the cross-ratio defined above
  agrees with the standard four-point cross-ratio on $\RP^1$, given by
  \begin{equation}
    \label{eq:cross_ratio_definition}
    [a, b; c, d] := \frac{(d - a)(c - b)}{(c - a)(d - b)}.
  \end{equation}
  The differences in (\ref{eq:cross_ratio_definition}) can be measured
  in any affine chart in $\RP^1$ containing $a, b, c, d$. Our
  convention is chosen so that if we identify $\RP^1$ with
  $\R \cup \{\infty\}$, we have $[0, \infty; 1, z] = z$.
\end{remark}

\begin{definition}
  Let $\Omega \subset \pflags$ be a proper domain. We define the
  function $C_\Omega:\Omega \times \Omega \to \R$ by
  \[
    C_\Omega(x,y) = \sup_{\xi_1, \xi_2 \in
      \Omega^*}\log|[\iota^*(\xi_1), \iota^*(\xi_2); \iota(x),
    \iota(y)]|.
  \]
\end{definition}
For any $g \in G$ and any proper domain $\Omega \subset \pflags$, we
have $(g\Omega)^* = (g\Omega^*)$. So $C_\Omega$ must satisfy the
$G$-invariance condition (\ref{eq:g_invariance}).

If $\Omega$ is a properly convex subset of $\P(V)$, and $\zeta$,
$\iota$, $\iota^*$ are the identity maps on $\PGL(V)$, $\P(V)$, and
$\P(V^*)$ respectively, then $C_\Omega$ agrees with the well-studied
\emph{Hilbert metric} on $\Omega$. More generally we have:

\begin{theorem}[\cite{zimmer2018proper}, Theorem 5.2]
  If $\Omega$ is open and bounded in an affine chart, then $C_\Omega$
  is a metric on $\Omega$ which induces the standard topology on
  $\Omega$ as an open subset of $\pflags$.
\end{theorem}

\begin{remark}
  This particular result in \cite{zimmer2018proper} is stated only for
  noncompact simple Lie groups, but the proof only assumes that $G$ is
  semisimple with no compact factor.
\end{remark}

Since taking duals of proper domains reverses inclusions, it follows
that if $\Omega_1 \subset \Omega_2$, then
$C_{\Omega_1} \ge C_{\Omega_2}$.  Our goal now is to sharpen this
inequality, and show:
\begin{prop}
  \label{prop:caratheodory_metric_shrinks}
  Let $\Omega_1$, $\Omega_2$ be proper domains in $\pflags$, such that
  $\overline{\Omega_1} \subset \Omega_2$.

  There exists a constant $\lambda > 1$ (depending on $\Omega_1$ and
  $\Omega_2$) so that for all $x,y \in \Omega_1$,
  \[
    C_{\Omega_1}(x, y) \ge \lambda \cdot C_{\Omega_2}(x,y).
  \]
\end{prop}

A consequence is the following, which in particular implies
\Cref{prop:compatible_system_in_flags_contracting}.
\begin{cor}
  \label{cor:exponential_path_shrinking}
  Let $\mc{G}$ be a $\Gamma$-graph for a relatively hyperbolic group
  $\Gamma$, and let $\{U_v\}$ be a $\mc{G}$-compatible system of open
  subsets of $\pflags$. If each $U_v$ is a proper domain, then there
  are constants $\lambda_1, \lambda_2 > 0$ so that for any
  $\mc{G}$-path $\{\alpha_n\}$ in the $\Gamma$-graph $\mc{G}$, the
  diameter of
  \[
    \alpha_1 \cdots \alpha_n \cdot U_{v_{n+1}}
  \]
  with respect to $C_{U_{v_1}}$ is at most
  \[
    \lambda_1 \cdot \exp(-\lambda_2 \cdot n).
  \]
\end{cor}
\begin{proof}
  For any open set $U \subset \pflags$ and $A \subset U$, we let
  $\mr{diam}_U(A)$ denote the diameter of $A$ with respect to the
  metric $C_U$. We choose a uniform $\eps > 0$ so that in some fixed
  metric on $\pflags$, every edge $(v, w)$ in $\mc{G}$, and every
  $\alpha \in T_v$, we have
  \[
    \alpha N(U_w, \eps) \subset U_v.
  \]
  Then for each vertex set $U_v$, we write $U_v^\eps = N(U_v, \eps)$.

  We take
  \[
    \lambda_1 = \max\{\mr{diam}_{U_v^\eps}(U_v)\}.
  \]
  \Cref{prop:caratheodory_metric_shrinks} implies that there exists
  $\lambda_v > 0$ such that for all $x, y \in U_v$, we have
  \[
    C_{U_v}(x,y) \ge \exp(\lambda_v) \cdot C_{U_v^\eps}(x, y).
  \]
  Take $\lambda_2 = \min_{v} \{\lambda_v\}$. We claim that for all
  $n \ge 1$, we have
  \[
    \mr{diam}_{U_{v_1}^\eps}(\alpha_1 \cdots \alpha_nU_{v_{n+1}}) \le \lambda_1
    \exp(-\lambda_2 \cdot (n - 1)).
  \]
  We prove the claim via induction on the length of the $\mc{G}$-path
  $\{\alpha_n\}$. For $n = 1$, the claim is true because
  $\alpha_1U_{v_2} \subset U_{v_1}$. For $n > 1$, we can assume
  \[
    \lambda_1 \exp(-\lambda_2(n - 2)) \ge
    \mr{diam}_{U_{v_2}^\eps}(\alpha_2 \cdots \alpha_n \cdot
    U_{v_{n+1}}).
  \]
  Then we have
  \begin{align*}
    \mr{diam}_{U_{v_2}^\eps}(\alpha_2 \cdots \alpha_n \cdot
    U_{v_{n+1}}) &= \mr{diam}_{\alpha_1U_{v_2}^\eps}(\alpha_1 \cdots \alpha_n \cdot
                   U_{v_{n+1}})\\
                 &\ge \mr{diam}_{U_{v_1}}(\alpha_1 \cdots \alpha_n \cdot
                   U_{v_{n+1}})\\
                 &\ge \exp(\lambda_2) \cdot
                   \mr{diam}_{U_{v_1}^\eps}(\alpha_1 \cdots \alpha_n
                   \cdot U_{v_{n+1}}).
  \end{align*}
  Finally, the claim implies the corollary because we know that
  \begin{align*}
    \mr{diam}_{U_{v_1}}(\alpha_1 \cdots \alpha_nU_{v_{n+1}})
    &\le \mr{diam}_{\alpha_1 U_{v_2}^\eps}(\alpha_1 \cdots \alpha_nU_{v_{n+1}}) \\
    &= \mr{diam}_{U_{v_2}^\eps}(\alpha_2 \cdots \alpha_nU_{v_{n+1}})\\
    &\le \lambda_1 \exp(-\lambda_2(n - 2)).
  \end{align*}
  So, we can replace $\lambda_1$ with $\lambda_1\exp(2\lambda_2)$ to
  get the desired result.
\end{proof}

We now proceed with the proof of
\Cref{prop:caratheodory_metric_shrinks}. We first observe that in the
special case where $\Omega_1, \Omega_2$ are properly convex subsets of
real projective space, one can show the desired result essentially via
the following:
\begin{prop}
  \label{prop:hilbert_metric_contracts}
  Let $a, b, c, d$ be points in $\RP^1$, arranged so that
  $a < b < c < d \le a$ in a cyclic ordering on $\RP^1$. Then there
  exists a constant $\lambda > 1$, depending only on the cross-ratio
  $[a, b; c, d]$, so that for all distinct $x, y \in (b, c)$, we have
  \[
    |\log[b,c; x,y]| \ge \lambda \cdot |\log[a,d; x, y]|.
  \]
\end{prop}
\Cref{prop:hilbert_metric_contracts} is a standard fact in real
projective geometry and can be verified by a computation. Note that we
allow the degenerate case $a = d$: in this situation the right-hand
side is identically zero for distinct $x, y \in (b,c)$. We allow no
other equalities among $a,b,c,d$, so the cross-ratio $[a,b; c,d]$ lies
in $\R - \{1\}$.

To apply \Cref{prop:hilbert_metric_contracts} to our situation, we
need to get some control on the behavior of the embeddings
$\iota:\pflags \to \P(V)$ and $\iota^*:\oppflags \to \P(V^*)$. We do
so in the next three lemmas below.

\begin{lem}
  \label{lem:subgroup_acts_secant_line}
  Let $x,y$ be distinct points in $\pflags$. There exists a
  one-parameter subgroup $g_t \subset G$ such that $\zeta(g_t)$ fixes
  $\iota(x)$ and $\iota(y)$, and acts nontrivially on the projective
  line $L_{xy}$ spanned by $\iota(x)$ and $\iota(y)$.
\end{lem}
\begin{proof}
  We can write $x = gP^+$ for some $g \in G$. Let $\mf{a}$ denote an
  abelian subalgebra of the Lie algebra $\mf{g}$ of $G$, such that for
  a maximal compact $K \subset G$, the exponential map $\mf{a} \to G$
  induces an isometric embedding $\mf{a} \to G/K$ whose image is a
  maximal flat in $G/K$.

  There is a conjugate $\mf{a}'$ of $\mf{a}$ such that the action of
  $\exp(\mf{a}')$ on $\pflags$ fixes both $x$ and $y$ (see
  \cite{eberlein}, Proposition 2.21.14). So, up to the action of $G$
  on $\pflags$, we can assume that $x$ is fixed by a standard
  parabolic subgroup $\stdpargroup$ conjugate to $P^+$, and that
  $x, y$ are both fixed by the subgroup $\exp(\mf{a})$.

  We choose $Z \in \mf{a}^+$ so that $\alpha(Z) \ne 0$ for all
  $\alpha \in \theta$. Then $g_t = \exp(tZ)$ is a 1-parameter subgroup
  of $G$ fixing $x$. As $t \to +\infty$, $g_t$ is
  $\stdpargroup$-divergent, with unique attracting fixed point $x$.

  Then \cite{ggkw2017anosov}, Lemma 3.7 implies that $\zeta(g_t)$ is
  $P_1$-divergent, where $P_1$ is the stabilizer of a line in $V$, and
  $\iota(x)$ is the unique one-dimensional eigenspace of $\zeta(g_t)$
  whose eigenvalue has largest modulus. And, since $\zeta(g_t)$ fixes
  $\iota(x)$ and $\iota(y)$, $\zeta(g_t)$ preserves $L_{xy}$, and acts
  nontrivially since the eigenvalues of $\zeta(g_t)$ on $\iota(x)$ and
  $\iota(y)$ must be distinct.
\end{proof}

\begin{lem}
  \label{lem:subgroup_acts_tangent_line}
  Let $L$ be any projective line in $\P(V)$ tangent to the image of
  the embedding $\iota:\pflags \to \P(V)$ at a point $\iota(x)$ for
  $x \in \pflags$. There exists a one-parameter subgroup $g_t$ of $G$
  so that $\zeta(g_t)$ acts nontrivially on $L$ with unique fixed
  point $\iota(x)$.
\end{lem}
\begin{proof}
  Fix a sequence $y_n \in \pflags$ such that $y_n \ne x$ and the
  projective line $L_n$ spanned by $\iota(x)$ and $\iota(y_n)$
  converges to $L$. By \Cref{lem:subgroup_acts_secant_line}, there
  exists $Z_n \in \mf{g}$ so that $\zeta(\exp(tZ_n))$ acts
  nontrivially on $L_n$, with fixed points $\iota(x)$ and
  $\iota(y_n)$.

  In the projectivization $\P(\mf{g})$, $[Z_n]$ converges to some
  $[Z]$. Since $\zeta:G \to \PGL(V)$ has finite kernel, there is an
  induced map $\zeta:\P(\mf{g}) \to \P(\mf{sl}(V))$, which satisfies
  \[
    \zeta([Z_n]) \to \zeta([Z]).
  \]
  A continuity argument shows that the one-parameter subgroup
  $\zeta(\exp(tZ))$ acts nontrivially on the line $L$, and has unique
  fixed point at $\iota(x)$.
\end{proof}

\begin{lem}
  \label{lem:line_intersections_open}
  Let $\Omega \subset \pflags$ be a proper domain, and let $L$ be a
  projective line in $\P(V)$ which is either spanned by two points in
  $\iota(\Omega)$, or is tangent to $\iota(\pflags)$ at a point
  $\iota(x)$ for $x \in \Omega$. Then
  \[
    W_L = \{v \in L : v = \iota^*(\xi) \cap L \textrm{ for } \xi \in
    \Omega^*\}
  \]
  is a nonempty open subset of $L$.
\end{lem}
\begin{proof}
  First observe that if $\xi \in \Omega^*$, then the intersection
  $\iota^*(\xi) \cap L$ is a singleton: since $\iota^*(\xi)$ is a
  hyperplane which does not contain any point in $\iota(\Omega)$, it
  must be in general position with $L$. So, since $\Omega^*$ is
  nonempty, $W_L$ contains at least one point $v$. Choose
  $\xi \in \Omega^*$ so that $\iota^*(\xi) \cap L = v$. We need to
  show that an open interval $I \subset L$ containing $v$ is contained
  in $W_L$.

  If $L$ is spanned by $x, y \in \iota(\Omega)$, then
  \Cref{lem:subgroup_acts_secant_line} implies that we can find a
  one-parameter subgroup $g_t \in G$ such that $\zeta(g_t)$ fixes $x$
  and $y$, and acts nontrivially on $L$. Since $\Omega^*$ is open, we
  can find $\eps > 0$ so that $g_t \cdot \xi \in \Omega^*$ for
  $t \in (-\eps, \eps)$. Since $x$ and $y$ are in $\iota(\Omega)$,
  $\iota^*(\xi)$ is transverse to both $x$ and $y$, so we have
  $v \ne x$, $v \ne y$. Then as $t$ varies from $-\eps$ to $\eps$,
  \[
    \iota^*(g_t \cdot \xi) \cap L = \zeta(g_t) \cdot v
  \]
  gives an open interval in $W_L$ containing $v$.

  A similar argument using \Cref{lem:subgroup_acts_tangent_line} shows
  that the claim also holds if $L$ is tangent to $\iota(\pflags)$ at a
  point in $\iota(\Omega)$.
\end{proof}

We can now prove a slightly weaker version of
\Cref{prop:caratheodory_metric_shrinks}, which we will then use to
show the stronger version.
\begin{lem}
  \label{lem:cpct_caratheodory_bound}
  Let $\Omega_1, \Omega_2$ be proper domains in $\pflags$, with
  $\overline{\Omega_1} \subset \Omega_2$, and let $K \subset \Omega_1$
  be compact. There exists a constant $\lambda > 1$ such that for all
  $x, y \in K$,
  \[
    C_{\Omega_1}(x,y) \ge \lambda \cdot C_{\Omega_2}(x,y).
  \]
\end{lem}
\begin{proof}
  Since $K$ is compact, it suffices to show that for fixed $x \in
  \Omega_1$, the ratio
  \[
    \frac{C_{\Omega_1}(x,y)}{C_{\Omega_2}(x,y)}
  \]
  is bounded below by some $\lambda > 1$ as $y$ varies in $K - \{x\}$.

  Let $y \in K - \{x\}$, and let $L_{xy}$ denote the projective line
  spanned by $\iota(x)$ and $\iota(y)$. Choose
  $\xi, \eta \in \overline{\Omega_2^*}$ so that
  \[
    C_{\Omega_2}(x,y) = \log|[\iota^*(\xi), \iota^*(\eta); \iota(x),
    \iota(y)]|.
  \]
  That is, if $v = \iota^*(\xi) \cap L_{xy}$,
  $w = \iota^*(\eta) \cap L_{xy}$, we have
  \[
    C_{\Omega_2}(x,y) = \log|[v, w; \iota(x), \iota(y)]| =
    \log\frac{|v - \iota(y)|\cdot|w - \iota(x)|}{|v - \iota(x)|\cdot|w
      - \iota(y)|},
  \]
  where the distances are measured in any identification of $L_{xy}$
  with $\RP^1 = \R \cup \{\infty\}$.

  We can choose an identification of $L_{xy}$ with
  $\R \cup \{\infty\}$ so that either $v < \iota(x) < \iota(y) < w$ or
  $v < \iota(x) < w < \iota(y)$. In either case, for any
  $v' \in (v, \iota(x)) \subset L_{xy}$, we have
  \[
    \log|[v', w; \iota(x), \iota(y)]| > \log|[v, w; \iota(x),
    \iota(y)]|.
  \]

  We know that $\overline{\Omega_2^*} \subset \Omega_1^*$, so
  $\xi, \eta$ lie in $\Omega_1^*$. Then
  \Cref{lem:line_intersections_open} implies that there exists
  $\xi' \in \Omega_1^*$ so that $v' = \iota^*(\xi') \cap L_{xy}$ lies
  in the interval $(v, \iota(x)) \subset L_{xy}$. See
  \Cref{fig:nested_flag_nbhds}.

  \begin{figure}[h]
    \centering
    \def\svgwidth{14cm}
    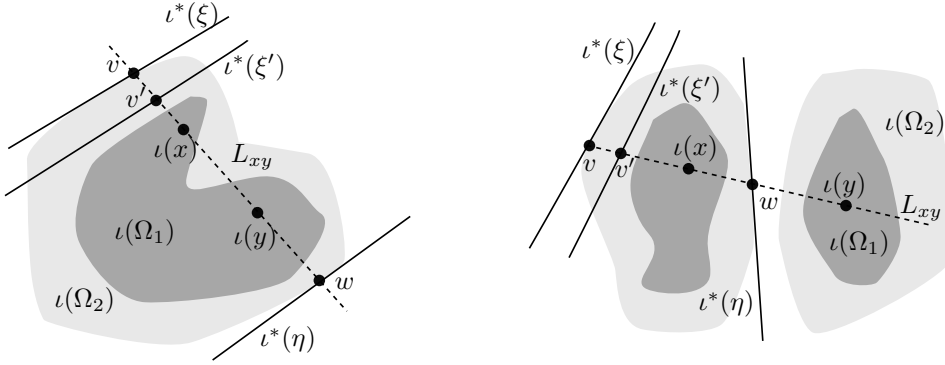
    \caption{We can always find $\xi' \in \Omega_1^*$ close to $\xi$
      so that the absolute value of the cross-ratio
      $[\iota^*(\xi), \iota^*(\eta); \iota(x), \iota(y)]$ increases
      when we replace $\xi$ with $\xi'$. In particular this is
      possible even when the sets $\iota(\Omega_1)$, $\iota(\Omega_2)$
      fail to be convex (left) or even connected
      (right).}\label{fig:nested_flag_nbhds}
  \end{figure}
  
  Then, we have
  \begin{align*}
    C_{\Omega_1}(x,y) &\ge \log|[\iota^*(\xi'), \iota^*(\eta) ;
                        \iota(x), \iota(y)]|\\
                      &= \log|[v', w; \iota(x), \iota(y)]\\
                      &> \log|[v, w; \iota(x), \iota(y)]\\
                      &= C_{\Omega_2}(x,y).
  \end{align*}
  This shows that $\frac{C_{\Omega_1}(x,y)}{C_{\Omega_2}(x,y)} > 1$
  for all $y \in K - \{x\}$. We still need to find some uniform
  $\lambda > 1$ so that
  $\frac{C_{\Omega_1}(x,y)}{C_{\Omega_2}(x,y)} \ge \lambda$ for all
  $y \in K - \{x\}$. To see this, suppose for the sake of a
  contradiction that for a sequence $y_n \in K - \{x\}$, we have
  \begin{equation}
    \label{eq:caratheodory_ratios_bounded}
    \frac{C_{\Omega_1}(x,y_n)}{C_{\Omega_2}(x, y_n)} \to 1.
  \end{equation}
  Since $K$ is compact, $y_n$ must converge to $x$. Up to subsequence,
  the sequence of projective lines $L_n$ spanned by $\iota(x)$ and
  $\iota(y_n)$ converges to a line $L$ tangent to $\iota(\pflags)$ at
  $\iota(x)$.

  For each $y_n$, choose $\xi_n$, $\eta_n \in \overline{\Omega_2^*}$
  so that
  \[
    C_{\Omega_2}(x, y_n) = \log|[\iota^*(\xi_n), \iota^*(\eta_n) ;
    \iota(x), \iota(y_n)]|.
  \]
  Let $v_n = \iota^*(\xi_n) \cap L_n$,
  $w_n = \iota^*(\eta_n) \cap L_n$. Then up to subsequence $\xi_n$
  converges to $\xi \in \overline{\Omega_2^*}$, $\eta_n$ converges to
  $\eta \in \overline{\Omega_2^*}$, and $v_n$, $w_n$ respectively
  converge to $v \in \iota^*(\xi) \cap L$ and
  $w \in \iota^*(\eta) \cap L$. Further, since $x$ is in $\Omega_2$,
  $\iota^*(\xi)$ and $\iota^*(\eta)$ are both transverse to
  $\iota(x)$. This implies that $\iota(x) \ne w$ and $\iota(x) \ne v$
  (although a priori we could have $v = w$). It also implies that
  neither $\iota^*(\xi)$ nor $\iota^*(\eta)$ contains the projective
  line $L$, and so
  \[
    v = \iota^*(\xi) \cap L, \qquad w = \iota^*(\eta) \cap L.
  \]
  Since $\xi \in \overline{\Omega_2^*} \subset \Omega_1^*$,
  \Cref{lem:line_intersections_open} implies that there exist
  $\xi', \eta' \in \Omega_1^*$ so that for some identification of $L$
  with $\R \cup \{\infty\}$, we have
  \[
    v < \iota^*(\xi') \cap L < \iota(x) < \iota^*(\eta') \cap L < w.
  \]
  Note that this is possible even if $v = w$, because then we can just
  identify both $v$ and $w$ with $\infty$. Let
  $v_n' = \iota^*(\xi') \cap L_n$, and let
  $w_n' = \iota^*(\eta') \cap L_n$. Respectively, $v_n'$ and $w_n'$
  converge to $v' = \iota^*(\xi') \cap L$ and
  $w' = \iota^*(\eta') \cap L$.

  This means that the cross-ratios $[v_n, v_n'; w_n', w_n]$ converge
  to $[v, v'; w', w] \in \R - \{1\}$, and in particular are bounded
  away from both $\infty$ and $1$ for all $n$.

  We choose identifications of $L_n$ with $\R \cup \{\infty\}$ so that
  $v_n < v_n' < \iota(x) < w_n' < w_n$. Since $y_n$ converges to $x$,
  for all sufficiently large $n$, we have $v_n' < \iota(y_n) <
  w_n'$. Then, \Cref{prop:hilbert_metric_contracts} implies that for
  all $n$, we have
  \[
    \log |[v_n', w_n'; \iota(x), \iota(y_n)]| \ge \lambda \cdot
    \log|[v_n, w_n; \iota(x), \iota(y_n)]|
  \]
  for some $\lambda > 1$ independent of $n$. But then since
  \[
    C_{\Omega_1}(x, y_n) \ge \log|[\iota^*(\xi'), \iota^*(\eta');
    \iota(x), \iota(y_n)]|,
  \]
  we have $C_{\Omega_1}(x,y_n) / C_{\Omega_2}(x,y_n) \ge \lambda$ for
  all $n$, contradicting (\ref{eq:caratheodory_ratios_bounded}) above.
\end{proof}

\begin{proof}[Proof of \Cref{prop:caratheodory_metric_shrinks}]
  We fix an open set $\Omega_{1.5}$ such that
  $\overline{\Omega_1} \subset \Omega_{1.5}$ and
  $\overline{\Omega_{1.5}} \subset \Omega_2$. Since
  $C_{\Omega_1}(x,y) \ge C_{\Omega_{1.5}}(x,y)$ for all
  $x, y \in \Omega_1$, we just need to see that there is some
  $\lambda > 1$ so that
  \[
    \frac{C_{\Omega_{1.5}}(x,y)}{C_{\Omega_2}(x,y)} \ge \lambda
  \]
  for all $x, y \in \overline{\Omega_1}$. This follows from
  \Cref{lem:cpct_caratheodory_bound}.
\end{proof}

\subsection{Contracting paths are $\pargroup$-divergent}

\begin{prop}
  \label{lem:graph_contracting_implies_contracting}
  Let $\mc{G}$ be a $\Gamma$-graph for a group $\Gamma \subset G$, and
  let $\{U_v\}$ be a $\mc{G}$-compatible system of open sets of
  $\pflags$ with each $U_v$ a proper domain.

  If $\alpha_n$ is a contracting $\mc{G}$-path, then the sequence
  \[
    \gamma_n = \alpha_1 \cdots \alpha_n
  \]
  is $\pargroup$-divergent with unique limit point $\xi$, where
  $\{\xi\} = \bigcap_{n=1}^\infty \gamma_n U_{v_{n+1}}$.
\end{prop}
\begin{proof}
  Consider the sequence of open sets
  \[
    \gamma_n \cdot U_{v_{n+1}}.
  \]
  Up to subsequence, $U_{v_{n+1}}$ is a fixed open set
  $U \subset \pflags$. By assumption $\{\alpha_n\}$ is a contracting
  path, so $\gamma_n \cdot U_{v_{n+1}}$ converges to a singleton
  $\{\xi\}$. So, we apply
  \Cref{prop:local_contracting_implies_global_contracting}.
\end{proof}


\section{A weaker criterion for EGF representations}
\label{sec:weak_egf}

In this section we assume that $P$ is a symmetric parabolic subgroup
of a semisimple Lie group $G$. Our goal is use the relative
quasigeodesic automaton and $\mc{G}$-compatible system of open sets
considered previously to prove a ``conical/peripheral''
characterization of EGF representations (we mentioned this
characterization earlier in \Cref{sec:egf_definition}).

The key idea of our ``conical/peripheral'' characterization is
contained in the following lemma. The point is that the
$\mc{G}$-compatible system gives us a way to exploit the connection
between $P$-divergence and nested subsets of $\flags$, discussed in
the previous section.
\begin{lem}
  \label{lem:automaton_divergent_limits}
  Let $(\Gamma, \mc{H})$ be a relatively hyperbolic pair and let
  $\rho:\Gamma \to G$ be a representation. Suppose there exists
  \begin{enumerate}
  \item a $\Gamma$-invariant closed set $\Lambda \subset \flags$ and a
    continuous equivariant surjective antipodal map
    $\phi:\Lambda \to \bgamh$, and
  \item a relative quasigeodesic automaton $\mc{G}$ and a
    $\mc{G}$-compatible system $\{U_v\}$ of open subsets of $\flags$,
    such that the $U_v$'s cover $\Lambda$ and each $U_v$ is a proper
    domain intersecting $\Lambda$ nontrivially.
  \end{enumerate}
  Then, for every sequence $\gamma_n \in \Gamma$ which is unbounded in
  the coned-off Cayley graph $\conedgraph$, the sequence
  $\rho(\gamma_n)$ is $P$-divergent, and every $P$-limit point of
  $\rho(\gamma_n)$ lies in $\Lambda$.
\end{lem}
\begin{proof}
  We will show that every subsequence of $\gamma_n$ has a
  $P$-contracting subsequence, so take an arbitrary subsequence of
  $\gamma_n$. By \Cref{lem:relative_automaton_quasidense_image}, we
  may assume that for a bounded sequence $b_n \in \Gamma$,
  $\gamma_nb_n$ is the endpoint of a finite $\mc{G}$-path
  $\{\alpha_m^{(n)}\}_{m=1}^{M_n}$. Up to subsequence $b_n$ is a
  constant $b$, independent of $n$.

  Let $\{v_m^n\}$ be the vertex path associated to
  $\{\alpha_m^{(n)}\}$. Up to subsequence $v_{M_n + 1}^n$ is a fixed
  vertex $v$, and $v_1^n$ is a fixed vertex $v'$. Let $U_{v'}^\eps$ be
  an $\eps$-neighborhood of $U_{v'}$, with $\eps$ chosen sufficiently
  small so that $U_{v'}^\eps$ is still a proper domain.

  The sequence $M_n$ must be unbounded, since the length of $\gamma_n$
  with respect to the coned-off Cayley graph metric is at most a fixed
  constant times $M_n$. \Cref{cor:exponential_path_shrinking} then
  implies that the diameter of
  \[
    \rho(\gamma_nb) \cdot U_v = \rho(\alpha_1^{(n)}) \cdots
    \rho(\alpha_{M_n}^{(n)}) U_v
  \]
  with respect to the metric $C_{U_{v'}^\eps}$ tends to zero,
  exponentially in $M_n$. Since this sequence of sets lies in the
  compact set $\overline{U_{v'}} \subset U_{v'}^\eps$, up to
  subsequence it must converge to a singleton $\{\xi\}$ in $G/P$. In
  fact $\xi$ must lie in $\Lambda$, because $\Lambda$ is compact and
  $\xi$ is the limit of a sequence of points in the sequence of
  nonempty closed sets
  $(\rho(\gamma_nb) \cdot \overline{U_v}) \cap \Lambda$. Then, since
  $\rho(\gamma_n) \cdot \rho(b)U_v$ converges to $\{\xi\}$,
  \Cref{prop:local_contracting_implies_global_contracting} implies
  that $\rho(\gamma_n)$ is $P$-divergent with unique $P$-limit $\xi$.
\end{proof}

An immediate corollary of the above (together with the construction in
section 6 of this paper) is the following:
\begin{cor}
  \label{cor:egf_rel_divergent}
  Let $(\Gamma, \mc{H})$ be a relatively hyperbolic pair, and let
  $\rho:\Gamma \to G$ be a $P$-EGF representation. Then $\rho$ is
  relatively $P$-divergent.
\end{cor}

\begin{lem}
  \label{lem:conical_convergence_compatible}
  Let $(\Gamma, \mc{H})$ be a relatively hyperbolic pair, let
  $\rho:\Gamma \to G$ be a representation, let
  $\Lambda \subset \flags$ be a closed $\Gamma$-invariant set, and let
  $\phi:\Lambda \to \bgamh$ be a continuous equivariant surjective
  antipodal map.

  Suppose that $\gamma_n \in \Gamma$ is a sequence converging to
  $z_+ \in \conbdry(\Gamma, \mc{H})$, such that $\rho(\gamma_n)$ is
  $P$-divergent and every $P$-limit point of $\rho(\gamma_n^{\pm 1})$
  lies in $\Lambda$. If $\gamma_n^{-1}$ converges to $z_- \in \bgamh$,
  then for every compact set $K \subset \Opp(\phi^{-1}(z_-))$ and
  every open $U$ containing $\phi^{-1}(z_+)$, for large enough $n$, we
  have $\rho(\gamma_n)K \subset U$.
\end{lem}
\begin{proof}
  Assume $\Gamma$ is non-elementary. It suffices to show that every
  subsequence of $\gamma_n$ has a further subsequence satisfying the
  desired property. So, we can freely extract subsequences throughout
  this proof.

  We assume $P$ symmetric, so $\rho(\gamma_n^{-1})$ is also
  $P$-divergent and has nonempty $P$-limit set. Let $\xi_\pm$ be a
  pair of flags in the $P$-limit sets of $\rho(\gamma_n^{\pm 1})$,
  respectively; by assumption we have $\xi_\pm \in \Lambda$. By
  \Cref{lem:klp_contraction_symmetry}, we have a subsequence so that
  $\rho(\gamma_n)$ converges to $\xi_+$ uniformly on compacts in
  $\Opp(\xi_-)$.

  Antipodality of $\phi$ implies that every compact subset of
  $\bgamh - \{\phi(\xi_-)\}$ is contained in
  $\phi(\Opp(\xi_-) \cap \Lambda)$. Then, by equivariance and
  continuity of $\phi$, we see that $\gamma_n$ must converge to
  $\phi(\xi_+)$ on compacts in $\bgamh - \{\phi(\xi_-)\}$. If we
  assume $\Gamma$ is non-elementary, this uniquely characterizes the
  points $\phi(\xi_\pm)$ as the limits of $\gamma_n^{\pm 1}$ in
  $\bgamh$. So, we see that $\rho(\gamma_n)$ converges uniformly to
  $\xi_+ \in \phi^{-1}(z_+)$ on every compact in
  $\Opp(\phi^{-1}(z_-)) \subset \Opp(\xi_-)$, as required.
\end{proof}

We recall the statement of our weaker characterization of EGF
representations here:

\conicalPeripheralEGF*
\begin{proof}
  To see the ``only if'' part, observe that if we know that $\phi$ is
  an EGF boundary extension, we can use the results of
  \Cref{sec:extended_convergence} to construct an automaton satisfying
  the hypotheses of \Cref{lem:automaton_divergent_limits}, which
  immediately implies that condition \ref{item:conical_extended_con}
  holds. Condition \ref{item:peripheral_extended_con} is immediate
  from the assumption that $\phi$ extends the convergence action of
  $\Gamma$ on $\bgamh$.

  So, we focus on the ``if'' part. Without loss of generality we can
  assume that for every conical $z \in \bgamh$, the set $C_z$ is equal
  to the set $\Opp(\phi^{-1}(z))$: for this choice of strata,
  antipodality of $\phi$ and the fact that antipodality is an open
  condition implies that \ref{item:strata_interior} still holds.
  
  Now, observe that if $\gamma_n$ is a sequence limiting conically to
  $z_+ \in \bgamh$, with $\gamma_n^{-1}$ converging to $z_-$, then
  \Cref{lem:conical_convergence_compatible}, together with condition
  \ref{item:conical_extended_con}, implies that the map
  $\phi:\Lambda \to \flags$ and strata $C_z$ satisfy both conditions
  \ref{item:conical_extended_convergence} and
  \ref{item:peripheral_extended_convergence} given at the beginning of
  \Cref{sec:extended_convergence}. So, by
  \Cref{prop:convergence_dynamics_implies_automaton}, we know that
  there is a relative quasigeodesic automaton $\mc{G}$ satisfying the
  hypotheses of \Cref{lem:automaton_divergent_limits}. Thus, the
  representation $\rho$ is relatively $P$-divergent, and for every
  sequence $\gamma_n$ in $\Gamma$ which is unbounded in $\conedgraph$,
  all of the $P$-limit points of $\rho(\gamma_n)$ are contained in
  $\Lambda$. But this is precisely condition \ref{item:rel_p_div} in
  \Cref{prop:weak_egf_to_egf_ii}, meaning the representation $\rho$
  and map $\phi:\Lambda \to \bgamh$ satisfy all the conditions in
  \Cref{prop:weak_egf_to_egf_ii} and $\rho$ is EGF.
\end{proof}

The arguments above also imply the following characterization of EGF
representations. This result is not needed anywhere else in the paper.
\begin{prop}
  \label{prop:weakest_egf}
  Let $(\Gamma, \mc{H})$ be a relatively hyperbolic pair, let
  $\rho:\Gamma \to G$ be a representation, and let $P \subset G$ be a
  symmetric parabolic subgroup.

  Then $\rho$ is EGF if and only if there is a closed $\rho$-invariant
  set $\Lambda \subseteq \flags$, a surjective equivariant antipodal
  map $\phi:\Lambda \to \bgamh$, and open subsets $C_z \subset \flags$
  for each $z \in \bgamh$ satisfying each of the following:
  \begin{enumerate}
    \custitem[\ref{item:strata_interior}] For every compact
    $K \subset \bgamh$, the set $\bigcap_{z \in K}C_z$ contains an
    open neighborhood of $\phi^{-1}(\bgamh - K)$.

    \custitem[(C2-con)] For any sequence $\gamma_n \in \Gamma$
    limiting conically to some point $z$ in $\bgamh$, with
    $\gamma_n^{-1} \to z_-$, any open set $U$ containing
    $\phi^{-1}(z)$, and any compact $K \subset C_{z_-}$, we have
    $\rho(\gamma_n) \cdot K \subset U$ for all sufficiently large $n$.
    
    \custitem[(C2-par)] For any parabolic point $p \in \bgamh$, any
    compact $K \subset C_p$, and any open set $U$ containing
    $\phi^{-1}(p)$, for all but finitely many $\gamma \in \Gamma_p$,
    we have $\rho(\gamma) \cdot K \subset U$.
  \end{enumerate}
\end{prop}
\begin{proof}
  The ``only if'' direction is immediate from the definition, so
  suppose we have a representation satisfying the hypotheses
  above. The results of \Cref{sec:extended_convergence} imply that
  there is a relative quasigeodesic automaton satisfying the
  hypotheses of \Cref{lem:automaton_divergent_limits}, so $\rho$ is
  relatively $P$-divergent, and all of the relative $P$-limit points
  lie in $\Lambda$. Thus \Cref{prop:weak_egf_to_egf_ii} applies and
  $\rho$ is EGF.
\end{proof}


\section{Relative stability}
\label{sec:relative_stability}
In this section we prove the main \emph{relative stability property}
for EGF representations (\Cref{thm:cusp_stable_stability}).

\subsection{Deformations of EGF representations}

In general, the set of EGF representations is not an open subset of
$\Hom(\Gamma, G)$. However, it is relatively open in a subspace of
$\Hom(\Gamma, G)$ where we restrict the deformations of the peripheral
subgroups appropriately. Roughly speaking, we want to consider
subspaces $\mc{W} \subset \Hom(\Gamma, G)$ where the \emph{dynamical}
behavior of the peripheral subgroups changes continuously under
deformation. That is, if $p$ is a parabolic point, and $\rho_0$ is a
representation of $\Gamma$ where $\rho_0(\Stab_\Gamma(p))$ attracts
points towards a particular subset $\Lambda_p \subset \flags$ at a
particular ``speed,'' then we want to work with small deformations
$\rho_t$ of $\rho_0$ where $\rho_t(\Stab_\Gamma(p))$ attracts points
towards a small deformation of $\Lambda_p$ at a similar ``speed.''

The precise condition is the following:

\begin{definition}
  \label{defn:peripheral_stability}
  Let $\rho_0:\Gamma \to G$ be an EGF representation with boundary
  extension $\phi:\Lambda \to \bgamh$, and let
  $\mc{W} \subset \Hom(\Gamma, G)$ contain $\rho_0$.

  We say that $\mc{W}$ is \emph{peripherally stable at
    $(\rho_0, \phi)$} if for every $p \in \parbdry(\Gamma, \mc{H})$,
  every open set $U$ containing $\phi^{-1}(p)$, every compact set
  $K \subset C_p$, and every cofinite set $T \subset \Stab_\Gamma(p)$
  such that
  \[
    \rho_0(T) \cdot K \subset U,
  \]
  there is an open set $\mc{W}' \subset \mc{W}$ containing $\rho_0$,
  such that for every $\rho' \in \mc{W}'$, we have
  \[
    \rho'(T) \cdot K \subset U.
  \]
\end{definition}

We restate the main result of the paper below:

\cuspStableStability*

\begin{remark}
  In \cite{bowditch1998spaces}, Bowditch explored the deformation
  spaces of geometrically finite groups $\Gamma \subset \PO(d,1)$, and
  gave an explicit discription of semialgebraic subspaces of
  $\Hom(\Gamma, \PO(d,1))$ in which small deformations of $\Gamma$ are
  still geometrically finite.

  Bowditch's deformation spaces are peripherally stable, so it seems
  desirable to find a general algebraic description of peripherally
  stable subspaces.

  Even in $\PO(d,1)$, the question is subtle, however. Bowditch also
  gives examples of geometrically finite representations
  $\rho:\Gamma \to \PO(d,1)$ (for $d \ge 4$) and deformations $\rho_t$
  of $\rho$ in $\Hom(\Gamma, \PO(d,1))$ such that the restriction of
  $\rho_t$ to each cusp group in $\Gamma$ is discrete, faithful, and
  parabolic, but $\rho_t$ is not even discrete; further examples exist
  where the deformation is discrete, but not geometrically finite.
\end{remark}

\begin{example}
  \label{ex:peripheral_stability}
  Let $B \in \SL(d, \R)$ be a $d$-dimensional Jordan block with
  eigenvalue 1 and eigenvector $v$, and let $A \in \SL(d + 2, \R)$ be
  the block matrix
  $\begin{psmallmatrix}B \\& 1\\ & & 1\end{psmallmatrix}$.

  Although $[v]$ is not quite an attracting fixed point of $A$, it is
  still an ``attracting subspace'' in the sense that if $K$ is any
  compact subset of $\RP^{d+1}$ which does not intersect a fixed
  hyperplane of $A$, then $A^n \cdot K$ converges to $\{[v]\}$. Via a
  ping-pong argument, one can use this ``attracting'' behavior to show
  that for some $k \ge 1$ and some $M \in \SL(d + 2, \R)$, the group
  $\Gamma$ generated by $\alpha = A^k$ and $\beta = MA^kM^{-1}$ is a
  discrete free group with free generators $\alpha, \beta$. The group
  $\Gamma$ is hyperbolic relative to the subgroups
  $\langle\alpha\rangle$, $\langle \beta \rangle$, and the inclusion
  $\Gamma \hookrightarrow \SL(d + 2, \R)$ is EGF with respect to
  $P_{1,d+1}$ (the stabilizer of a line in a hyperplane in
  $\R^{d+2}$).

  Here, there are peripherally stable deformations of $\Gamma$ which
  change the Jordan block decomposition of $A$. For instance, consider
  a continuous path $A_t:[0, 1] \to \SL(d+2, \R)$ given by
  $A_t = \begin{psmallmatrix}B_t\\& 1 \\ & & 1\end{psmallmatrix}$,
  where $B_0 = B$ and $B_t$ is a \emph{diagonalizable} matrix in
  $\SL(d, \R)$. For small values of $t$, the group $\Gamma_t$
  generated by $\alpha_t = A_t^k$ and $\beta$ is still discrete and
  freely generated by $\alpha_t$ and $\beta$---since the
  ``attracting'' fixed points of $A_t$ deform continuously with $t$,
  the same exact ping-pong setup works for all small $t \ge 0$. And
  indeed the path in $\Hom(\Gamma, \SL(d+2, \R))$ determined by the
  path $A_t$ is a peripherally stable subspace.

  On the other hand, consider the path
  $A_t' = \begin{psmallmatrix}B\\& e^t\\ & &
    e^{-t}\end{psmallmatrix}$, and let $\alpha_t' = A_t'^k$. In this
  case the corresponding subspace of $\Hom(\Gamma, \SL(d+2, \R))$ is
  \emph{not} peripherally stable: while the group generated by
  $\alpha_t'$ is still discrete, the attracting fixed points of $A_t'$
  do \emph{not} deform continuously in $t$. So, there is no way to use
  the ping-pong setup for $\Gamma$ to ensure that
  $\Gamma_t' = \langle \alpha_t', \beta\rangle$ is a discrete group.
\end{example}

\begin{example}
  Here is a somewhat more interesting example of a
  \emph{non}-peripherally stable deformation. Let $M$ be a
  finite-volume noncompact hyperbolic 3-manifold, with holonomy
  representation $\rho:\pi_1M \to \PSL(2, \C)$ (so there is an
  identification $M = \H^3 / \rho(\pi_1M)$). Then $\pi_1M$ is
  hyperbolic relative to the collection $\mc{C}$ of conjugates of cusp
  groups (each of which is isomorphic to $\Z^2$), and the
  representation $\rho$ is geometrically finite (in particular, EGF).

  In this case, for \emph{any} sufficiently small nontrivial
  deformation $\rho'$ of $\rho$ in the character variety
  $\Hom(\pi_1M, \PSL(2, \C)) / \PSL(2, \C)$, the restriction of
  $\rho'$ to some cusp group $C \in \mc{C}$ either fails to be
  discrete or has infinite kernel. So $\Hom(\pi_1M, \PSL(2, \C))$ is
  \emph{not} peripherally stable, because any sufficiently small
  deformation of $\rho$ inside of a peripherally stable subspace must
  have discrete image and finite kernel on each $C \in \mc{C}$. This
  is true despite the fact that arbitrarily small deformations of
  $\rho$ are holonomy representations of complete hyperbolic structures on Dehn
  fillings of $M$ (so in particular, they are discrete).
\end{example}

The main ingredient in the proof of \Cref{thm:cusp_stable_stability}
is the relative quasigeodesic automaton $\mc{G}$ and the associated
$\mc{G}$-compatible system of open sets $\{U_v\}$ we constructed in
\Cref{prop:convergence_dynamics_implies_automaton}. The following
proposition is immediate from the definition of peripheral stability:

\begin{prop}
  Let $\rho:\Gamma \to G$ be an EGF representation with boundary
  extension $\phi$, and let $\mc{W} \subset \Hom(\Gamma, G)$ be a subspace
  which is peripherally stable at $(\rho, \phi)$.

  If $\mc{G}$ is a relative quasigeodesic automaton for $\Gamma$, and
  $\{U_v\}$ is a $\mc{G}$-compatible system of open subsets of
  $\flags$ for $\rho(\Gamma)$, then there is an open subset
  $\mc{W}' \subset \mc{W}$ containing $\rho$ such that for every
  $\rho' \in \mc{W}'$, $\{U_v\}$ is also a $\mc{G}$-compatible system
  of open sets for $\rho'(\Gamma)$.
\end{prop}

\Cref{thm:cusp_stable_stability} then follows from a kind of converse
to \Cref{prop:convergence_dynamics_implies_automaton}: we will show
that we can reconstruct a map extending the convergence dynamics of
$\Gamma$ from the $\mc{G}$-compatible system $\{U_v\}$.

\subsection{An equivariant map on conical limit points}

For the rest of this section, we let $\rho:\Gamma \to G$ be a
representation which is EGF with respect to a symmetric parabolic
subgroup $P \subset G$. We let $\phi:\Lambda \to \bgamh$ be a boundary
extension for $\rho$, and assume that $\mc{W} \subset \Hom(\Gamma, G)$
is peripherally stable at $(\rho, \phi)$. We also let $Z$ be a compact
subset of $\bgamh$, and let $V \subset \flags$ be an open subset
containing $\phi^{-1}(Z)$. We again fix a finite subset
$\Pi \subset \parbdry(\Gamma, \mc{H})$, containing one point from
every $\Gamma$-orbit in $\parbdry(\Gamma, \mc{H})$.

Using \Cref{prop:convergence_dynamics_implies_automaton}, we can find
a relative quasigeodesic automaton $\mc{G}$ and $\mc{G}$-compatible
system $\{U_v\}$ of open subsets of $\flags$ for $\rho(\Gamma)$. Using
\Cref{prop:compact_bowditch_nbhd}, we can ensure that for any
$z \in Z$, there is a $\mc{G}$-path $\{\alpha_n\}$ limiting to $z$
(with vertex path $\{v_n\}$) so that $U_{v_1}$ is contained in $V$.

For each $p \in \Pi$, we also fix a compact set
$K_p \subset \bgamh - \{p\}$ such that
$\Gamma_p \cdot K_p = \bgamh - \{p\}$, and assume that the automaton
$\mc{G}$ has been constructed to satisfy
\Cref{prop:parabolic_outward_edge_nbhds}.

Antipodality of the map $\phi$ implies that for each $z \in \bgamh$,
each fiber $\phi^{-1}(z)$ is a closed subset of some affine chart in
$\flags$. So, we can also assume that $U_v$ is a proper domain for
each vertex $v$ of $\mc{G}$. In fact, by way of the following lemma,
we can assume even more:

\begin{lem}
  \label{lem:compatible_system_transverse}
  Let $\rho$ be an EGF representation with boundary map
  $\phi:\Lambda \to \bgamh$.

  For any $\delta > 0$, we can find a relative quasigeodesic automaton
  $\mc{G}$ with $\mc{G}$-compatible system $\{U_v\}$ of open sets in
  $\flags$ as in \Cref{prop:convergence_dynamics_implies_automaton},
  so that for any $x, y \in \bgamh$ with $d(x,y) > \delta$, if
  $\phi^{-1}(x) \subset U_v$ and $\phi^{-1}(y) \subset U_w$, then
  $\overline{U_v}$ and $\overline{U_w}$ are opposite.
\end{lem}
\begin{proof}
  We choose $\eps > 0$ so that if $d(v, w) > \delta/2$ for
  $v, w \in \bgamh$, then the closed $\eps$-neighborhoods of
  \[
    \phi^{-1}(v), \qquad \phi^{-1}(w)
  \]
  are opposite. This is possible for a fixed pair $v, w \in \bgamh$
  since antipodality is an open condition, and $\phi^{-1}(v)$,
  $\phi^{-1}(w)$ are opposite compact sets. Then we can pick a uniform
  $\eps$ for all pairs since the the subset
  $\{(u, v) \in (\bgamh)^2: d(u,v) \ge \delta/2\}$ is compact.

  Shrinking $\eps$ if necessary, we can also use uniform continuity of
  the map $\phi:\Lambda \to \bgamh$ to ensure that if $\xi, \xi'$ and
  $\nu, \nu'$ are flags in $\Lambda$ satisfying
  $d(\phi(\xi), \phi(\nu)) > \delta$ and
  $d(\xi, \xi'), d(\nu, \nu') < \eps$, then
  $d(\phi(\xi'), \phi(\nu')) > \delta/2$.

  Consider $\mc{G}$-compatible systems of open subsets $\{U_v\}$ and
  $\{W_v\}$ for the action of $\Gamma$ on $\flags$ and $\bgamh$,
  coming from \Cref{prop:convergence_dynamics_implies_automaton}. We
  can ensure that for each vertex $a$ we have
  $U_a \subset N(\phi^{-1}(w), \eps)$ for some $w \in W_a$.

  Now, let $x, y \in \bgamh$ satisfy $d(x,y) > \delta$, and assume
  that $\phi^{-1}(x) \subset U_a$ and $\phi^{-1}(y) \subset U_b$ for
  vertices $a, b$. From the above we know that
  \[
    U_a \subset N(\phi^{-1}(v), \eps), \qquad U_b \subset
    N(\phi^{-1}(w), \eps)
  \]
  for $v \in W_a$ and $w \in W_b$. Thus we have flags
  $\xi \in \phi^{-1}(x)$ and $\xi' \in \phi^{-1}(v)$ satisfying
  $d(\xi, \xi') < \eps$, and flags
  $\nu \in \phi^{-1}(y), \nu' \in \phi^{-1}(w)$ satisfying
  $d(\nu, \nu') < \eps$. From our choice of $\eps$, we have that
  $d(v, w) > \delta/2$. So our choice of $\eps$ also ensures that the
  closures of $N(\phi^{-1}(w), \eps)$ and $N(\phi^{-1}(v), \eps)$ are
  opposite, and therefore so are $\overline{U_a}$ and
  $\overline{U_b}$.
\end{proof}

Using cocompactness of the action of $\Gamma$ on the space of distinct
pairs in $\bgamh$, we know that there exists some fixed
$\delta > 0$ such that for any distinct
$z_1, z_2 \in \bgamh$, we can find some
$\gamma \in \Gamma$ such that $d(\gamma z_1, \gamma z_2) >
\delta$. Then, in light of \Cref{lem:compatible_system_transverse}, we
can make the following assumption:
\begin{assumption}
  \label{ass:opposite_open_sets}
  For any $z_1, z_2 \in \bgamh$ satisfying $d(z_1, z_2) > \delta$, if
  $\phi^{-1}(z_1) \subset U_v$ and $\phi^{-1}(z_2) \subset U_w$ for
  $v,w$ vertices of $\mc{G}$, then $\overline{U_v}$ and
  $\overline{U_w}$ are opposite.
\end{assumption}

With our relative quasigeodesic automaton $\mc{G}$ and compatible
system of open sets $\{U_v\}$ fixed, we now choose an open subset
$\mc{W}' \subset \mc{W}$ so that for any $\rho' \in \mc{W}'$,
$\{U_v\}$ is also a $\mc{G}$-compatible system for the action of
$\rho'(\Gamma)$ on $\flags$. Our main goal for the rest of this
section is to show that any $\rho' \in \mc{W}'$ is an EGF
representation. So, we fix some $\rho' \in \mc{W}'$.

Let $\mr{Path}(\mc{G})$ denote the set of infinite $\mc{G}$-paths.
\Cref{prop:compatible_system_in_flags_contracting} implies that every
path in $\mr{Path}(\mc{G})$ is contracting for the $\rho'$-action, so
we have a map
\[
  \psi_{\rho'}:\mr{Path}(\mc{G}) \to \flags,
\]
where the path $\{\alpha_n\}$ maps to the unique element of
\[
  \bigcap_{n=1}^\infty \rho'(\alpha_1) \cdots \rho'(\alpha_n)
  U_{v_{n+1}}.
\]
\begin{lem}
  \label{lem:equivariant_conical_map}
  The map $\psi_{\rho'}:\mr{Path}(\mc{G}) \to \flags$ induces an
  equivariant map
  \[
    \psi_{\rho'}:\conbdry \Gamma \to \flags.
  \]
\end{lem}
\begin{proof}
  We first need to see that $\psi_{\rho'}$ is well-defined, i.e. that
  if $z$ is a conical limit point and $\{\alpha_n\}$, $\{\beta_n\}$
  are $\mc{G}$-paths limiting to $z$, then
  $\psi_{\rho'}(\{\alpha_n\}) = \psi_{\rho'}(\{\beta_n\})$.

  Let
  \[
    \gamma_n = \alpha_1 \cdots \alpha_n, \qquad \eta_m = \beta_1
    \cdots \beta_m.
  \]
  By \Cref{lem:codings_hausdorff_close}, the Hausdorff distance
  between the sets $\{\gamma_n\}$ and $\{\eta_m\}$ is
  finite. \Cref{lem:graph_contracting_implies_contracting} implies
  that $\rho'(\gamma_n)$ and $\rho'(\eta_n)$ are both $P$-divergent
  sequences and each have a unique $P$-limit point in $\flags$, given
  by $\psi_{\rho'}(\{\alpha_n\})$, $\psi_{\rho'}(\{\beta_m\})$,
  respectively. Combining the bound on Hausdorff distance with
  \cite[Lemma 4.23]{klp2017anosov}, we see that the set of $P$-limit
  points of $\rho'(\gamma_n)$ and $\rho'(\eta_n)$ must agree and
  therefore $\psi_{\rho'}(\{\alpha_n\}) = \psi_{\rho'}(\{\beta_m\})$.

  The argument for equivariance is nearly the same. Fix
  $g \in \Gamma$. Then if $\{\alpha_n\}$ is a $\mc{G}$-path limiting
  to $z$, and $\{\beta_m\}$ is a $\mc{G}$-path limiting to $gz$,
  \Cref{lem:codings_hausdorff_close} also tells us that the sequences
  \[
    \{g\alpha_1 \cdots \alpha_n\}, \qquad \{\beta_1 \cdots \beta_m\}
  \]
  have finite Hausdorff distance in $\cay(\Gamma)$. Moreover, the same
  argument as above tells us that $\rho'$ maps each of these sequences
  to a $P$-divergent sequence, with unique $P$-limit points given
  respectively by
  \[
    \rho'(g)\psi_{\rho'}(\{\alpha_n\}), \qquad
    \psi_{\rho'}(\{\beta_m\}).
  \]
  Again, the bound on Hausdorff distance means that these $P$-limit
  points actually agree.
\end{proof}

Later we will see that $\psi_{\rho'}$ is both continuous and
injective.

\subsection{Extending $\psi_{\rho'}$ to parabolic points}

We want to extend the map $\psi_{\rho'} : \conbdry\Gamma \to \flags$
to the entire Bowditch boundary $\bgamh$. To do so, we need to
view $\psi_{\rho'}$ as a map to the set of closed subsets of $\flags$.

The first step is to define $\psi_{\rho'}$ on the finite set
$\Pi \subset \parbdry\Gamma$. The idea is to define $\psi_{\rho'}(p)$
for $p \in \Pi$ to be a set of ``accumulation points'' for orbits in
$\flags$ for the $\rho'$-action of $\Gamma_p$. Since $\rho'(\Gamma_p)$
need not be $P$-divergent, these accumulation points are \emph{not}
generically independent of the choice of orbit. However, the system of
open sets attached to our automaton (together with the peripheral
stability condition) will allow us to make a choice which is
reasonably well-behaved.

The specific construction is the following. For any vertex $v$ in
$\mc{G}$, we consider the set
\[
  B_v = \bigcup_{(v, y) \textrm{ edge in } \mc{G}} U_y.
\]
Then, for each $p \in \Pi$, we pick a parabolic vertex $w$ so that
$p_w = p$. We define $\Lambda_p'$ to be the closure of the set of
accumulation points of sequences of the form
$\rho'(\gamma_n) \cdot x$, for $x \in B_w$ and $\gamma_n$ distinct
elements of $\Gamma_{p}$. Part (3) of
\Cref{prop:convergence_dynamics_implies_automaton} guarantees that
$B_w \subset C_p$, and $\mc{G}$-compatibility of the system $\{U_v\}$
for the $\rho'(\Gamma)$-action on $\flags$ implies that
$\Lambda_{p}' \subset U_w$. By construction, $\Lambda_{p}'$ is
invariant under the action of $\rho(\Gamma_p)$.

Next, given a parabolic point $q \in \parbdry \Gamma$, we write
$q = g \cdot p$ for $p \in \Pi$, and then define
\[
  \psi_{\rho'}(q) := \rho'(g) \Lambda_p'.
\]
Since $\Lambda_p'$ is $\Gamma_p$-invariant and $\Gamma_p$ is exactly
the stabilizer of $p$, this does not depend on the choice of coset
representative in $g\Gamma_p$. Moreover $\psi_{\rho'}$ is still
$\rho'$-equivariant.

In addition, if $v$ is any parabolic vertex with parabolic point
$p_v = g \cdot p$ for $p \in \Pi$, part (2) of
\Cref{prop:convergence_dynamics_implies_automaton} ensures that
$B_v = B_w$ for any parabolic vertex $w$ with $p_w = p$. So,
$\rho'(g) \cdot \Lambda_p'$ is exactly the closure of the set of
accumulation points of the form $\rho'(g\gamma_n)\cdot x$ for
sequences $\gamma_n \in \Gamma_p$ and $x \in B_v$. Then
$\mc{G}$-compatibility implies that
$\psi_{\rho'}(p_v) = \rho'(g)\Lambda_p'$ is a subset of $U_v$.

\begin{remark}
  There is a natural topology on the space of closed subsets of
  $\flags$, induced by the Hausdorff distance arising from some (any)
  choice of metric on $\flags$. We emphasize that the map
  $\psi_{\rho'}$ is \emph{not} necessarily continuous with respect to
  this topology.
\end{remark}

Ultimately we want to use $\psi_{\rho'}$ to define a map extending the
convergence dynamics of $\Gamma$, so we will need to also define the
sets $C_z'$ for each $z \in \bgamh$. For now, we only define
$C_p'$ for $p \in \Pi$: this will be the set
\[
  \bigcup_{ \gamma \in \Gamma_p} \rho'(\gamma) B_w.
\]
We can immediately observe:
\begin{prop}
  \label{prop:convergence_dynamics_deformed_parabolics}
  $C_p'$ is $\rho'(\Gamma_p)$-invariant. Moreover, for any infinite
  sequence $\gamma_n \in \Gamma_p$, any compact $K \subset C_p'$, and
  any open $U \subset \flags$ containing $\Lambda_p'$, for
  sufficiently large $n$, $\rho'(\gamma_n) \cdot K$ lies in $U$.
\end{prop}
\begin{proof}
  First, $\Gamma_p$-invariance follows directly from the
  definition. Fix a compact $K \subset C_p'$ and an open
  $U \subset \flags$ containing $\Lambda_p'$. Then $K$ is contained in
  finitely many sets $\rho'(\gamma) B_w$ for $\gamma \in \Gamma_p$, so
  any accumulation point of $\rho'(\gamma_n) x$ for $x \in K$ and
  $\gamma_n \in \Gamma_p$ lies in $\Lambda_p'$. In particular, for
  sufficiently large $n$, $\rho'(\gamma_n)x$ lies in $U$, and since
  $K$ is compact we can pick $n$ large enough so that
  $\rho'(\gamma_n)x \in U$ for all $x \in K$.
\end{proof}

We next want to use $\psi_{\rho'}$ to define an antipodal extension
from a subset of $\flags$ to $\bgamh$.

\begin{lem}
  \label{lem:initial_vertex_neighborhood}
  For any $z \in \bgamh$, if $\{\alpha_n\}$ is a $\mc{G}$-path
  limiting to $z$ with corresponding vertex path $\{v_n\}$, then
  $\phi^{-1}(z)$ and $\psi_{\rho'}(z)$ are both subsets of $U_{v_1}$.
\end{lem}
\begin{proof}
  If $z$ is a conical limit point, then this follows immediately from
  \Cref{prop:bowditch_graphs_geodesic} and the definition of
  $\psi_{\rho'}$. On the other hand, if $z$ is a parabolic point, then
  $z = \alpha_1 \cdots \alpha_N p_v$, where $v$ is a parabolic vertex
  at the end of the vertex path $\{v_n\}$. By part (3) of
  \Cref{prop:convergence_dynamics_implies_automaton}, we have
  $p_v \in W_v$ and thus $\phi^{-1}(p_v) \subset U_v$. By
  $\rho$-equivariance of $\phi$ we have
  \[
    \phi^{-1}(z) = \rho(\alpha_1 \cdots \alpha_N) \phi^{-1}(p_v),
  \]
  so by $\mc{G}$-compatibility we have $\phi^{-1}(z) \subset
  U_{v_1}$. On the other hand, we have constructed $\psi_{\rho'}$ so
  that $\psi_{\rho'}(p_v) \subset U_v$, so $\rho'$-equivariance of
  $\psi_{\rho'}$ and $\mc{G}$-compatibility also show that
  $\psi_{\rho'}(z) \subset U_{v_1}$.
\end{proof}

\begin{lem}
  \label{lem:boundary_map_well_defined}
  For any two distinct points $z_1, z_2$ in $\bgamh$, the sets
  \[
    \psi_{\rho'}(z_1), \qquad \psi_{\rho'}(z_2)
  \]
  are opposite (in particular disjoint).
\end{lem}
\begin{proof}
  We know that for any distinct $z_1, z_2 > 0$, we can find
  $\gamma \in \Gamma$ so that $d(\gamma z_1, \gamma z_2) >
  \delta$. So, since $\psi_{\rho'}$ is $\rho'$-equivariant, we just
  need to show that if $z_1, z_2 \in \bgamh$ satisfy
  $d(z_1, z_2) > \delta$, then $\psi_{\rho'}(z_1)$ is opposite to
  $\psi_{\rho'}(z_2)$.

  Let $\{\alpha_n\}$, $\{\beta_n\}$ be $\mc{G}$-paths respectively
  limiting to points $z_1, z_2 \in \bgamh$ with
  $d(z_1, z_2) > \delta$, with corresponding vertex paths $\{v_n\}$
  and $\{w_n\}$. By \Cref{lem:initial_vertex_neighborhood}, we must
  have $\phi^{-1}(z_1) \subset U_{v_1}$ and
  $\phi^{-1}(z_2) \subset U_{w_2}$, so under
  \Cref{ass:opposite_open_sets}, we know that $U_{v_1}$ and $U_{w_1}$
  are opposite. But then we are done since
  \Cref{lem:initial_vertex_neighborhood} also implies that
  $\psi_{\rho'}(z_1) \subset U_{v_1}$ and
  $\psi_{\rho'}(z_2) \subset U_{w_1}$.
\end{proof}

\subsection{The boundary set of the deformed representation}

We define our candidate boundary set $\Lambda' \subset \flags$ by
\[
  \Lambda' = \bigcup_{z \in \bgamh} \psi_{\rho'}(z).
\]
We then have an equivariant map
\[
  \phi':\Lambda' \to \bgamh,
\]
where $\phi'(x) = z$ if $x \in
\psi_{\rho'}(z)$. \Cref{lem:boundary_map_well_defined} implies that
$\phi'$ is well-defined and antipodal. It is necessarily both
surjective and $\rho'$-equivariant, and its fibers are either
singletons or translates of the sets $\Lambda_p'$ for $p \in \Pi$.

It now remains to verify the properties of the candidate set
$\Lambda'$ and the map $\phi'$ needed to show that $\phi'$ is an EGF
boundary extension.
\begin{lem}
  \label{lem:deformed_automaton_nonempty_int}
  For every vertex $v$ of $\mc{G}$, the intersection $\Lambda' \cap
  U_v$ is nonempty.
\end{lem}
\begin{proof}
  The construction in \Cref{sec:extended_convergence} ensures that
  every vertex of the automaton $\mc{G}$ has at least one outgoing
  edge. In particular this means that for a given vertex $v$, there is
  an infinite $\mc{G}$-path whose first vertex is $v$. This
  $\mc{G}$-path limits to a conical limit point $z$, and
  \Cref{lem:initial_vertex_neighborhood} implies that
  $\psi_{\rho'}(z)$ is a (nonempty) subset of both $U_v$ and
  $\Lambda'$.
\end{proof}

\begin{lem}
  For any $z \in Z$, we have $\phi'^{-1}(z) \subset V$.
\end{lem}
\begin{proof}
  Recall that we used \Cref{prop:compact_bowditch_nbhd} to construct
  our automaton so that for any $z \in Z$, there is a $\mc{G}$-path
  limiting to $z$ with vertex path $\{v_n\}$ such that
  $U_{v_1} \subset V$. Then \Cref{lem:initial_vertex_neighborhood}
  implies $\phi'^{-1}(z) \subset V$.
\end{proof}

\begin{lem}
  \label{lem:deformed_map_proper}
  $\Lambda'$ is compact.
\end{lem}
\begin{proof}
  Fix a sequence $y_n \in \Lambda'$, and let $x_n = \phi'(y_n)$. Since
  $\bgamh$ is compact, up to subsequence $x_n \to x$. We want to
  see that a subsequence of $y_n$ converges to some $y \in
  \Lambda'$. We consider two possibilities:
  \begin{aside}{Case 1: $x$ is a parabolic point}
    We can write $x = g \cdot p$, where $p \in \Pi$. Let $w$ be a
    parabolic vertex with $p_w = p$, and consider the compact set
    $K_p \subset \bgamh - \{p\}$, chosen so that
    $\Gamma_p K_p = \bgamh - \{p\}$. If $x_n = x$ for infinitely many
    $n$, then $y_n \in \Lambda'_x \subseteq \Lambda'$, and we are done
    since $\Lambda_x'$ is compact by construction. So, assume
    otherwise, and choose $\gamma_n \in \Gamma_p$ so that
    $z_n = \gamma_n^{-1}g^{-1}x_n \in K_p$.

    We have assumed (using \Cref{prop:parabolic_outward_edge_nbhds})
    that the automaton $\mc{G}$ has been constructed so that there is
    always a $\mc{G}$-path limiting to $z_n$ whose first vertex $v_n$
    is connected to $w$ by an edge $(w,
    v_n)$. \Cref{lem:initial_vertex_neighborhood} implies that
    $\phi'^{-1}(z_n)$ lies in $U_n$, which is contained in $C_p'$ by
    definition.
    
    Then using \Cref{prop:convergence_dynamics_deformed_parabolics},
    we know that up to subsequence,
    \[
      \rho'(\gamma_n) \phi'^{-1}(z_n) = \rho'(\gamma_n)
      \phi'^{-1}(\gamma_n^{-1} g^{-1}x_n)
    \]
    converges to a compact subset of $\Lambda_p'$, which means that
    \[
      y_n \in \rho'(g)\rho'(\gamma_n) \phi'^{-1}(\gamma_n^{-1} g^{-1} x_n)
    \]
    subconverges to a point in $\rho'(g) \Lambda_p'$.
  \end{aside}
  \begin{aside}{Case 2: $x$ is a conical limit point}
    We want to show that any sequence in $\phi'^{-1}(x_n)$ limits to
    $\phi'^{-1}(x)$, so fix any $\eps > 0$. Using
    \Cref{cor:exponential_path_shrinking}, we can choose $N$ so that
    if $\{\alpha_m\}$ is any $\mc{G}$-path limiting to $x$, with
    corresponding vertex path $\{v_m\}$, then the diameter of
    \[
      \rho'(\alpha_1 \cdots \alpha_N) U_{v_{N+1}}
    \]
    is less than $\eps$ with respect to a metric on $U_{v_1}$. We fix
    such a $\mc{G}$-path $\{\alpha_m\}$. Then, we use
    \Cref{lem:close_points_share_prefixes} to see that for
    sufficiently large $n$, there is a $\mc{G}$-path $\{\beta^n_m\}$
    limiting to $x_n$ with $\beta_i = \alpha_i$ for $i \le N$. Thus
    the Hausdorff distance (with respect to $C_{U_{v_1}}$) between
    $\phi'^{-1}(x_n)$ and $\phi'^{-1}(x)$ is at most $\eps$. Since
    $\phi'^{-1}(x_n)$ and $\phi'^{-1}(x)$ both lie in the compact set
    $\rho'(\alpha_1)\overline{U_{v_2}} \subset U_{v_1}$, this proves
    the claim.
  \end{aside}
\end{proof}

\begin{lem}
  $\phi'$ is continuous and proper.
\end{lem}
\begin{proof}
  Since $\Lambda'$ is compact, we just need to show continuity. Fix
  $y \in \Lambda'$ and a sequence $y_n \in \Lambda'$ approaching
  $y$. We want to show that $\phi'(y_n)$ approaches $\phi'(y) = x$.

  Suppose otherwise. Since $\bgamh$ is compact, up to a
  subsequence $z_n = \phi'(y_n)$ approaches $z \ne x$. Using the
  equivariance of $\phi'$, and cocompactness of the $\Gamma$-action on
  distinct pairs in $\bgamh$, we may assume that
  $d(x, z) > \delta$. For sufficiently large $n$, we have
  $d(x, z_n) > \delta$ as well. Then, as in the proof of
  \Cref{lem:boundary_map_well_defined}, by
  \Cref{ass:opposite_open_sets} we know that for any vertices $v,w$ in
  $\mc{G}$ such that $U_v$ contains $\psi_{\rho'}(x)$ and $U_w$
  contains $\psi_{\rho'}(z_n)$, the intersection
  $\overline{U_v} \cap \overline{U_w}$ is empty.

  But by definition of $\phi'$, we have
  \[
    y \in \psi_{\rho'}(x) \subset U_v, \qquad y_n \in
    \psi_{\rho'}(z_n) \subset U_w
  \]
  for vertices $v, w$ in $\mc{G}$. This contradicts the fact that
  $y_n \to y$.
\end{proof}

\subsection{Extending the convergence action}

We have constructed a surjective $\rho'$-equivariant antipodal map
$\phi':\Lambda' \to \bgamh$. To prove that $\rho'$ is EGF, we will
apply \Cref{prop:weak_egf_to_egf_ii}, so we need to define repelling
strata $C_z'$ and check that they satisfy condition
\ref{item:weak_strata_interior}. For each conical $z \in \bgamh$, let
$C_z' = \Opp((\phi')^{-1}(z))$. We have already defined the strata
$C_p'$ for parabolic $p \in \Pi$, so for arbitrary parabolic $q$, let
$g \in \Gamma$ so that $q = gp$, and define $C_q' = \rho'(g)C_p'$. The
choice of $g$ does not matter since we have defined the sets $C_p'$ to
be $\Gamma_p$-invariant.

\begin{lem}
  The sets $C_z'$ satisfy condition \ref{item:weak_strata_interior};
  that is, for every $z \in \bgamh$, we have
  $\Lambda - \phi^{-1}(z) \subset C_z'$.
\end{lem}
\begin{proof}
  For conical $z \in \bgamh$ this is immediate from antipodality. For
  parabolic points, it suffices to show that the condition holds for
  $p \in \Pi$. Consider the compact set $K_p \subset \bgamh - \{p\}$
  satisfying $\Gamma_pK_p = \bgamh - \{p\}$. We observed in the proof
  of \Cref{lem:deformed_map_proper} that $C_p'$ contains
  $\phi'^{-1}(K_p)$. But then since $C_p'$ is
  $\rho'(\Gamma_p)$-invariant (by
  \Cref{prop:convergence_dynamics_deformed_parabolics}), we have
  \[
    \rho'(\Gamma_p) \cdot \phi'^{-1}(K_p) = \phi'^{-1}(\bgamh - \{p\})
    \subset C_p'.
  \]
\end{proof}

The final step in the proof of \Cref{thm:cusp_stable_stability} is to
show that conditions \ref{item:parabolic_flag_convergence} and
\ref{item:rel_p_div} in \Cref{prop:weak_egf_to_egf_ii} also
hold. Condition \ref{item:parabolic_flag_convergence} is immediate
from \Cref{prop:convergence_dynamics_deformed_parabolics}, and the
equivariance in the definition of the strata $C_q'$ for $q$
parabolic. Finally, to see that condition \ref{item:rel_p_div} holds,
observe that the relative quasi-geodesic automaton $\mc{G}$, the map
$\phi'$, and the $\mc{G}$-compatible system $\{U_v\}$ satisfy the
hypotheses of \Cref{lem:automaton_divergent_limits}. The conclusion of
this lemma gives us the desired condition. This concludes the proof of
\Cref{thm:cusp_stable_stability}.

\begin{remark}
  \label{rem:conical_pts_singletons}
  The definition of the set $\Lambda'$ and the map $\phi'$ immediately
  imply that the \emph{fibers} of the deformed boundary extension
  $\phi':\Lambda' \to \bgamh$ satisfy the conclusions of
  \Cref{prop:conical_pts_singletons}: the fiber over each conical
  limit point is a singleton, and the fiber over each parabolic point
  $p$ is the closure of the accumulation sets of $\Gamma_p$-orbits in
  $C_p'$. So, we obtain \Cref{prop:conical_pts_singletons} by taking
  $\mc{W}$ to be the singleton $\{\rho\}$, and following the proof of
  \Cref{thm:cusp_stable_stability} (using $C_p$ for $C_p'$
  throughout).
\end{remark}










\appendix

\section{Contraction dynamics on flag manifolds}
\label{sec:contraction_appendix}

Let $V$ be a real vector space, and let $A_n$ be a sequence of
elements of $\PGL(V)$. It is sometimes possible to determine the
global dynamical behavior of $A_n$ on $\P(V)$ by considering the
action of $A_n$ on a small open subset of $\P(V)$: if there is an open
subset $U \subset \P(V)$ such that $A_n \cdot U$ converges to a point
in $\P(V)$, then in fact there is a dense open subset
$U_- \subset \P(V)$ (the complement of a hyperplane) on which $A_n$
converges to the same point, uniformly on compacts.

A similar statement holds for the action of $A_n$ on Grassmannians
$\Gr(k, V)$. These claims can be proved by considering the behavior
of the singular value gaps of $A_n$ as $n \to \infty$.

In this appendix we give a general result along these lines, where we
take sequences of group elements $g_n \in G$ for a semisimple Lie
group $G$ with no compact factor and trivial center, and consider the
limiting behavior of $g_n$ on open subsets of some flag manifold
$\pflags$, where $\pargroup$ is a parabolic subgroup.

\localContractingImpliesGlobalContracting*

We will prove \Cref{prop:local_contracting_implies_global_contracting}
by reducing it to the case where $G = \PGL(d,\R)$ and
$\pargroup = P_1$ is the stabilizer of
$[e_1] \in \RP^{d-1} \simeq G/P_1$. In this situation,
$\pargroup$-divergence can be understood in terms of the behavior of
the \emph{singular value gaps} of the sequence $g_n$:

\begin{prop}
  Suppose that $G = \PGL(d, \R)$, and let $\pargroup = P_1 \subset G$
  be the stabilizer of a line in $\R^d$. A sequence $g_n \in G$ is
  $P_1$-divergent if and only if
  \[
    \frac{\sigma_1(g_n)}{\sigma_2(g_n)} \to \infty,
  \]
  where $\sigma_i(g_n)$ is the $i$th-largest singular value of $g_n$.
\end{prop}

For convenience, we give a proof of
\Cref{prop:local_contracting_implies_global_contracting} in this
special case.
\begin{lem}
  \label{lem:local_contracting_implies_sv_gaps}
  Let $g_n$ be a sequence in $\PGL(d,\R)$, and suppose that for a
  nonempty open subset $U \subset \RP^{d-1}$, $g_nU$ converges to a
  point in $\RP^{d-1}$. Then, the singular value gap
  \[
    \frac{\sigma_1(g_n)}{\sigma_2(g_n)}
  \]
  tends to $\infty$ as $n \to \infty$.
\end{lem}
\begin{proof}
  It suffices to show that any subsequence of $g_n$ has a subsequence
  which satisfies the property. Using the Cartan decomposition of
  $\PGL(d,\R)$, we can write
  \[
    g_n = k_na_nk_n',
  \]
  for $k_n, k_n' \in K = \PO(d)$ and $a_n$ a diagonal matrix whose
  diagonal entries are $\sigma_1, \ldots, \sigma_d$. Up to subsequence
  $k_n$ and $k_n'$ converge respectively to $k, k' \in K$. For
  sufficiently large $n$, $k_n' U \cap k' U$ is nonempty, so by
  replacing $U$ with $k'U$ we can assume that $k_n' = \mr{id}$ for all
  $n$. Furthermore, if $k_na_n U$ converges to a point
  $z \in \RP^{d-1}$, then $a_nU$ converges to $k^{-1}z$.

  So, $a_nU$ converges to a point, and since $a_n$ is a diagonal
  matrix, the gap between the moduli of its largest and second-largest
  eigenvalues must be unbounded.
\end{proof}

To prove the general case of
\Cref{prop:local_contracting_implies_global_contracting}, we take an
irreducible representation $\zeta:G \to \PGL(V)$ coming from
\Cref{thm:flag_irrep_embedding}, so that $\pargroup$ maps to the
stabilizer of a line $\ell$ in $V$, $\oppgroup$ maps to the stabilizer
of a hyperplane $H$ in $V$, and $g\pargroup g^{-1}$,
$h\oppgroup h^{-1}$ are opposite if and only if $\zeta(g)\ell$,
$\zeta(h)H$ are transverse. As in \cref{sec:flag_convergence}, this
determines embeddings $\iota:\flags \to \P(V)$ and
$\iota^*:\oppflags \to \P(V^*)$ by
\[
  \iota(g\pargroup) = \zeta(g)\ell, \qquad \iota^*(g\oppgroup) =
  \zeta(g)H.
\]

The representation $\zeta$ additionally has the property that for any
sequence $g_n \in G$, the singular value gaps
\[
  \sigma_1(\zeta(g_n)) / \sigma_2(\zeta(g_n))
\]
are unbounded if and only if $g_n$ is $\pargroup$-divergent (see
\cite{ggkw2017anosov}, Lemma 3.7).

\begin{proof}[Proof of
  \Cref{prop:local_contracting_implies_global_contracting}]
  By \cite{zimmer2018proper}, Lemma 4.7, there exist flags
  $\xi_1, \ldots, \xi_D \in U$ so that lifts of $\iota(\xi_i)$ give a
  basis of $V$. Since $g_n \cdot U$ converges to a point in $\flags$,
  the set
  \[
    \{\zeta(g_n) \cdot \iota(\xi_i) : 1 \le i \le D\}
  \]
  converges to a single point in $\P(V)$.

  This means that we can fix lifts $\tilde{\iota(\xi_i)} \in V$ so
  that, up to a subsequence, $\zeta(g_n)$ takes the projective
  $(D-1)$-simplex
  \[
    \left[\sum_{i=1}^D \lambda_i \tilde{\iota(\xi_i)} : \lambda_i >
      0\right] \subset \P(V)
  \]
  to a point. This simplex is an open subset of $\P(V)$. Now we can
  apply \Cref{lem:local_contracting_implies_sv_gaps} to see that the
  sequence $g_n$ is $\pargroup$-divergent.

  We now just need to check that $\xi$ is the unique $P^+$-limit point
  of $g_n$. Choose any subsequence of $g_n$. Then any
  $P^+$-contracting subsequence $g_m$ of this subsequence satisfies
  \[
    g_m|_{\Opp(\xi_-)} \to \xi'
  \]
  uniformly on compacts for some $\xi_- \in \oppflags$ and
  $\xi' \in \pflags$. But since $\Opp(\xi_-)$ is open and dense, it
  intersects $U$ nontrivially and thus $\xi' = \xi$.
\end{proof}


\bibliographystyle{alpha}
\bibliography{references}

\newcommand{\etalchar}[1]{$^{#1}$}
\begin{thebibliography}{{Wan}23b}

\bibitem[Bal14]{ballas2014}
Samuel~A. Ballas.
\newblock Deformations of noncompact projective manifolds.
\newblock {\em Algebr. Geom. Topol.}, 14(5):2595--2625, 2014.

\bibitem[Bal21]{ballas2021}
Samuel~A. Ballas.
\newblock Constructing convex projective 3-manifolds with generalized cusps.
\newblock {\em J. Lond. Math. Soc. (2)}, 103(4):1276--1313, 2021.

\bibitem[BCL20]{bcl2020generalized}
Samuel~A. Ballas, Daryl Cooper, and Arielle Leitner.
\newblock Generalized cusps in real projective manifolds: classification.
\newblock {\em J. Topol.}, 13(4):1455--1496, 2020.

\bibitem[BCM12]{BCM12}
Pierre-Louis {Blayac}, Micka{\"e}l {Crampon}, and Ludovic {Marquis}.
\newblock {Finitude g{\'e}om{\'e}trique en g{\'e}om{\'e}trie de Hilbert + an
  erratum/addendum}.
\newblock {\em arXiv e-prints}, page arXiv:1202.5442, February 2012.

\bibitem[BDL15]{bdl2015convex}
Samuel Ballas, Jeffrey Danciger, and Gye-Seon Lee.
\newblock Convex projective structures on non-hyperbolic three-manifolds.
\newblock {\em Geometry \& Topology}, 22, 08 2015.

\bibitem[Ben04]{benoist04convexes}
Yves Benoist.
\newblock Convexes divisibles. {I}.
\newblock In {\em Algebraic groups and arithmetic}, pages 339--374. Tata Inst.
  Fund. Res., Mumbai, 2004.

\bibitem[Ben06a]{benoist2006convexes}
Yves Benoist.
\newblock Convexes divisibles. {IV}. {S}tructure du bord en dimension 3.
\newblock {\em Invent. Math.}, 164(2):249--278, 2006.

\bibitem[Ben06b]{benoist2006convexeshyp}
Yves Benoist.
\newblock Convexes hyperboliques et quasiisom\'{e}tries.
\newblock {\em Geom. Dedicata}, 122:109--134, 2006.

\bibitem[BH99]{bh1999metric}
Martin~R. Bridson and Andr\'{e} Haefliger.
\newblock {\em Metric spaces of non-positive curvature}, volume 319 of {\em
  Grundlehren der Mathematischen Wissenschaften [Fundamental Principles of
  Mathematical Sciences]}.
\newblock Springer-Verlag, Berlin, 1999.

\bibitem[BM20]{bm2020properly}
Samuel Ballas and Ludovic Marquis.
\newblock Properly convex bending of hyperbolic manifolds.
\newblock {\em Groups Geom. Dyn.}, 14(2):653--688, 2020.

\bibitem[Bob19]{bobb2019convex}
Martin~D. Bobb.
\newblock Convex projective manifolds with a cusp of any non-diagonalizable
  type.
\newblock {\em J. Lond. Math. Soc. (2)}, 100(1):183--202, 2019.

\bibitem[Bog97]{bogopolski97infinite}
O.~V. Bogopol'ski\u{\i}.
\newblock Infinite commensurable hyperbolic groups are bi-{L}ipschitz
  equivalent.
\newblock {\em Algebra i Logika}, 36(3):259--272, 357, 1997.

\bibitem[Bow98]{bowditch1998spaces}
Brian~H. Bowditch.
\newblock Spaces of geometrically finite representations.
\newblock {\em Ann. Acad. Sci. Fenn. Math.}, 23(2):389--414, 1998.

\bibitem[Bow99]{bowditch1999convergence}
B.~H. Bowditch.
\newblock Convergence groups and configuration spaces.
\newblock In {\em Geometric group theory down under ({C}anberra, 1996)}, pages
  23--54. de Gruyter, Berlin, 1999.

\bibitem[Bow12]{bowditch2012relatively}
B.~H. Bowditch.
\newblock Relatively hyperbolic groups.
\newblock {\em Internat. J. Algebra Comput.}, 22(3):1250016, 66, 2012.

\bibitem[BPS19]{bps2019anosov}
Jairo Bochi, Rafael Potrie, and Andr\'{e}s Sambarino.
\newblock Anosov representations and dominated splittings.
\newblock {\em J. Eur. Math. Soc. (JEMS)}, 21(11):3343--3414, 2019.

\bibitem[BV23]{BV}
Pierre-Louis {Blayac} and Gabriele {Viaggi}.
\newblock {Divisible convex sets with properly embedded cones}.
\newblock {\em arXiv e-prints}, page arXiv:2302.07177, February 2023.

\bibitem[{Cho}10]{choi2011}
Suhyoung {Choi}.
\newblock {The convex real projective orbifolds with radial or totally geodesic
  ends: The closedness and openness of deformations}.
\newblock {\em arXiv e-prints}, page arXiv:1011.1060, November 2010.

\bibitem[CLM20]{clm2020convex}
Suhyoung Choi, Gye-Seon Lee, and Ludovic Marquis.
\newblock Convex projective generalized {D}ehn filling.
\newblock {\em Ann. Sci. \'{E}c. Norm. Sup\'{e}r. (4)}, 53(1):217--266, 2020.

\bibitem[CLM22]{CLM2}
Suhyoung Choi, Gye-Seon Lee, and Ludovic Marquis.
\newblock Deformation spaces of {C}oxeter truncation polytopes.
\newblock {\em J. Lond. Math. Soc. (2)}, 106(4):3822--3864, 2022.

\bibitem[CLT18]{clt2018deforming}
Daryl Cooper, Darren Long, and Stephan Tillmann.
\newblock Deforming convex projective manifolds.
\newblock {\em Geom. Topol.}, 22(3):1349--1404, 2018.

\bibitem[CM14]{cm2014finitude}
Micka\"{e}l Crampon and Ludovic Marquis.
\newblock Finitude g\'{e}om\'{e}trique en g\'{e}om\'{e}trie de {H}ilbert.
\newblock {\em Ann. Inst. Fourier (Grenoble)}, 64(6):2299--2377, 2014.

\bibitem[CZZ21]{czz2021cusped}
Richard {Canary}, Tengren {Zhang}, and Andrew {Zimmer}.
\newblock {Cusped Hitchin representations and Anosov representations of
  geometrically finite Fuchsian groups}.
\newblock {\em arXiv e-prints}, page arXiv:2103.06588, March 2021.

\bibitem[DGK18]{DGKgeomded}
Jeffrey Danciger, Fran\c{c}ois Gu\'{e}ritaud, and Fanny Kassel.
\newblock Convex cocompactness in pseudo-{R}iemannian hyperbolic spaces.
\newblock {\em Geom. Dedicata}, 192:87--126, 2018.

\bibitem[DGK{\etalchar{+}}21]{dgklm2021convex}
Jeffrey {Danciger}, Fran{\c{c}}ois {Gu{\'e}ritaud}, Fanny {Kassel}, Gye-Seon
  {Lee}, and Ludovic {Marquis}.
\newblock {Convex cocompactness for Coxeter groups}.
\newblock {\em arXiv e-prints}, page arXiv:2102.02757, February 2021.

\bibitem[DGK24]{dgk2017convex}
Jeffrey Danciger, Fran\c~cois Gu\'eritaud, and Fanny Kassel.
\newblock Convex cocompact actions in real projective geometry.
\newblock {\em Ann. Sci. \'Ec. Norm. Sup\'er. (4)}, 57(6):1753--1843, 2024.

\bibitem[DS05]{ds2005tree}
Cornelia Dru\c{t}u and Mark Sapir.
\newblock Tree-graded spaces and asymptotic cones of groups.
\newblock {\em Topology}, 44(5):959--1058, 2005.
\newblock With an appendix by Denis Osin and Mark Sapir.

\bibitem[Ebe96]{eberlein}
Patrick~B. Eberlein.
\newblock {\em Geometry of nonpositively curved manifolds}.
\newblock Chicago Lectures in Mathematics. University of Chicago Press,
  Chicago, IL, 1996.

\bibitem[Ger12]{gerasimov2012floyd}
Victor Gerasimov.
\newblock Floyd maps for relatively hyperbolic groups.
\newblock {\em Geom. Funct. Anal.}, 22(5):1361--1399, 2012.

\bibitem[GGKW17]{ggkw2017anosov}
Fran\c{c}ois Gu\'{e}ritaud, Olivier Guichard, Fanny Kassel, and Anna Wienhard.
\newblock Anosov representations and proper actions.
\newblock {\em Geom. Topol.}, 21(1):485--584, 2017.

\bibitem[GM08]{gm2008dehn}
Daniel Groves and Jason~Fox Manning.
\newblock Dehn filling in relatively hyperbolic groups.
\newblock {\em Israel J. Math.}, 168:317--429, 2008.

\bibitem[GP13]{gp2013quasi}
Victor Gerasimov and Leonid Potyagailo.
\newblock Quasi-isometric maps and {F}loyd boundaries of relatively hyperbolic
  groups.
\newblock {\em J. Eur. Math. Soc. (JEMS)}, 15(6):2115--2137, 2013.

\bibitem[GW12]{gw2012anosov}
Olivier Guichard and Anna Wienhard.
\newblock Anosov representations: domains of discontinuity and applications.
\newblock {\em Invent. Math.}, 190(2):357--438, 2012.

\bibitem[GW24]{GW2024}
Antonin {Guilloux} and Theodore {Weisman}.
\newblock Limits of limit sets in rank-one symmetric spaces.
\newblock 2024.

\bibitem[Hel01]{helgason}
Sigurdur Helgason.
\newblock {\em Differential geometry, {L}ie groups, and symmetric spaces},
  volume~34 of {\em Graduate Studies in Mathematics}.
\newblock American Mathematical Society, Providence, RI, 2001.
\newblock Corrected reprint of the 1978 original.

\bibitem[IZ22]{iz2022structure}
Mitul {Islam} and Andrew {Zimmer}.
\newblock {The structure of relatively hyperbolic groups in convex real
  projective geometry}.
\newblock {\em arXiv e-prints}, page arXiv:2203.16596, March 2022.

\bibitem[Kap07]{kapovich2007}
Michael Kapovich.
\newblock Convex projective structures on {G}romov-{T}hurston manifolds.
\newblock {\em Geom. Topol.}, 11:1777--1830, 2007.

\bibitem[KKL19]{kkl2019structural}
Michael {Kapovich}, Sungwoon {Kim}, and Jaejeong {Lee}.
\newblock {Structural stability of meandering-hyperbolic group actions}.
\newblock {\em arXiv e-prints}, page arXiv:1904.06921, April 2019.

\bibitem[KL23]{kl2018relativizing}
Michael Kapovich and Bernhard Leeb.
\newblock Relativizing characterizations of {A}nosov subgroups, {I}.
\newblock {\em Groups Geom. Dyn.}, 17(3):1005--1071, 2023.
\newblock With an appendix by Gregory A. Soifer.

\bibitem[KLP17]{klp2017anosov}
Michael Kapovich, Bernhard Leeb, and Joan Porti.
\newblock Anosov subgroups: dynamical and geometric characterizations.
\newblock {\em Eur. J. Math.}, 3(4):808--898, 2017.

\bibitem[Kna02]{knapp}
Anthony~W. Knapp.
\newblock {\em Lie groups beyond an introduction}, volume 140 of {\em Progress
  in Mathematics}.
\newblock Birkh\"{a}user Boston, Inc., Boston, MA, second edition, 2002.

\bibitem[Lab06]{labourie2006anosov}
Fran\c{c}ois Labourie.
\newblock Anosov flows, surface groups and curves in projective space.
\newblock {\em Invent. Math.}, 165(1):51--114, 2006.

\bibitem[LLS26]{LLS}
Gye-Seon Lee, Jaejeong Lee, and Florian Stecker.
\newblock Anosov triangle reflection groups in {$\rm{SL}(3,\Bbb R)$}.
\newblock {\em J. Differential Geom.}, 132(1):203--253, 2026.

\bibitem[MMW22]{MMW1}
Kathryn Mann, Jason~Fox Manning, and Theodore Weisman.
\newblock Stability of hyperbolic groups acting on their boundaries.
\newblock Preprint, \href{http://arxiv.org/abs/2206.14914}{arXiv:2206.14914},
  2022.

\bibitem[MMW24]{MMW2}
Kathryn {Mann}, Jason~Fox {Manning}, and Theodore {Weisman}.
\newblock {Topological stability of relatively hyperbolic groups acting on
  their boundaries}.
\newblock {\em arXiv e-prints}, page arXiv:2402.06144, February 2024.

\bibitem[{Rie}21]{riestenberg2021quantified}
J.~Maxwell {Riestenberg}.
\newblock {A quantified local-to-global principle for Morse quasigeodesics}.
\newblock {\em arXiv e-prints}, page arXiv:2101.07162, January 2021.

\bibitem[Tra13]{tran2013relations}
Hung~Cong Tran.
\newblock Relations between various boundaries of relatively hyperbolic groups.
\newblock {\em Internat. J. Algebra Comput.}, 23(7):1551--1572, 2013.

\bibitem[Tuk94]{tukia1994convergence}
Pekka Tukia.
\newblock Convergence groups and {G}romov's metric hyperbolic spaces.
\newblock {\em New Zealand J. Math.}, 23(2):157--187, 1994.

\bibitem[Tuk98]{tukia1998conical}
Pekka Tukia.
\newblock Conical limit points and uniform convergence groups.
\newblock {\em J. Reine Angew. Math.}, 501:71--98, 1998.

\bibitem[Wan23a]{wang2023a}
Tianqi Wang.
\newblock Anosov representations over closed subflows.
\newblock {\em Trans. Amer. Math. Soc.}, 376(9):6177--6214, 2023.

\bibitem[{Wan}23b]{wang2023b}
Tianqi {Wang}.
\newblock {Notions of Anosov representation of relatively hyperbolic groups}.
\newblock {\em arXiv e-prints}, page arXiv:2309.15636, September 2023.

\bibitem[Wei23a]{weisman2023dynamical}
Theodore Weisman.
\newblock Dynamical properties of convex cocompact actions in projective space.
\newblock {\em Journal of Topology}, 16(3):990--1047, 2023.

\bibitem[Wei23b]{WeismanExamples}
Theodore Weisman.
\newblock Examples of extended geometrically finite representations.
\newblock 2023.

\bibitem[Wol20]{wolf2020}
Adva Wolf.
\newblock {\em Convex {P}rojective {G}eometrically {F}inite {S}tructures}.
\newblock ProQuest LLC, Ann Arbor, MI, 2020.
\newblock Thesis (Ph.D.)--Stanford University.

\bibitem[Yam04]{yaman2004topological}
Asli Yaman.
\newblock A topological characterisation of relatively hyperbolic groups.
\newblock {\em J. Reine Angew. Math.}, 566:41--89, 2004.

\bibitem[Zhu21]{zhu2019relatively}
Feng Zhu.
\newblock Relatively dominated representations.
\newblock {\em Ann. Inst. Fourier (Grenoble)}, 71(5):2169--2235, 2021.

\bibitem[Zim18]{zimmer2018proper}
Andrew~M. Zimmer.
\newblock Proper quasi-homogeneous domains in flag manifolds and geometric
  structures.
\newblock {\em Ann. Inst. Fourier (Grenoble)}, 68(6):2635--2662, 2018.

\bibitem[Zim21]{zimmer2021projective}
Andrew Zimmer.
\newblock Projective {A}nosov representations, convex cocompact actions, and
  rigidity.
\newblock {\em J. Differential Geom.}, 119(3):513--586, 2021.

\bibitem[ZZ22]{zz1}
Feng {Zhu} and Andrew {Zimmer}.
\newblock {Relatively Anosov representations via flows I: theory}.
\newblock {\em arXiv e-prints}, page arXiv:2207.14737, July 2022.

\end{thebibliography}

\end{document}